\documentclass[review]{elsarticle}

\usepackage{lineno,hyperref}
\usepackage{amsmath}
\usepackage{amssymb}
\newtheorem{remark}{Remark}

\usepackage{tikz}
\modulolinenumbers[5]

\begin{document}

\begin{frontmatter}

\title{A mesh-free particle  method  for  continuum modeling of granular flow}


\author[TUKL1]{Sudarshan Tiwari\corref{cor1}}
\ead{tiwari@mathematik.uni-kl.de}
\cortext[cor1]{Corresponding Author}
\author[TUKL]{Axel Klar}
\ead{klar@mathematik.uni-kl.de}

\address[TUKL1]{Fachbereich Mathematik, TU Kaiserslautern,
         Gottlieb-Daimler-Strasse, 
         67663 Kaiserslautern, Germany}
\address[TUKL]{Fachbereich Mathematik, TU Kaiserslautern,
         Gottlieb-Daimler-Strasse, 
         67663 Kaiserslautern, Germany\\
         Fraunhofer ITWM Kaiserslautern, 67663 Kaiserslautern, Germany} 
%

\begin{abstract}
Based on the continuum model for granular media developed in   Dunatunga et al. we propose  a mesh-free generalized finite difference method
 	for the simulation of granular flows. The model is given by an elasto-viscoplastic model with a yield criterion using the $\mu(I)$ rheology from Jop et al. The numerical procedure is based on a mesh-free particle method with a least squares approximation of the derivatives in the balance equations combined with the numerical algorithm developed in Dunatunga et al. to compute the plastic stresses.
 	The method is numerically tested and verified for  several numerical experiments including granular column collapse and rigid body motion in granular materials. For comparison a nonlinear microscopic model from Lacaze et al.  is implemented and results are compared to 
 	the those obtained from the continuum model for granular column collapse and rigid body coupling to granular flow. 

\end{abstract}

\begin{keyword}
\texttt{Mesh-free particle method,   granular flow, elasto-viscoplastic model}
\MSC[2010] 76T25, 76M28
\end{keyword}

\end{frontmatter}

\linenumbers

\section{Introduction}
 \label{introduction}

Granular materials have been modelled using  microscopic and continuum  approaches in many works.
The most direct approach is the microscopic one modelling each  sand particle individually and describing its motion by 
Newton-type equations.
In the context of granular flow the most popular of these microscopic approaches has been termed DEM
\cite{CS79,CB02,LPK08,LLH,SDW},

On the other hand several different  approaches have been developed on the macroscopic level.  A rather recent modelling approach is based on the 
 $\mu(I)$-rheology, see \cite{C05, JFP06}.  This offers an efficient ansatz  for a rheological model for sand with advantages compared to classical Drucker-Prager models  and has been extended in a variety of recent works, see, for example, 
\cite{DK15,LSP,JLP, KAMRIN2010167, S87, BSSG}.

 There exist a great variety of computational methods for the above discussed continuum equations.
 We concentrate on meshfree methods and their use for granular flow simulation.
 A very popular method is the  material point method, see  \cite{DBD,MAMM15}. For its  application in the present context  we refer to  \cite{DK15}.
 However, for continuum mechanics many other meshfree approaches exist.
 For example those based on generalized finite differences  (GFDM)\cite{OK} avoid the additional background mesh in the MPM method
 and may be more efficient in some contexts.

 In the present paper we use the same set of equations as in \cite{DK15} but implement them in the framework of a GFDM method, see \cite{SKT}, showing that the approach developed in \cite{DK15} can be flexibly used in different meshfree CFD implementations. Moreover, for comparison, a microscopic nonlinear Hertz-Mindlin model based on the considerations in \cite{LPK08}
 is implemented in the same numerical framework.
 We note that the parameters  of the  microscopic model are physically consistent with the macroscopic ones.

 The resulting implementation is numerically verified and validated for several test cases.  A classical well-investigated test  case used in many papers is the collapse of a granular column, see  \cite{LH,LSP}. We will use this test-case as our main benchmark to compare microscopic and macroscopic implementations.
  Moreover, as a second test case, we consider  a moving rigid body in sand. In particular, we investigate the case of  a 
sphere falling into  a box of sand. We refer to 
 \cite{NGM10} for other approaches investigating such a coupling of granular flow and rigid body motion.
 
 
\section{Equations}
 \label{macroeq}
 
 We consider   the continuum  equations in  Lagrangian form
\begin{eqnarray}
\label{lagrange}	\dot{\bf x} &=& {\bf v}\\
\dot \rho  &=& - \rho \nabla\cdot {\bf v} \nonumber\\
\dot {\bf v} &=& \frac{1}{\rho} \nabla \cdot \sigma +  {\bf g} \nonumber
\end{eqnarray}
where $({\bf x},t) \in \mathbb{R}^3 \times \mathbb{R}^+$ is space-time and ${\bf v} = {\bf v}({\bf x},t)$ the velocity vector. Moreover, $\rho= \rho({\bf x},t)$ is the mass density, ${\bf g}$ the  gravitational force and $\sigma$ the Cauchy stress tensor. 

We define the velocity gradient 
\begin{equation}
\label{L}
L= \nabla {\bf v} = \left(\frac{\partial v^{(i)}}{\partial x^{(i)}} \right)_{ij}, 
\end{equation}
where ${\bf v} = (v^{(1)}, v^{(2)}, v^{(3)})$, ${\bf x} =  (x^{(1)}, x^{(2)}, x^{(3)})$. The  the strain-rate tensor $D$,  as well as the spin tensor $W$ are defind by ,
\begin{equation}
\label{DW}
D= \frac{1}{2} (L+ L^T), \; W =  \frac{1}{2} (L- L^T),
\end{equation}
compare \cite{DK15}.

To determine the stress tensor $\sigma= \sigma( \rho,L)$ we follow  the approach in \cite{JFP06} and \cite{DK15}
and use  the $\mu(I)$-rheology from \cite{JFP06}
and an hypoelastic-plastic approach  as in \cite{DK15} to close  equations (\ref{lagrange}).

\subsection{The stress-strain relation}

To obtain the relation $\sigma= \sigma(\rho, L)$ we 
define first (in three spatial dimensions)
\begin{align}
p = -\frac{1}{3} \mbox{tr} \sigma
\end{align}
and the strain-deviator 
\begin{align}
\sigma_0=  \sigma + p \mathbb I .
\end{align}
Moreover, define the shear stress and the elastic part of the strain rate tensor as 
\begin{align}
\bar \tau  = \frac{\vert \sigma_0 \vert}{\sqrt{2}}
\; \mbox{and } \;
	D_e= D-D_p,
\end{align}
where the plastic part of the strain rate tensor is given by 
\begin{align}
\label{dp}
D_p = \frac{1}{\sqrt{2}}  \dot {\bar{ \gamma} } (\bar \tau)
\frac{\sigma_0}{\vert \sigma_0 \vert},
\end{align}
where 
$$
\vert \sigma_0 \vert  = \sqrt{tr(\sigma_0^2)}.
$$
The plastic shear strain rate  $\dot {\bar{ \gamma} } (\bar \tau)$ is defined in the next subsection.

\subsection{The constitutive relation ($\mu(I)$-rheology)}

To define a constitutive law means 
to define   a   relation  $\dot{\bar \gamma} = \dot{\bar \gamma} (\bar \tau)$ leading to a relation
$D_p =  D_p (\sigma ) $.
It  is given by
the  considerations from  \cite{JFP06}. Let $\rho_s$ and $d$ be the grain density and mean particle size.
Define  the inertial number or normalized flow rate $I$ as
\[
I =   \frac{\dot{\bar \gamma}  d \sqrt{\rho_s}}{\sqrt{p}}  =  \frac{\dot{\bar \gamma} }{\xi \sqrt{p}} I_0
\]
with constants $\rho_s,d,I_0$ and $\xi = \frac{I_0}{d \sqrt{ \rho_s}}$
or
\[
\dot{\bar \gamma}  =  \xi \sqrt{p}  \frac{I }{I_0}.
\]
Then, the $\mu(I)$-rheology is a  relation between the friction coefficient $\mu$  and $I$ given by 
\begin{align}
\mu = \mu(I)  = \begin{cases} \mu_s + \frac{\mu_2 - \mu_s}{I_0/I + 1},& I>0 \\
\mu  \le \mu_s , &I=0.
\end{cases}
\end{align}

Since 
\[
\mu   = \frac{\bar \tau}{p} = \frac{\vert \sigma_0 \vert }{p \sqrt{2} },
\] 
we have 
\begin{align}
\bar \tau (\dot{\bar \gamma})= \begin{cases} p \left( \mu_s + \frac{\mu_2 - \mu_s}{\frac{\xi \sqrt{p}}{\dot {\bar \gamma}} + 1}\right) ,& \dot{\bar \gamma} >0 \\
\bar \tau \le  p \mu_s , & \dot{\bar \gamma} =0
\end{cases}
\end{align}
and the (well-defined) inverse
\begin{align}
\label{tau}
\dot{\bar \gamma}	(\bar \tau,p)= \begin{cases}   \sqrt{p} \xi    \frac{\bar \tau-\mu_s p }{\mu_2 p - \bar \tau} ,& \bar \tau > \mu_s  p\\
0  , &\bar \tau  \le  \mu_s p
\end{cases}.
\end{align}

This gives $D_p=D_p(\sigma)$.

\subsection{Determination of $\sigma$}
\label{Kamrin}
$\sigma$ is finally obtained by the   following definitions, see \cite{DK15}.
Let $\rho_C$ be a critical density.
For $\rho< \rho_C$ define $\sigma=0$.
For $\rho >\rho_C$ use 
\begin{align}
	\sigma_{\triangle} =  \dot \sigma - W \cdot \sigma + \sigma \cdot W
\end{align}
and
\begin{align}
	\sigma_{\triangle} =  \frac{E}{1+\nu} \left[   (D - D_p) + \frac{\nu}{1- 2 \nu} \mbox{tr}  (D-D_p) \mathbb I  \right]
\end{align}
with Youngs modulus $E = 3 K_C (1-2 \nu)$, compressive bulk modulus $K_C$ and Poisson ratio $\nu$.
This leads to  a closed ODE for $\sigma$ depending on $L$
\begin{align}
\label{ode}
	\dot \sigma  = \frac{E}{1+\nu} \left[   (D - D_p) + \frac{\nu}{1- 2 \nu} \mbox{tr}  (D-D_p) \mathbb I  \right]+ W \cdot \sigma -\sigma \cdot W
\end{align}
and finally to a closed set of equations.


 \section{Coloumb constitutive model with $\mu(I)$-rheology}
\label{Jop}

A simplification of the above equations  is obtained by considering the incompressible Navier-Stokes  equations in  Lagrangian form
\begin{eqnarray}
\dot{\bf x} &=& {\bf v} \\
\nabla\cdot {\bf v} &=& 0 \\
\dot{ \bf v} &=& \frac{1}{\rho_0} \nabla \sigma +  {\bf g} ,
\end{eqnarray}
where $\rho_0$ is the constant mass density, other quantities are the same as above and $\sigma$  is the Cauchy stress tensor, compare \cite{JFP06}.
In this simplified  case,  the granular material is described as an incompressible fluid with the internal stress tensor given by 

\begin{equation}
\sigma = \sigma_0 - p\mathbb I, \quad \mbox{where} \quad \sigma_0  = \eta \left( \frac{L + L^T}{2}\right) = \eta D,
\end{equation}
where  the effective  viscosity 
 $\eta$ is chosen as 
\begin{equation}
\eta= \eta(|D|,p) 
\end{equation}
with  
\[
\eta(|D|,p) = \frac{\mu(I) p \sqrt{2}}{|D|} \quad \mbox{and} \quad I = \sqrt{2}|D|  \frac{d \sqrt{\rho_s}}{\sqrt{p}} .
\]

As before  $\mu(I)$ is given by 
\begin{equation}
\mu(I) = \mu_s + \frac{\mu_2 - \mu_s}{I_0/I + 1}, 
\end{equation}
where $\mu_s = \mu_s(p), \mu_2, I_0$ are determined from experiments.
An important property of this constitutive law is that the effective viscosity $\eta$ diverges to infinity when the shear rate $|D|$ approaches  zero, since in this case $I$ goes to zero and $\mu(I)$ goes to $\mu_s$.
Note that we have as before
\[
\bar \tau = \frac{|\sigma_0| }{\sqrt{2}}=  \frac{\eta(|D|,p)}{\sqrt{2}} \vert D \vert =   \mu(I) p .
\] 

 \begin{remark}
 \label{comp}
 
 Compressible model with $\mu(I)$-rheology.
Based on the compressible balance laws
\begin{eqnarray}	\dot{\bf x} &=& {\bf v}\\
	\dot \rho  &=& - \rho \nabla\cdot {\bf v}  \nonumber \\
	\dot {\bf v} &=& \frac{1}{\rho} \nabla \cdot \sigma +  {\bf g} \nonumber	
\end{eqnarray}
one can derive a compressible plasticity model with $\mu(I)$ rheology, see \cite{BSSG}.
In this case 
$$
\sigma = -p \mathbb I + \eta D_0,
$$
where 
$$ D_0 = \frac{1}{2} (L+L^T) - \frac{1}{d} \; \mbox{div}\; {\bf v} \;  \mathbb I $$
and $p =p(\rho,L) $ and $\eta (\rho,L) $ are depending on $\rho$ and $L$    by a nonlinear relation.

In \cite{BSSG} it is also shown that the  above incompressible equations can be  obtained  with a suitable limit procedure.
\end{remark}


\section{The microscopic model}
\label{Lac}
For comparison we implemented a   non-linear Hertz-Mindlin model as used in  \cite{LPK08} with parameters defined by the  parameters of the macroscopic equations. With the Heaviside function $H$ and the interaction force ${\bf F} = {\bf F}({\bf x}, {\bf v})$  the model is given by a spring-damper model of the following form.
For $i = 1, \ldots, N$ we have 
\begin{eqnarray}
\dot{\bf x}_i &=& {\bf v}_i\\
\dot{\bf v}_i &=& \frac{1}{m} \sum_{j \neq i} {\bf F}({\bf x}_i - {\bf x}_j, {\bf v}_i - {\bf v} _j)  + {\bf g}
\nonumber,
\end{eqnarray}
where $N$ is the total number granular particles and ${\bf F} $ is the interaction force $${\bf F} = H( R_i+R_j - \Vert {\bf x}_i - {\bf x}_j \Vert ) \left( F_n {\bf n} + F_t {\bf t} \right)$$ where $R_i$, $R_j$ are radii of granulars ${\bf x}_i$ and ${\bf x}_j$ respectively,  ${\bf n}$ is the normal vector  given by ${\bf n}= \frac{({\bf x}_i - {\bf x}_j)}{\Vert{\bf x}_i - {\bf x} _j\Vert} $ 
and $ {\bf t}$  the tangential vector  given by ${\bf t}= \frac{ {\bf v}_t}{\Vert {\bf v}_t \Vert } $, where ${\bf v}_t = {\bf v} - <{\bf v}, {\bf n}> {\bf n}$. 
The normal force is decomposed into elastic force $F_n^e$ and dissipative force $F_n^{diss}$ and given by 
\begin{equation}
F_n = F_n^e - F_n^{diss} 
\label{normal_coef}
\end{equation}
with
\begin{equation}
F_n^e = k_n\left(R_i + R_j - \Vert {\bf x}_i - {\bf x}_j\Vert \right)^{3/2} 
\end{equation}
and
\begin{equation}
F_n^{diss} = \gamma_n (R_i + R_j - \Vert{\bf x}_i - {\bf x}_j\Vert)^{1/4}({\bf v}_i - {\bf v}_j)\cdot {\bf n}.
\end{equation}
Here $$k_n =  \frac{4}{3} \frac{G}{1-\nu} \sqrt{\frac{R}{2}}$$ with $G = \frac{E}{2(1+\nu)}$ and $R = R_i $ for all $i = 1, \ldots, N$. 
The dissipative coefficient $\gamma_n$  is given by 
\[
\gamma_n = \beta(e) \sqrt{\frac{5}{2} m k_n} , 
\]
where $\beta(e) = -\frac{\log(e)}{\sqrt{(\log(e)^2 + \pi}}$ and $m$ is the mass of a granular particle. Here $e$ is the coefficient of restitution chosen as  $e = 0.5$.

The tangential force is given as 
\begin{eqnarray} 
F_t =
\left\{ 
\begin{array}{l}  
F_t ^e - F_t^{diss}, 
\quad \mbox{if    }  |F_t^e - F_t^{diss} |  \le \mu F_n
\\ \nonumber
\mu F_n,  \qquad \qquad \mbox{if} |F_t^e - F_t^{diss} |  > \mu F_n.
\end{array}
\right.
\label{tangential_coef}
\end{eqnarray}
Here, the Coulomb friction is $\mu $ =  $0.3815 kg/s$ and 
$$F_t^e = k_t (R_i + R_j - \Vert x_i - x_j\Vert)^{1/2} \delta^t$$ and 
$$F_t^{diss} = \gamma_t (R_i + R_j - \|{\bf x}_i - {\bf x}_j\|)^{1/4}({\bf v}_i - {\bf v}_j)\cdot {\bf t}.$$
 $\delta^t$ is the tangential displacement vector
$$
\delta^t =t_c  ({\bf v}_i - {\bf v}_j) \cdot {\bf t},
$$
where $t_c$ is a constant related to the contact time.

The tangential spring constant $k_t$ and the dissipative constant $\gamma_t$ are given by 
$$k_t =  4\frac{G}{2-\nu}  \sqrt{\frac{R}{2}}, \; \; \gamma_t =2 \beta(e) \sqrt{\frac{5}{12}m  k_t}.$$

\subsection{Rigid body motion}
\label{rigid}
For our second example, the above equations have to be coupled to rigid body motion. We consider the forces acting on the rigid body in case of continuum and microscopic model.

\subsubsection{Rigid body motion  for the continuum model}
The rigid body motion is given by the Newton-Euler equations:
\begin{equation}
M\frac{d { \bf V} }{dt} = {\mathcal{ {  F} } } + M {\bf g} , \; \;
[I]\cdot \frac{d  { \bf \omega} }  {dt} + {{ \bf \omega} }\times ({[I]\cdot}{ { \bf \omega} }) =  {\mathcal{ { T} } },
\label{euler_newton}
\end{equation}
where $M$ is the total mass of the body with center of mass ${ \bf X}_c$, ${ \bf V} $ is the translational velocity and $ {\bf \omega}$ is the angular velocity of the rigid body. 
$\mathcal{ {  F} }$ is the translation force, $ \mathcal{ { T} } $ is the torque  and $[I]$ is the moment of inertia.  
The center of mass of the rigid body is obtained by  
\begin{equation}
\frac{d { \bf X}_c}{dt} = { \bf V}.   
\label{pos_velo_rigid}
\end{equation}
Finally, the velocity of a point on the surface $\partial S$ of the rigid body is given by   $ {{\bf V}_S} ={{ \bf V}}+{ {\bf  \omega} } \times({ { \bf x} } -{ {\bf X}_c}) ,~{ { \bf x} }\in \partial S$.

The force $\mathcal{ { F} }$ and torque ${  \mathcal{ { T} } }$, that the granular particles exert on the  rigid body, 
are computed according to 
\begin{equation}
\mathcal{ { F }} = \int_{\partial S} (-\sigma \cdot{ { \bf n} }_{s}) dA, \; \;
\mathcal{ { T }}= \int_{\partial S} ( { \bf x} - {\bf X}_c)\times(-\sigma \cdot{ { \bf n} }_{s}) dA,
\label{force_torque_bgk}
\end{equation}
where $\sigma$ is the Cauchy stress tensor computed from the continuum model and ${\bf n}_s$ is the unit normal vector on the surface pointing towards the granular materials. 


\subsubsection{Rigid body motion  for the microscopic model}
In this case we compute the forces acting on the   rigid body from its neighboring sand particles. On the boundaries of the  rigid body, we generate  
boundary particles.
The motion of the rigid body is then obtained by solving
\begin{eqnarray}
	M\frac{d {\bf V}}{dt} &=& \sum_{i}  \sum_{j} {{ { \bf F} } } ({\bf x}_i - {\bf x}_j, {\bf v}_i - {\bf v}_j)  + M \bf{g},
\end{eqnarray}
for all boundary particles ${\bf x}_i$ on the rigid body.  Note that ${\bf x}_j$ are the neighboring sand particles of boundary particle ${\bf x}_i$, where the boundary particles are excluded from the neighbor list. 

The force ${{ \bf F} } $ is computed as for the microscopic interaction model described in the previous section.


\section{Numerical algorithms}

We consider the continuum equations and discuss space and time discretization of the meshfree particle method and the computation of the stress tensor for the model
in Section \ref{macroeq}.
The numerical treatment of  the  incompressible continuum model is discussed in a short remark.

\subsection{Space and time discretization}

Let $\Omega$ be the computational domain with boundary $\Gamma$. We first approximate 
the boundary $\Gamma$ of the domain by a set of discrete points called  boundary particles.   In a second step we approximate the 
interior of the computational domain using an second  set of  so called interior points or interior particles. Note that the particles used  for solving the continuum models are not real particles in the sense of a microscopic models, but   numerical grid points, which move with the flow  and carry all necessary information like density and  velocity. 
 Let $x_i, i=1, \ldots N$ be the grid points with $N$  the total number of grid points.  The discrete version of  the system of equations in Lagrangian form (\ref{lagrange}) is 
\begin{eqnarray}
\label{char1}
\dot{\bf x}_i &=& {\bf v}_i \\
\dot \rho_i  &=& - \rho_i (\nabla\cdot {\bf v})_i  \nonumber \\
\dot {\bf v}_i &=& \frac{1}{\rho_i} (\nabla \cdot \sigma)_i +  {\bf g} \nonumber  \label{momentum1},
\end{eqnarray}
for $i = 1, \ldots, N$. 

Approximating the spatial derivatives on the right hand sides of system of equations (\ref{char1})
at each point ${\bf x}_i$,
we obtain a  system of Ordinary Differential Equations. The ODE system can be solved by any standard ODE solver. We simply use an  explicit Euler method for the time integration. 
The spatial derivatives at each point are approximated using the values on  its surrounding cloud of points and a  weighted least squares method, see the next subsection.
The computation of the  Cauchy stress tensor $\sigma$ is described in subsection \ref{stress}.

\subsection{Least squares approximation of the derivatives}

Consider  the problem of approximating the spatial 
derivatives of a function 
$ f(t, {\bf x}) $  at the particle position ${\bf x}= {\bf x}_i$ in terms of the values of 
its neighboring points ${\bf x}_{j}$.
In order to restrict the number of  points
we  associate a weight 
function $ w_{j} = w(\Vert{\bf x}_{j} - {\bf x}\Vert; h)) $  with
small compact support, where  the smoothing length
$ h $ determines the size of the support.
We consider 
a Gaussian weight function depending on the distance between the central particle ${\bf x} ={\bf x}_i$ and its neighbor ${\bf x}_{j}$ in the following form  
\begin{eqnarray} 
w(\Vert{\bf x}_{j} - {\bf x}\Vert; h) =
\left\{ 
\begin{array}{l}  
exp (- \alpha \frac{\Vert {\bf x}_{j} - {\bf x}  \Vert^2 }{h^2} ), 
\quad \mbox{if    }  \frac{\Vert {\bf x}_{j}- {\bf x}  \Vert}{h} \le 1 
\\ 
0,  \qquad \qquad \mbox{else},
\end{array}
\right.
\label{tiwweight}
\end{eqnarray}
with $ \alpha $ a user defined positive constant
chosen here as 
$\alpha = 6.25$.  
The size of the support $h$ defines
a set of neighboring particles  
$ P({\bf x}, h) = \{ {\bf x}_{j} : j=1,2,\ldots,m \} $.
We choose  the constant $h$ equal to $3$ times the initial spacing of the grid points.

In the above presented equations and the construction of the Cauchy tensor only the first order partial derivatives are involved. Thus, we use a first order Taylor expansion of $f({\bf x}_{j})$ around ${\bf x}={\bf x}_i $ 
\begin{equation}
 f({\bf x}_{j}) =  f({\bf x}) + ({\bf x}_{j}- {\bf x}) \nabla f ({\bf x}) 
  + e_{j}, 
\label{taylor}
\end{equation}
for $j=1, \ldots, m$, where $e_{j}$ is the error in the  expansion.   
The unknown ${\bf a} = \nabla f({\bf x})$ is now computed by minimizing the errors $e_{j}$ for 
$j=1, \ldots, m$.  
The system of equations can be re-written  as 
\begin{equation}
\label{sys}
{ \bf e } = { \bf b} - H {\bf a}, 
\end{equation}
where ${\bf e} = [e_{1}, \ldots, e_{m}]^T$,  ${\bf b} = [f_{1} - f, \ldots, f_{m} - f]^T $ 
and
\begin{eqnarray}
H = \left( \begin{array}{ccc}
dx_1 &  ~dy_1 & ~ dz_1   \\
\vdots  &\vdots & \vdots   \\
dx_m &  ~dy_m  & ~dz_m  
 \end{array} \right)
\label{matrixM}
\end{eqnarray}
with $dx_j = x_{j} - x, \;  dy_j = y_{j}-y, \; dz_j = z_{j} - z$. 

Imposing ${\bf e} =0 $ in (\ref{sys}) results in an overdetermined linear system of algebraic equations, which in general has no solution. The unknown
${\bf a}$ is therefore obtained from the weighted least squares method by minimizing the quadratic form 
\begin{equation}
J = \sum_{j=1}^m w_{j} e_{j}^2 = (H {\bf a}  - {\bf  b})^T W (H {\bf a} - {\bf b}),
\end{equation}
where 
\begin{eqnarray*}
W=\left( \begin{array}{cccc}
w_1 & 0 & \cdots& 0 \\
\vdots & \vdots & \cdots & \vdots \\
0 & 0 & \cdots & w_m  
\end{array} \right).
\end{eqnarray*}

The minimization of $J$ formally yields 
\begin{equation}
{\bf a} = (H^TWH)^{-1}(H^TW){\bf b}. 
\end{equation}
For further details on the spatial discretization and higher order approximation of derivatives, we refer to \cite{SKT,TKR,TKH}.

\begin{remark}{Numerical algorithm for the simplified macroscopic model.}
	The spatial discretization of the  incompressible model  proceeds similarly. 
	 For the treatment of the incompressibility constraint, we refer to \cite{TK07, DTKB08, TKH}.
\end{remark}

\begin{remark}
When grid points or particles move, after some time steps they may form a close cluster or may disperse from each other. If two points are getting very close to each other, we remove them and replace a single new point in the middle of them. If particles  disperse from each other, they form  holes in the computational domain. In such a case there are not enough  neighboring points in order to approximate the derivatives and the numerical scheme becomes unstable. Therefore, we have to add points in such a case.  At the  newly created points one has to interpolate the flow quantities from their  nearest neighboring points. For details, we refer  \cite{TKR}.

\end{remark}

\subsection{Stress computation}
\label{stress}
To compute the stress for the model in Subsection \ref{Kamrin} we  closely follow the procedure developed in \cite{DK15}. We  include the computations for the sake of completeness.
$\sigma$ is numerically found by the following algorithm.
Given quantities are  the constants $\mu_s$, $\mu_2$, $\xi,I_0$ and $\nu, E, \rho_c$.
One time step of the  so called elastic trial step algorithm is given by the following.

Suppose  $\rho^n, v^n,L^n, \sigma^n$ are known from the previous time step.
Then

\begin{itemize}
	\item Determine $\rho^{n+1},v^{n+1} $ from an explicit solve of the balance equations (\ref{lagrange}).
	\item Determine $L^{n+1}$ and $D^{n+1},W^{n+1}$ from (\ref{L}) and (\ref{DW}).
	\item Solve (\ref{ode}) using $\sigma^n, D^{n+1},W^{n+1}$ with $D_p =0$ for one time step $\Delta t $ with the explicit Euler scheme.
	The result is called $\sigma^{tr}$.
	\item 
	This gives $p^{tr}$ and $\sigma_0^{tr}$ and $\bar \tau^{tr}$.
	\item 
	If $\rho^{n+1}< \rho_C$ or $p^{tr} \le 0$, set  $\sigma^{n+1} =0$.
\end{itemize}

\begin{itemize}
	\item
	If  $\rho^{n+1}> \rho_C$ and $p^{tr} > 0$,
	determine an explicit expression for  $\dot {\bar \gamma}$ depending on  $\bar \tau^{n+1}$ and $p^{tr}$ using (\ref{tau}). Determine an explicit expression for  $D_p^{tr}$ depending on  $\dot {\bar \gamma}$ and $\sigma_0^{n+1}$ using (\ref{dp}).
	Set $p^{n+1} = p^{tr}$.
\end{itemize}

Now we use the fact that 
\begin{equation}
\sigma^{n+1} - \sigma^n = \sigma^{tr} - \sigma^n - \Delta t \frac{E}{1+\nu} \left[    D_p^{tr} + \frac{\nu}{1- 2 \nu} \mbox{tr} D_p^{tr} \mathbb I \right]
\end{equation}
or 
\begin{equation}
\sigma^{n+1}  = \sigma^{tr} - \Delta t 2 G D_p^{tr}
\end{equation}
with
$G = \frac{E}{2(1+\nu)}$
and therefore
\begin{equation}
\sigma^{n+1}  - p^{tr}= \sigma^{tr} -p^{tr}- \Delta t 2 G \frac{1}{\sqrt{2}}  \dot {\bar{ \gamma} } (\bar \tau^{n+1})\frac{\sigma_0^{n+1}}{\vert \sigma_0^{n+1} \vert}.
\end{equation}
This  gives 
\begin{equation}
\sigma_0^{n+1}  \left(1+ \Delta t 2 G \frac{1}{\sqrt{2}}  \dot {\bar{ \gamma} } (\bar \tau^{n+1})\frac{1}{\vert \sigma_0^{n+1} \vert}
\right) = \sigma_0^{tr} 
\end{equation}

Taking absolute values this leads to 
\begin{equation}
\vert \sigma_0^{n+1} \vert = \vert \sigma_0^{tr} \vert - \sqrt{2} G \Delta t \dot{\bar \gamma}
\end{equation}
or
\begin{align}
\label{1}
\bar \tau^{n+1} = \bar \tau^{tr} - G \Delta t \dot{\bar \gamma}
\end{align}
and
\begin{align}
\label{22}
\frac{\sigma_0^{n+1}}{\sigma_0^{tr}}  = \frac{\vert\sigma_0^{n+1} \vert  }{\vert\sigma_0^{tr} \vert } =  \frac{\bar \tau^{n+1}   }{\bar \tau^{tr}  } 
\end{align}
or
\begin{align}
\label{2}
\frac{\sigma^{n+1}+p^{tr}}{\sigma_0^{tr}}  = \frac{\vert\sigma_0^{n+1} \vert  }{\vert\sigma_0^{tr} \vert } =  \frac{\bar \tau^{n+1}   }{\bar \tau^{tr}  } .
\end{align}

These considerations yield the final step of the algorithm as 
\begin{itemize}
	\item
	Now compute  $\sigma^{n+1}$ as  $\sigma^{tr}$ in case $\bar \tau^{tr} > \mu_s  p^{tr}$. Moreover, determine first $ \bar \tau^{n+1}$
	from (\ref{1})  using the explicit expression for $\dot{\bar \gamma} $ and then $\sigma^{n+1}$ from (\ref{2}),
	if $\bar \tau^{tr} \le \mu_s  p^{tr}$ .
\end{itemize}

Note that 
equation (\ref{1}) is equivalent to 
\begin{align*}
\label{3}
\bar \tau^{n+1} = \bar \tau^{tr} - G \Delta t  \sqrt{p^{tr}} \xi    \frac{\bar \tau^{n+1}-\mu_s p^{tr} }{\mu_2 p^{tr} - \bar \tau^{n+1}}
\end{align*}	
This is equivalent to 
\begin{align*}
\bar \tau^{n+1} (\mu_2 p^{tr} - \bar \tau^{n+1})= \bar \tau^{tr} (\mu_2 p^{tr} - \bar \tau^{n+1})- G \Delta t  \sqrt{p^{tr}} \xi   (\bar \tau^{n+1}-\mu_s p^{tr})
\end{align*}	
or
\begin{align*}
0 = (\bar \tau^{n+1})^2  - \bar \tau^{n+1} (\mu_2 p^{tr} + \bar \tau^{tr} + G \Delta t  \sqrt{p^{tr}} \xi  ) +
\bar \tau^{tr} \mu_2 p^{tr} + G \Delta t  \sqrt{p^{tr}} \xi  p^{tr} \mu_s
\end{align*}	
or
\begin{align*}
0 = (\bar \tau^{n+1})^2  - \bar \tau^{n+1} B + H.
\end{align*}	
Here, following \cite{DK15}, the negative root has to be chosen, i.e.
\begin{align}
\bar \tau^{n+1} = \frac{1}{2} (B - \sqrt{B^2-4H}) =\frac{2H}{B + \sqrt{B^2-4H}}
\end{align}

\section{Numerical results}

\subsection{Physical situation and parameters}

We consider two types of numerical experiments. First, the collapse of a sand column with different geometries is considered.
This is a classical, well investigated test-case, see, for example  \cite{LH,LSP}.

Second, we consider the more complicated case of a free falling disc into a box filled with sand,
compare \cite{NGM10}.
For this second case the interaction of the free falling disc with the sand particles has to be considered, see subsection  \ref{rigid}.
We compare in both cases the results of the microscopic, the macroscopic and the simplified  model.

In both cases we use the following values for the parameters of the macroscopic equations
$E=  10^9 P$, $\nu=0.3$,  $\xi = 1.123  \sqrt{\frac{m}{kg}}$ (or $I_0= 0.32, d =   3 \cdot 10^{-3} m $),
$\mu_2 = 0.6435$ and $ \mu_s = 0.3819$.
The  material density is given by   $\rho_0 = \Phi \rho_S \approx 0.6 \cdot 2.450$ =   $ 1500 kg m^{-3} $.
Moreover, we use 
${\bf g} = (0, 0,g) m s^{-2} $ with  $g = -9.81  m s^{-2} $.

The above density  leads to a mass $m= \frac{4}{3}\pi R^3 \rho_S =\frac{4}{3}\pi 27 \cdot 10^{-9} 10^3 \cdot 2.45 \sim 3  \cdot 10^{-4}kg  $ which is required for the microscopic equations.
Additionally, we use  for the microscopic model the parameter $\mu = 0.3815$.
%
%
%

\subsection{Collapsing $2D$ sand column}

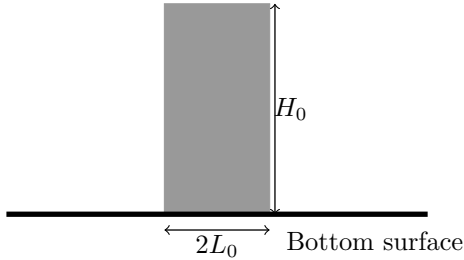
\begin{figure}[ht]
	\begin{center}
		\begin{tikzpicture} [scale=0.7]	      
	        \draw[line width=2pt,color=black,opacity=1] (0,0) -- (8,0) ;
	        \node[color=black] at (7,-0.5) {Bottom surface };
	        \draw [fill=black,opacity=0.4] (3,0) rectangle (5,4);
		\draw[line width=0.5pt,<->] (3,-0.3) -- (5,-0.3);
		\node[color=black] at (4,-0.6) {$2 L_0$ };
		\draw[line width=0.5pt, <->] (5.1,0) -- (5.1,4);
		\node[color=black] at (5.4,2) {$H_0$ };
		\end{tikzpicture}
	\caption{Initial geometry of a sand column of height $H_0$ and width $2 L_0$.}
 \label{Initial_column}
  \end{center}
  \end{figure}

 We considered a two dimensional collapsing column of sand as shown  in Fig. \ref{Initial_column}, where the sand has initial width $2 L_0$ and height $H_0$. The aspect ratio is defined as $a = H_0/L_0$. In our numerical experiments, we have used different aspect rations keeping $L_0$ constant  and varying $H_0$. 

\begin{figure}
	\centering
	\includegraphics[keepaspectratio=true, angle=0, width=0.32\textwidth]{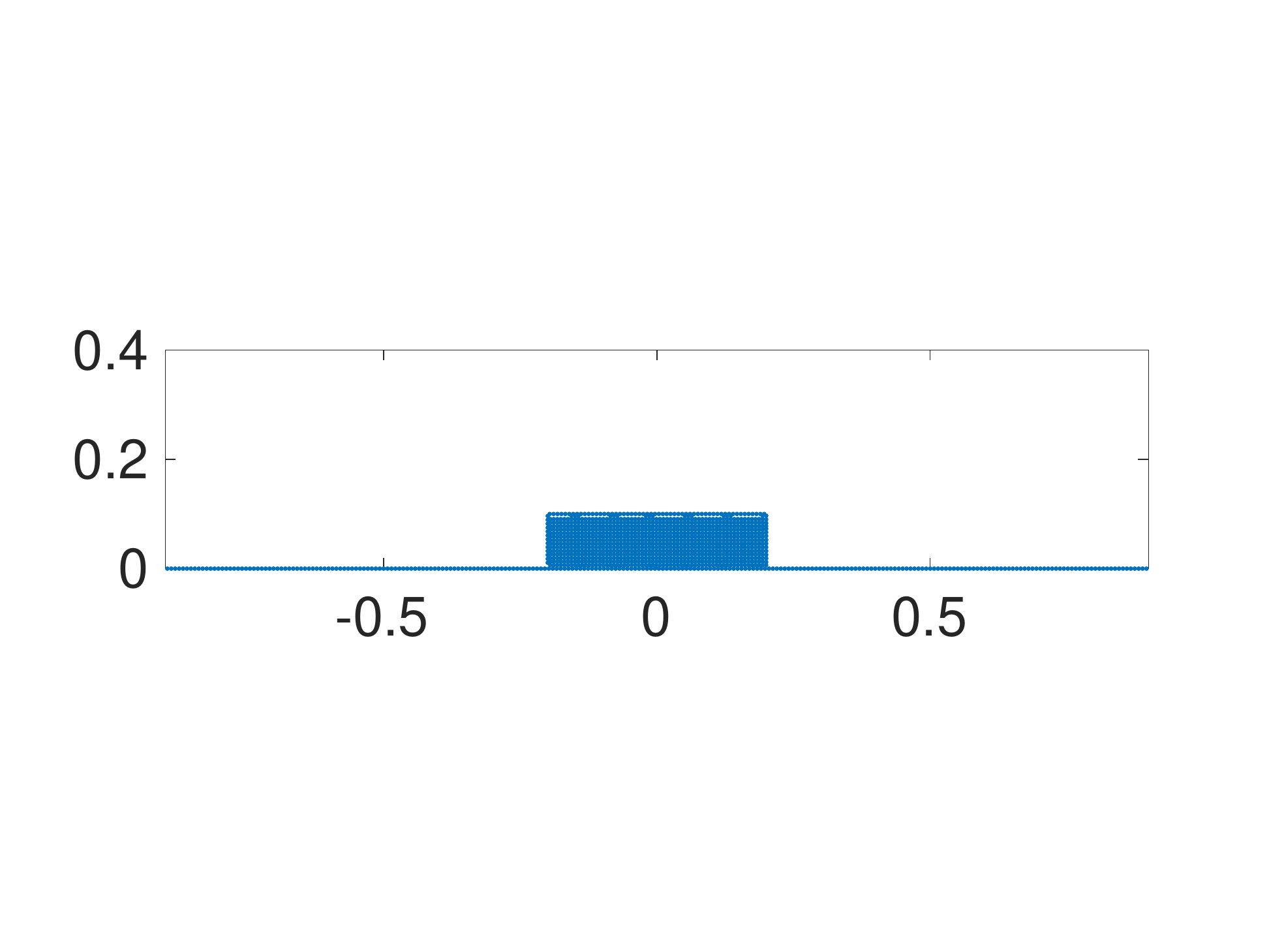} 
	\includegraphics[keepaspectratio=true, angle=0, width=0.32\textwidth]{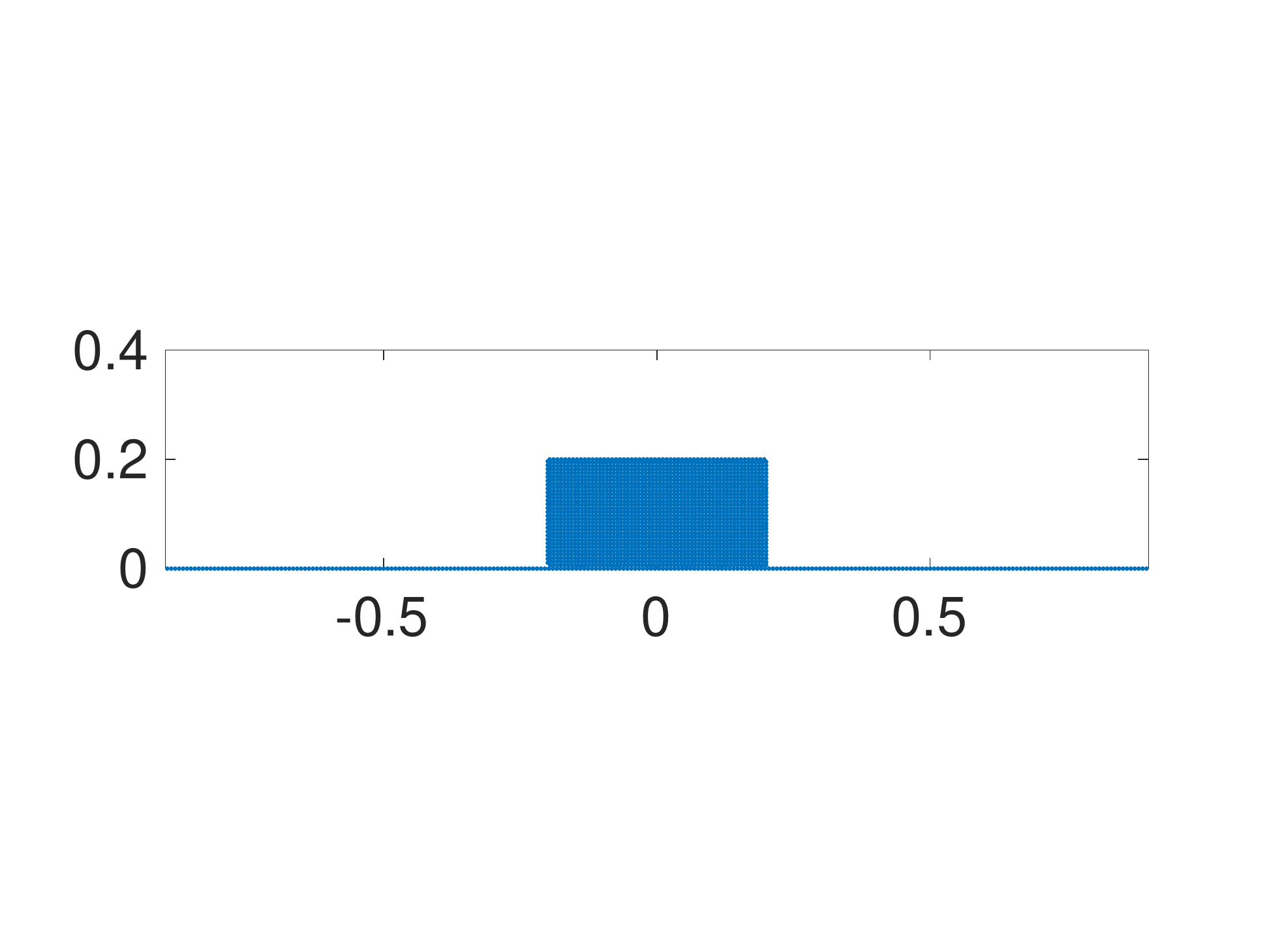} 
	\includegraphics[keepaspectratio=true, angle=0, width=0.32\textwidth]{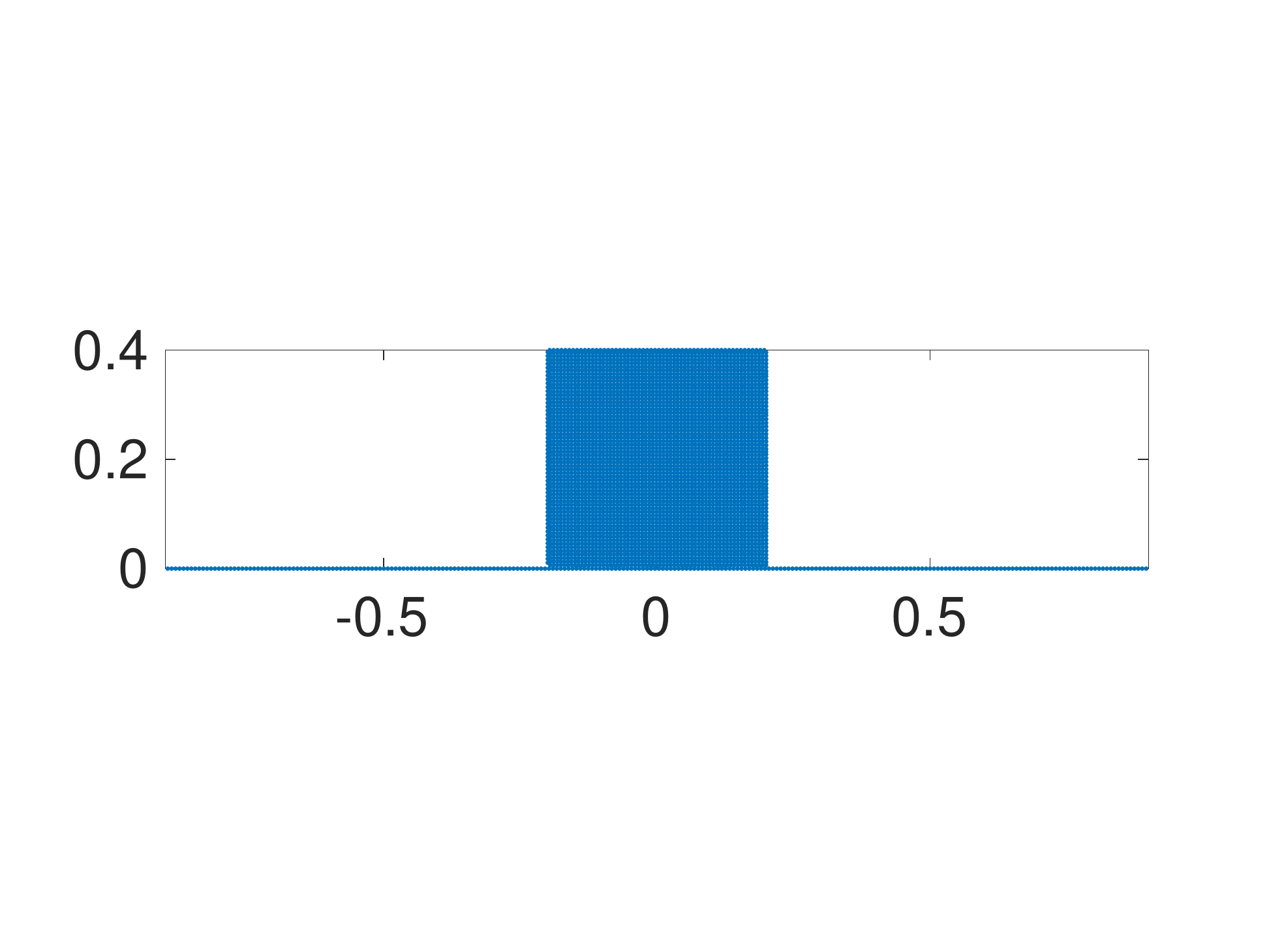} 
	\centering
	\caption{Initial values for $a = 0.5$ (left),  $a = 1$ ( middle), $a = 2$ (right) at $ t = 0$. }	
	\label{column_t0}
\end{figure}

\begin{figure}
	\centering
		\includegraphics[keepaspectratio=true, angle=0, width=0.32\textwidth]{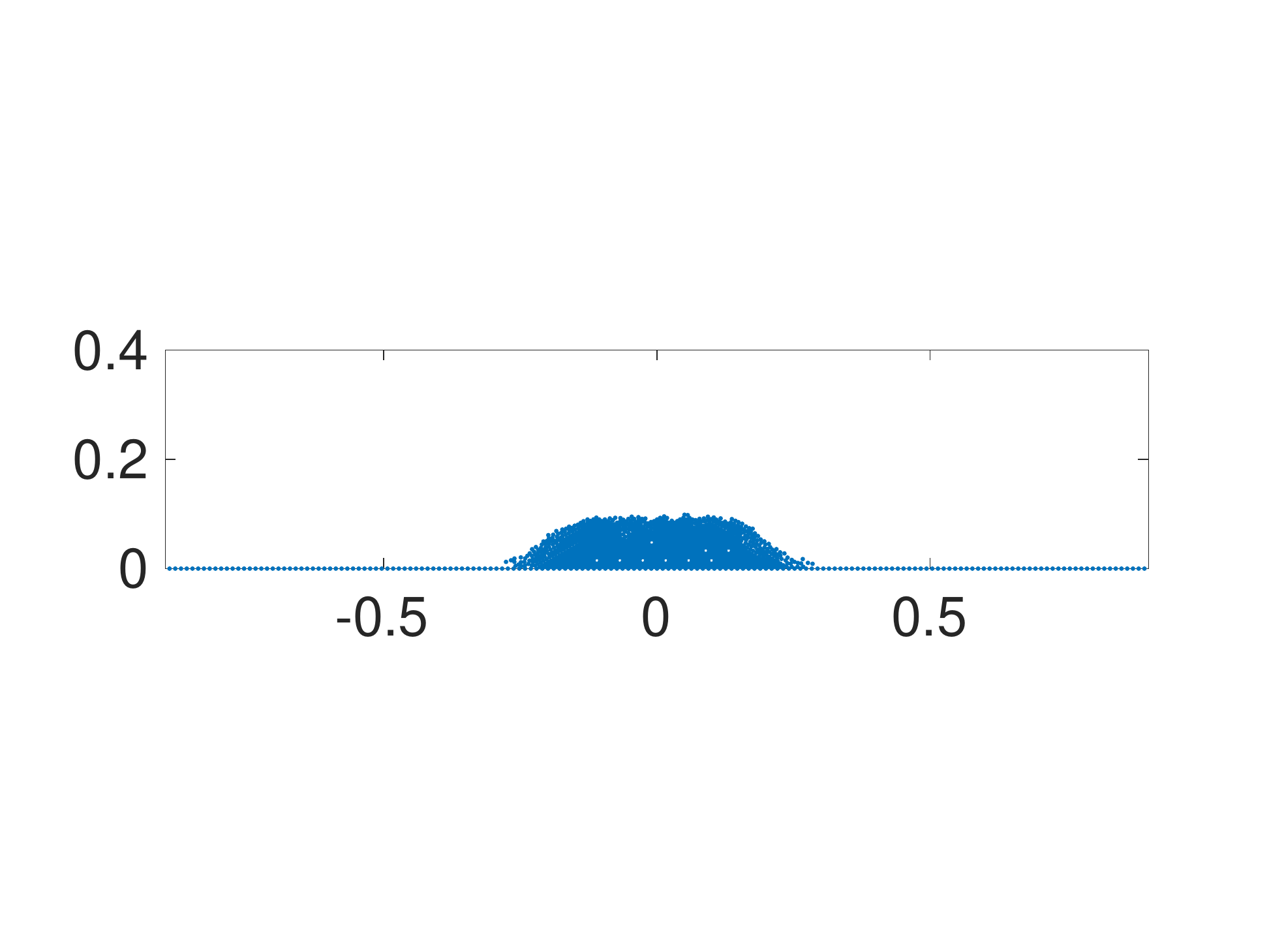} 
	\includegraphics[keepaspectratio=true, angle=0, width=0.32\textwidth]{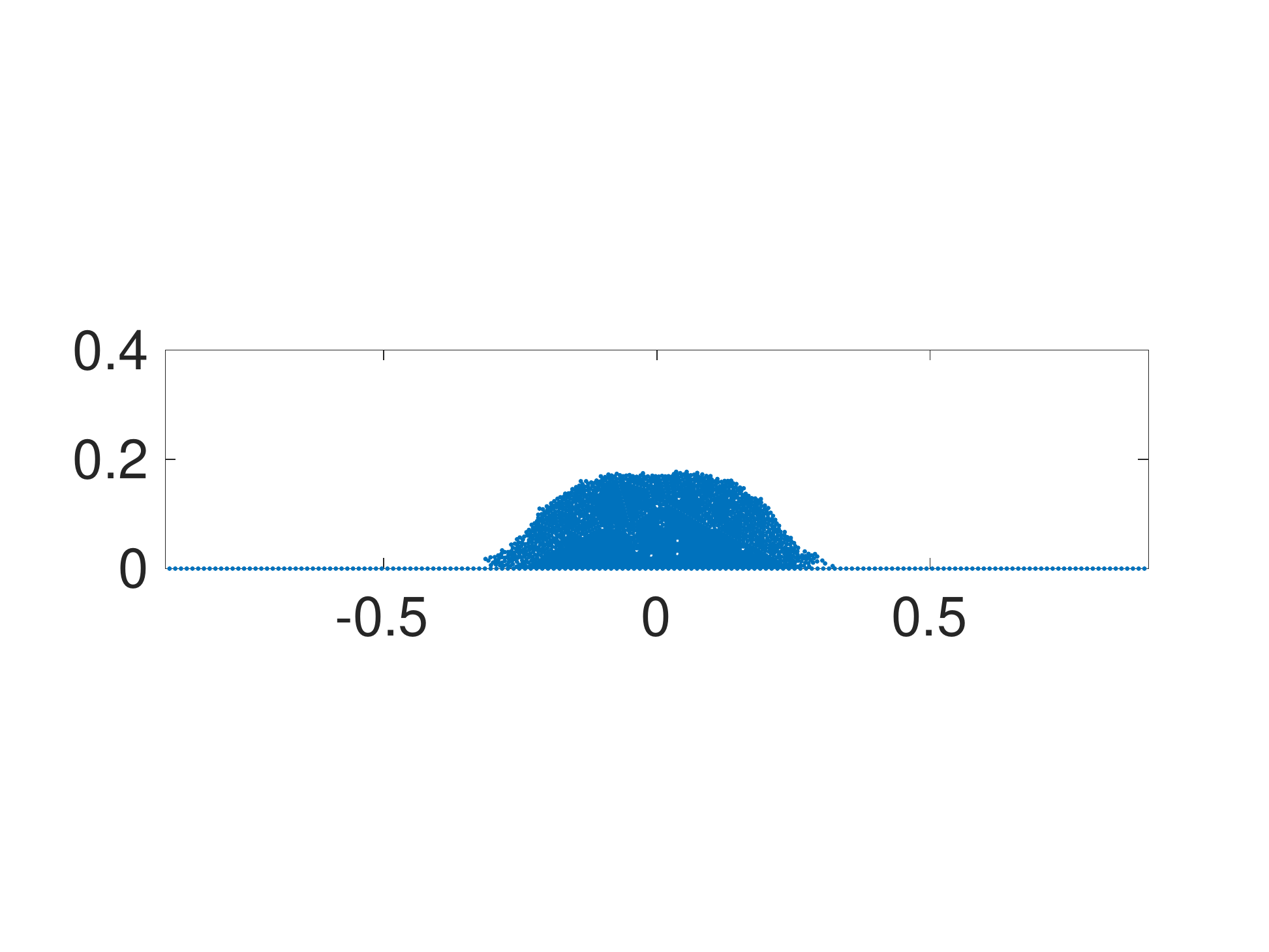} 
	\includegraphics[keepaspectratio=true, angle=0, width=0.32\textwidth]{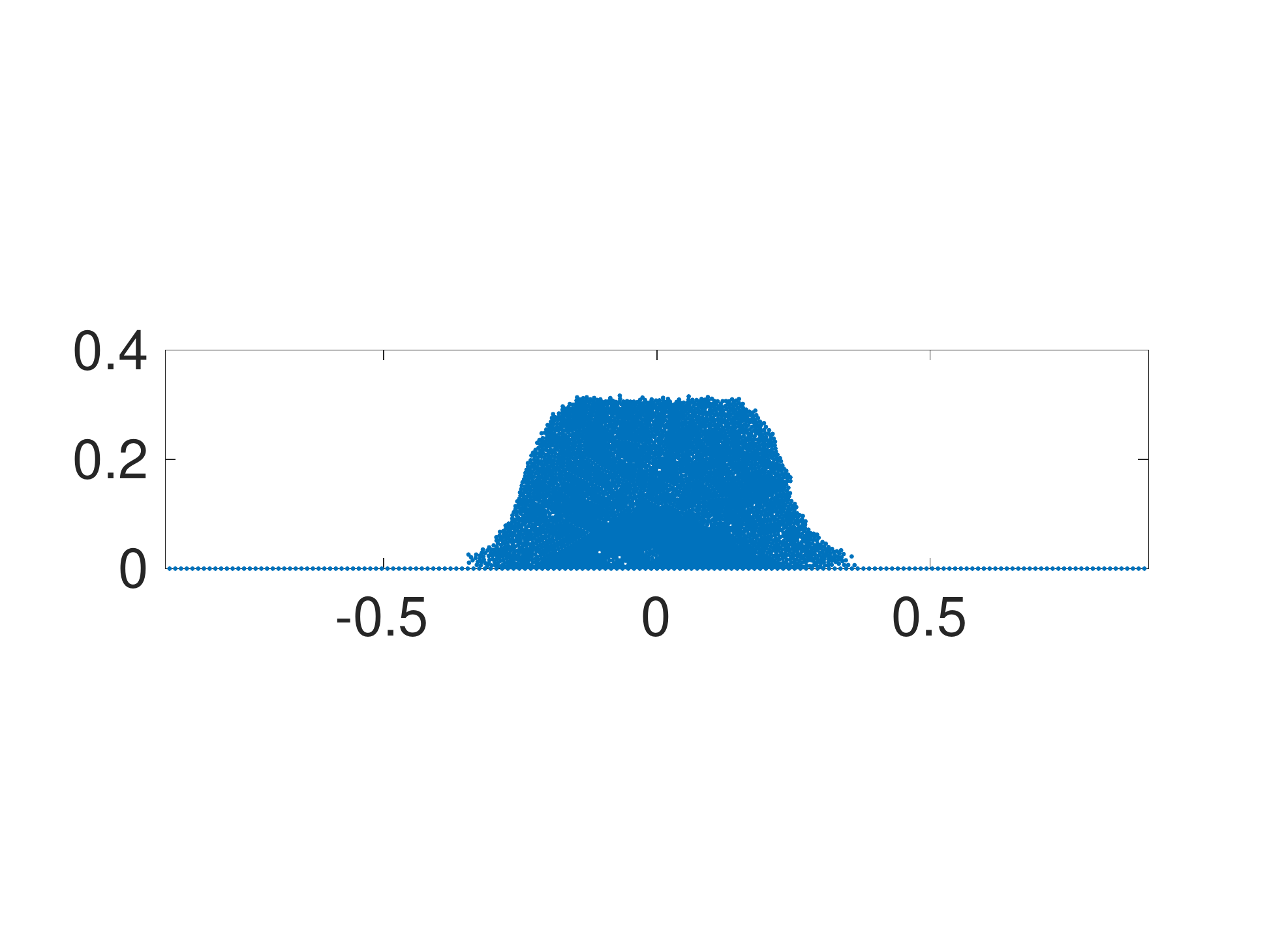} \\
	\vspace{-1cm}
	\includegraphics[keepaspectratio=true, angle=0, width=0.32\textwidth]{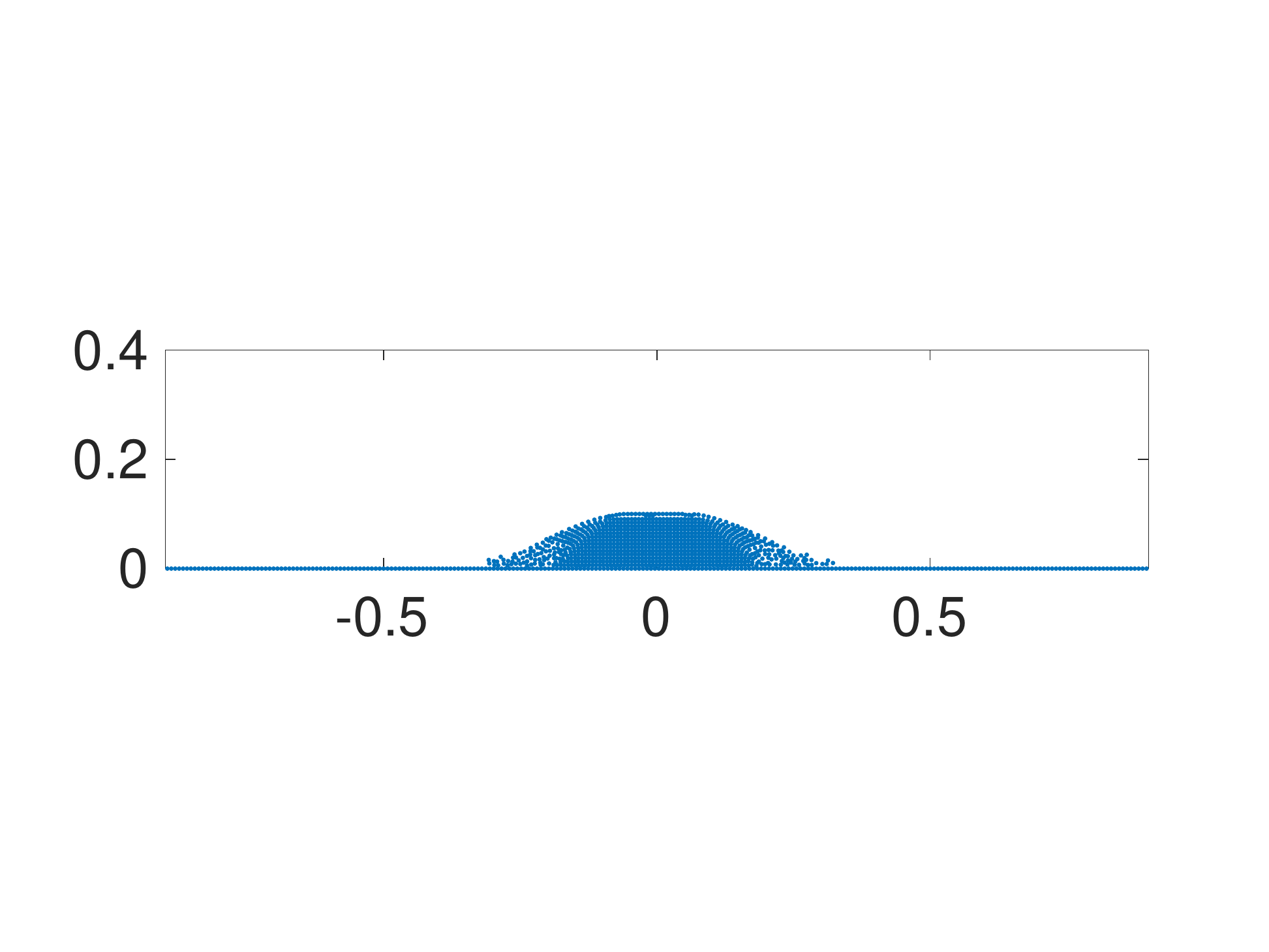} 
	\includegraphics[keepaspectratio=true, angle=0, width=0.32\textwidth]{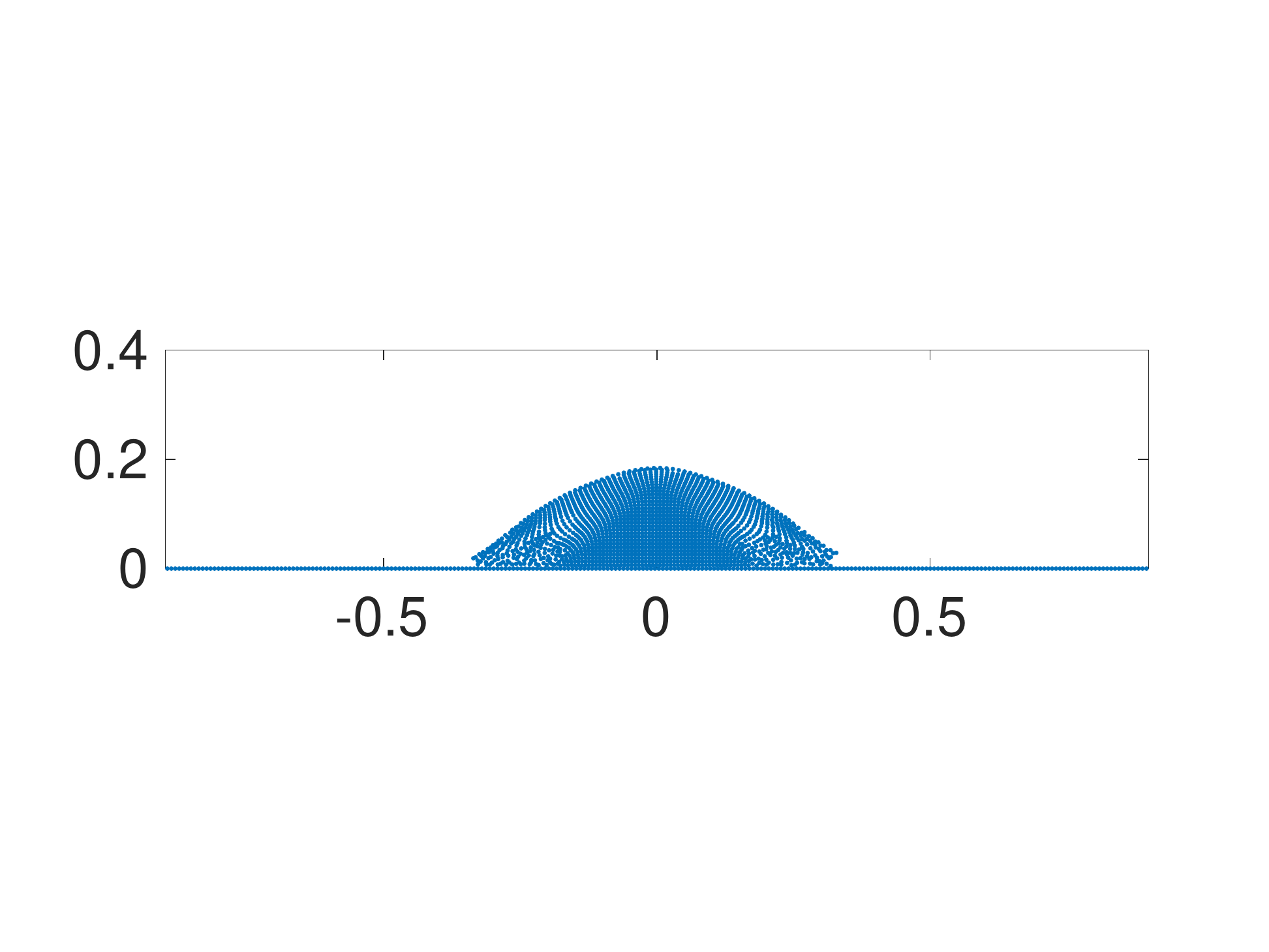} 
	\includegraphics[keepaspectratio=true, angle=0, width=0.32\textwidth]{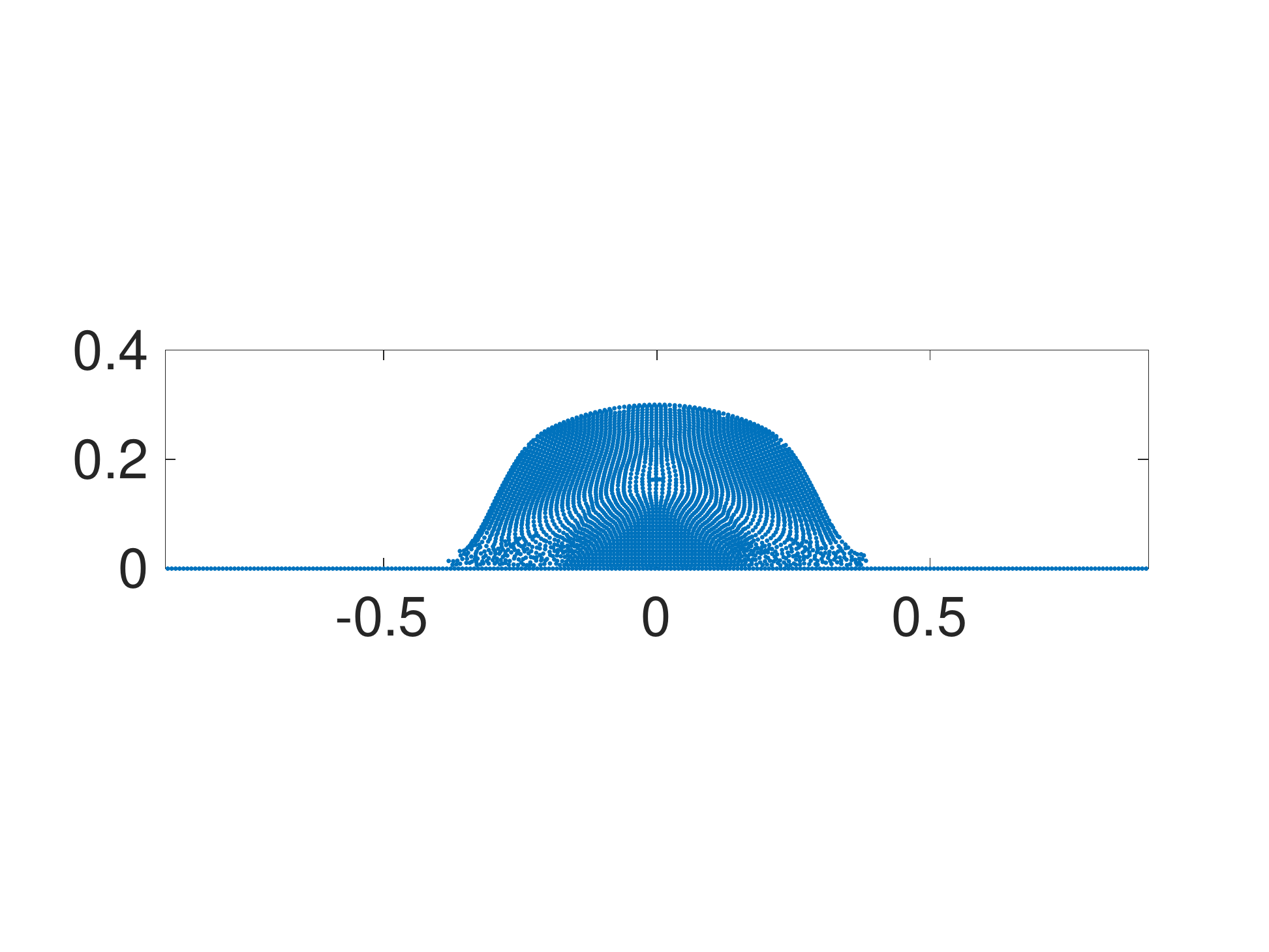} 	\\
	\vspace{-1cm}
	\includegraphics[keepaspectratio=true, angle=0, width=0.32\textwidth]{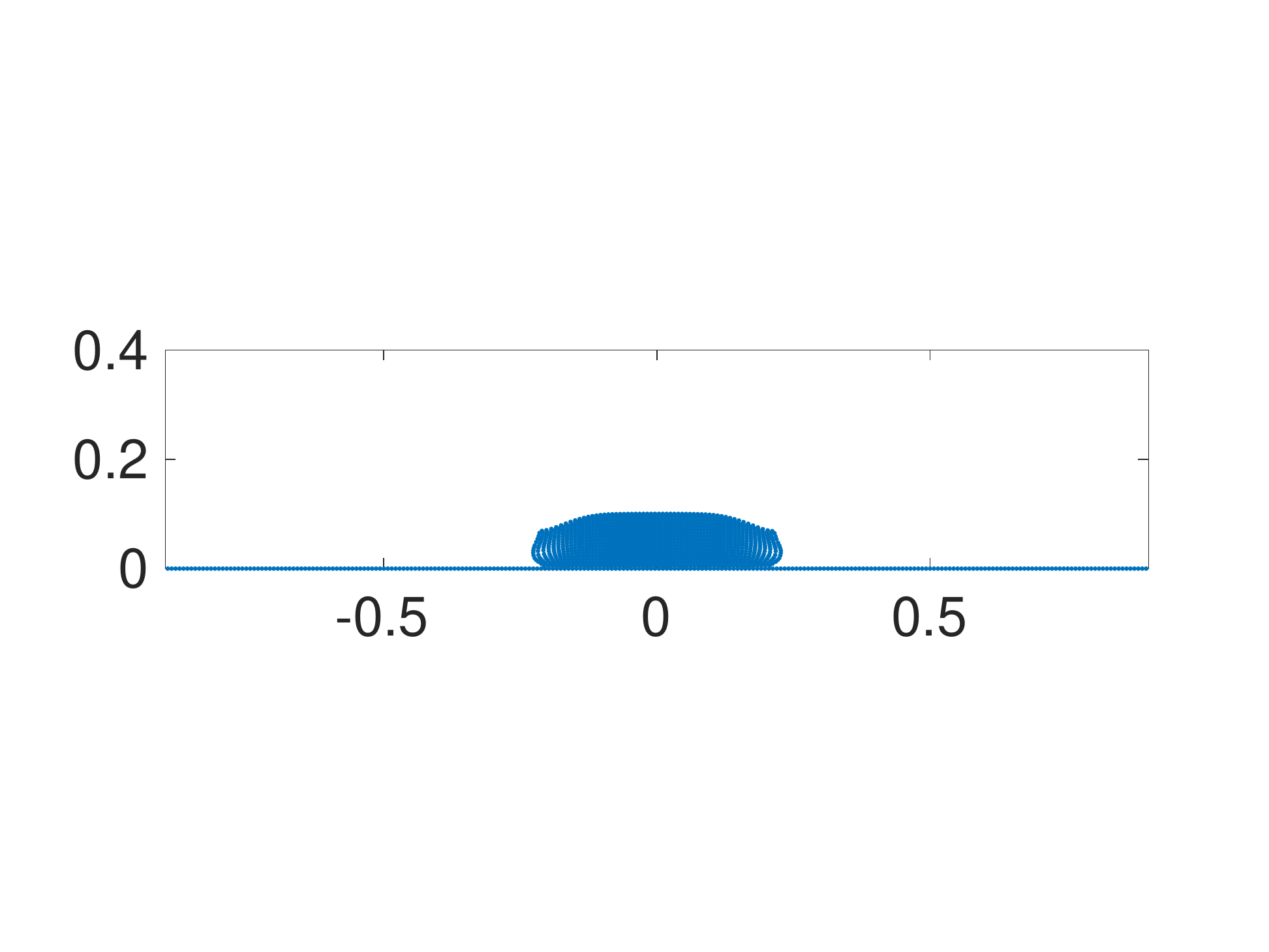} 
	\includegraphics[keepaspectratio=true, angle=0, width=0.32\textwidth]{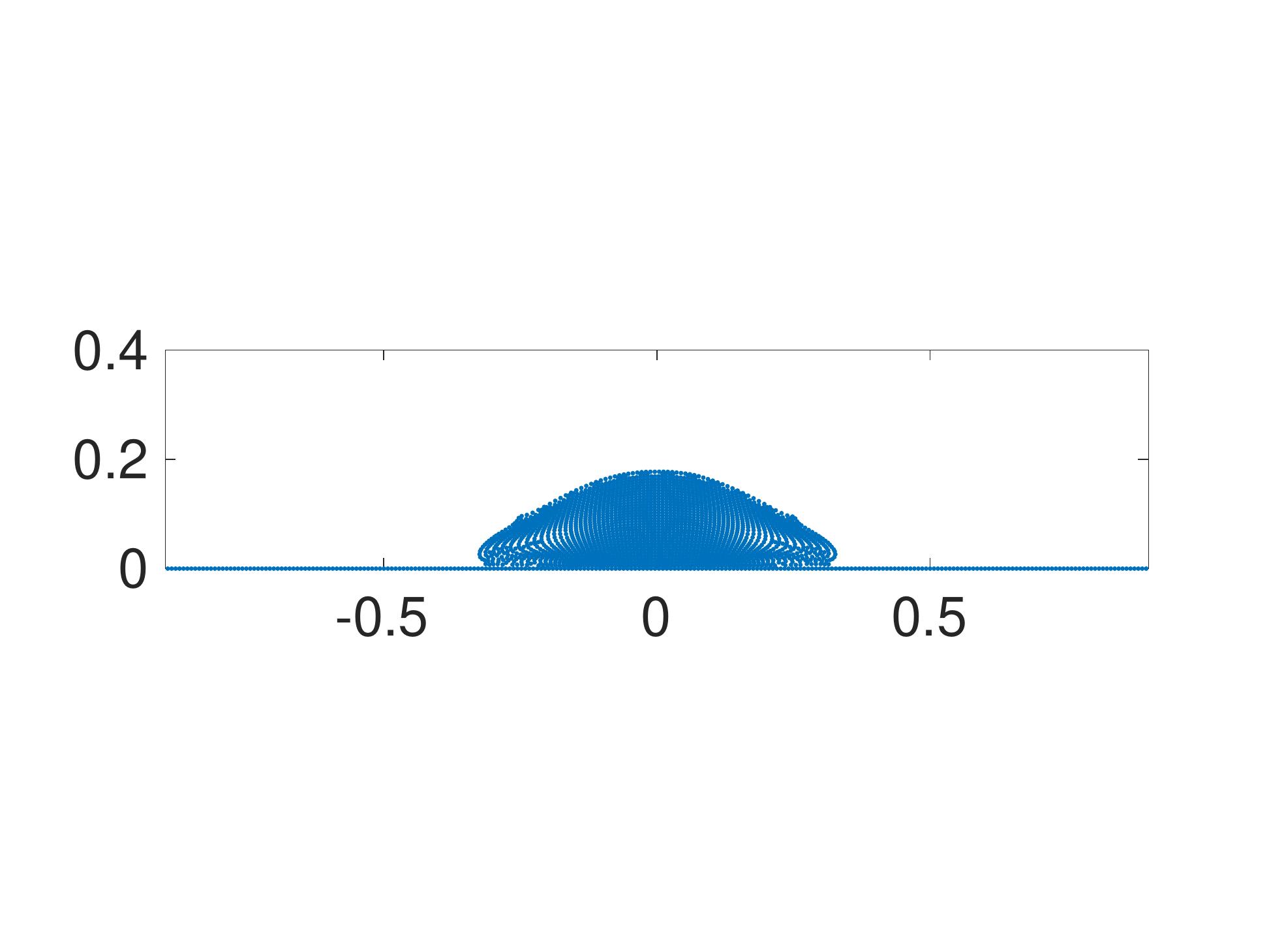} 
	\includegraphics[keepaspectratio=true, angle=0, width=0.32\textwidth]{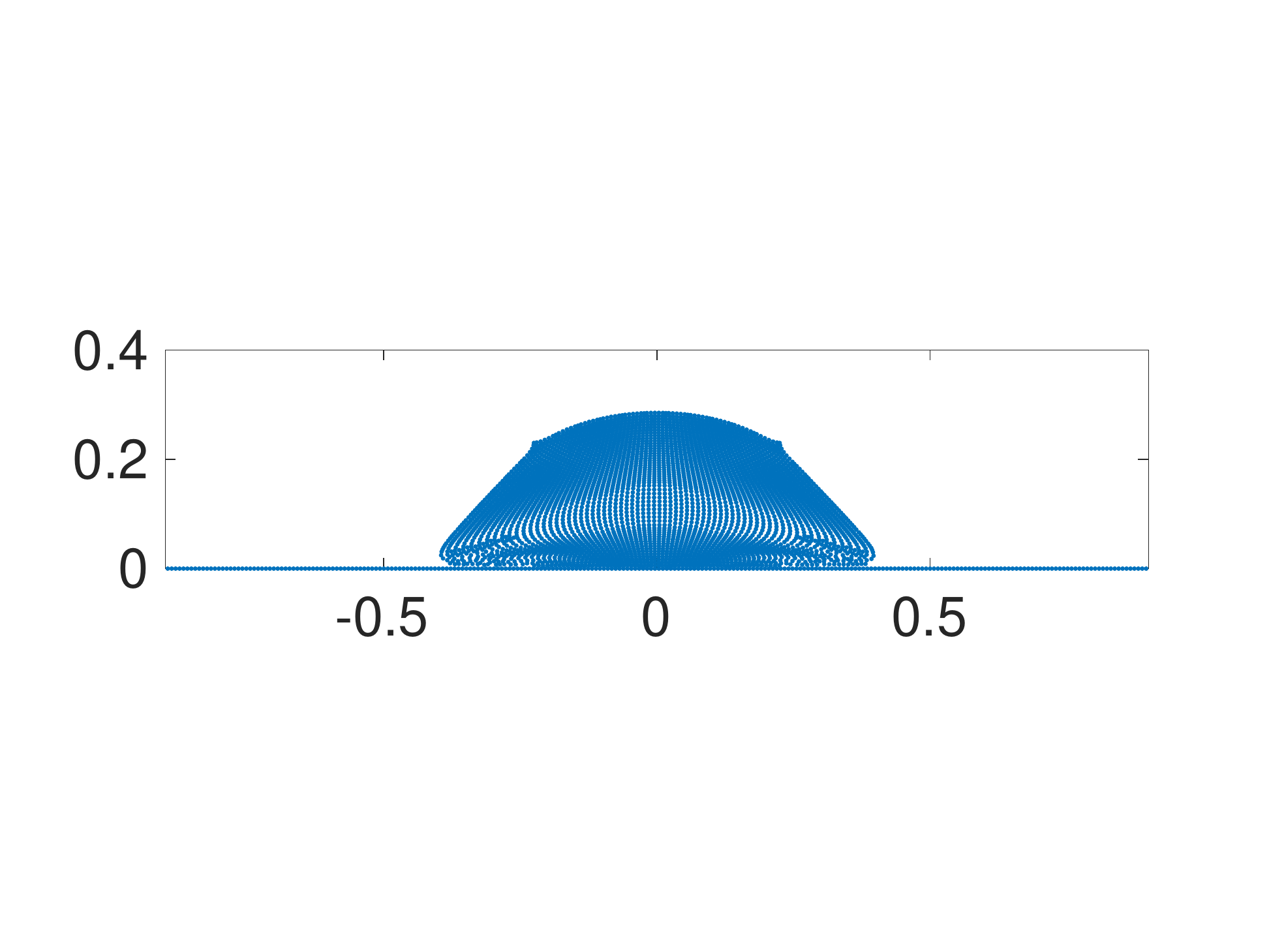} 	
	\caption{Time development at $ t = 0.2$ for $a = 0.5$ (left),  $a = 1$ ( middle), $a = 2$ (right). First row: micoscopic model (Section \ref{Lac}), second row: plasticity model (Section \ref{macroeq}),third row: Coloumb constitutive model (Section \ref{Jop})}	
	\label{column_compare_t0dot2}
	\centering
\end{figure}

\begin{figure}
	\centering
		\includegraphics[keepaspectratio=true, angle=0, width=0.32\textwidth]{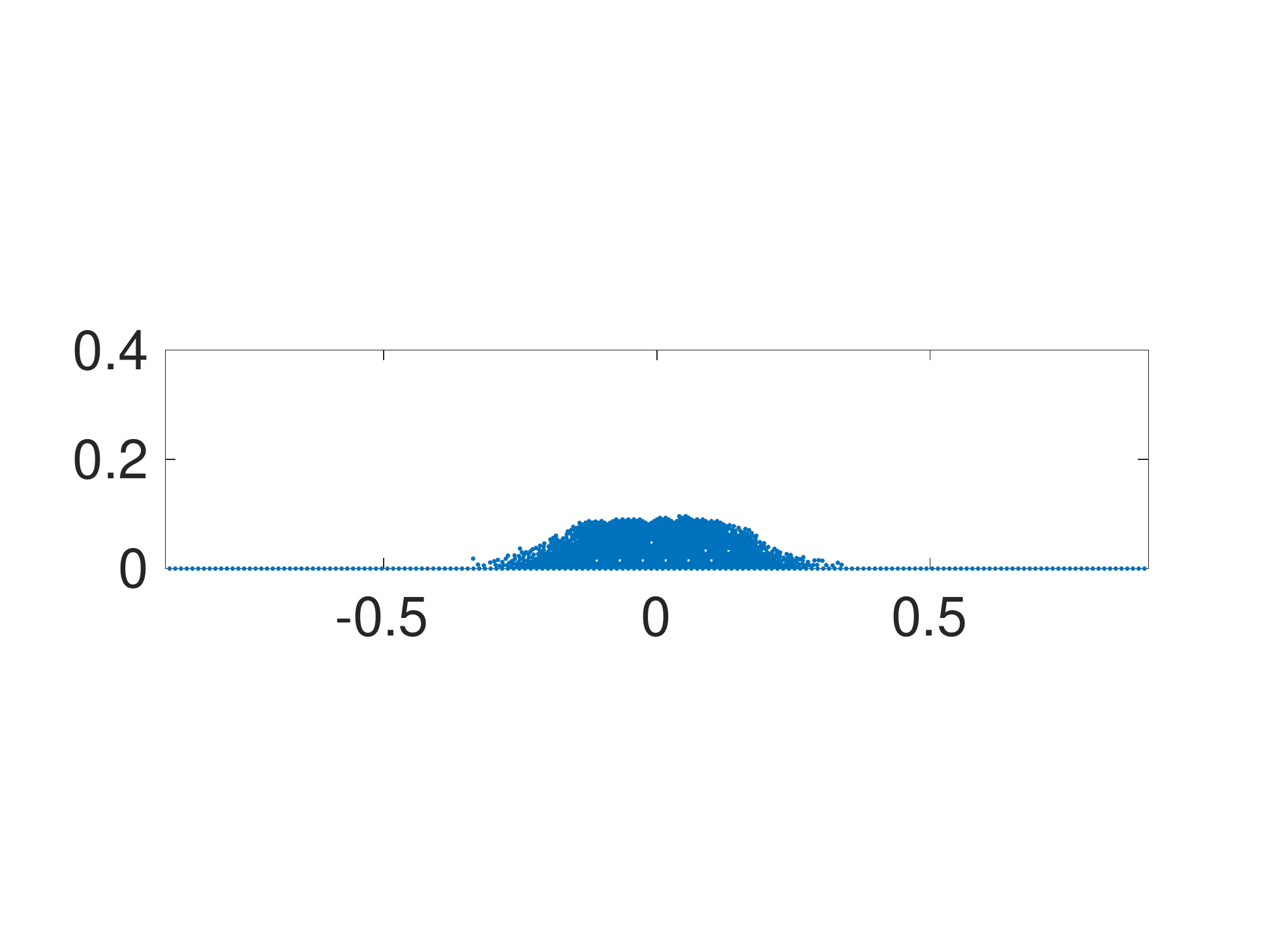} 
	\includegraphics[keepaspectratio=true, angle=0, width=0.32\textwidth]{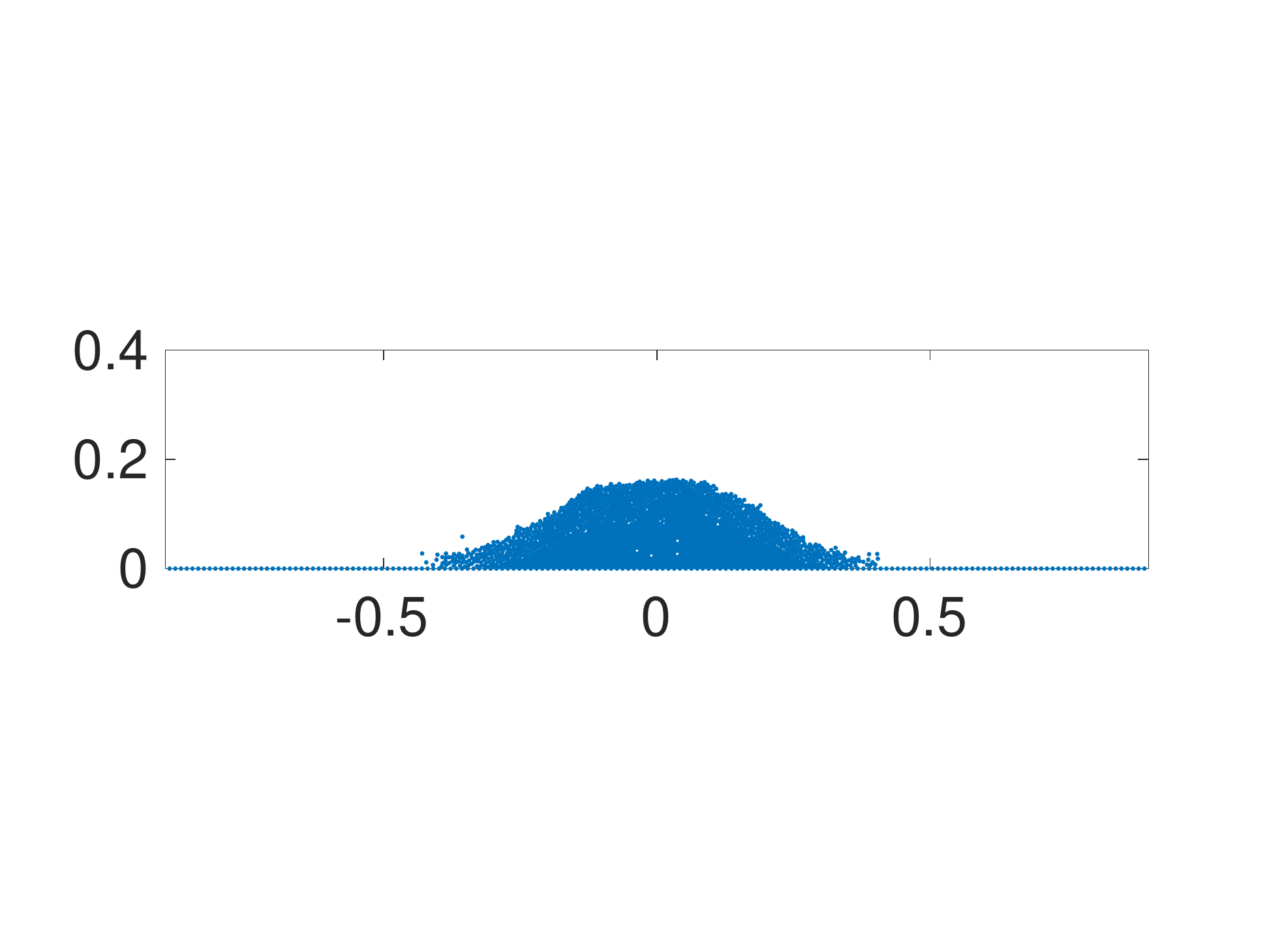} 
	\includegraphics[keepaspectratio=true, angle=0, width=0.32\textwidth]{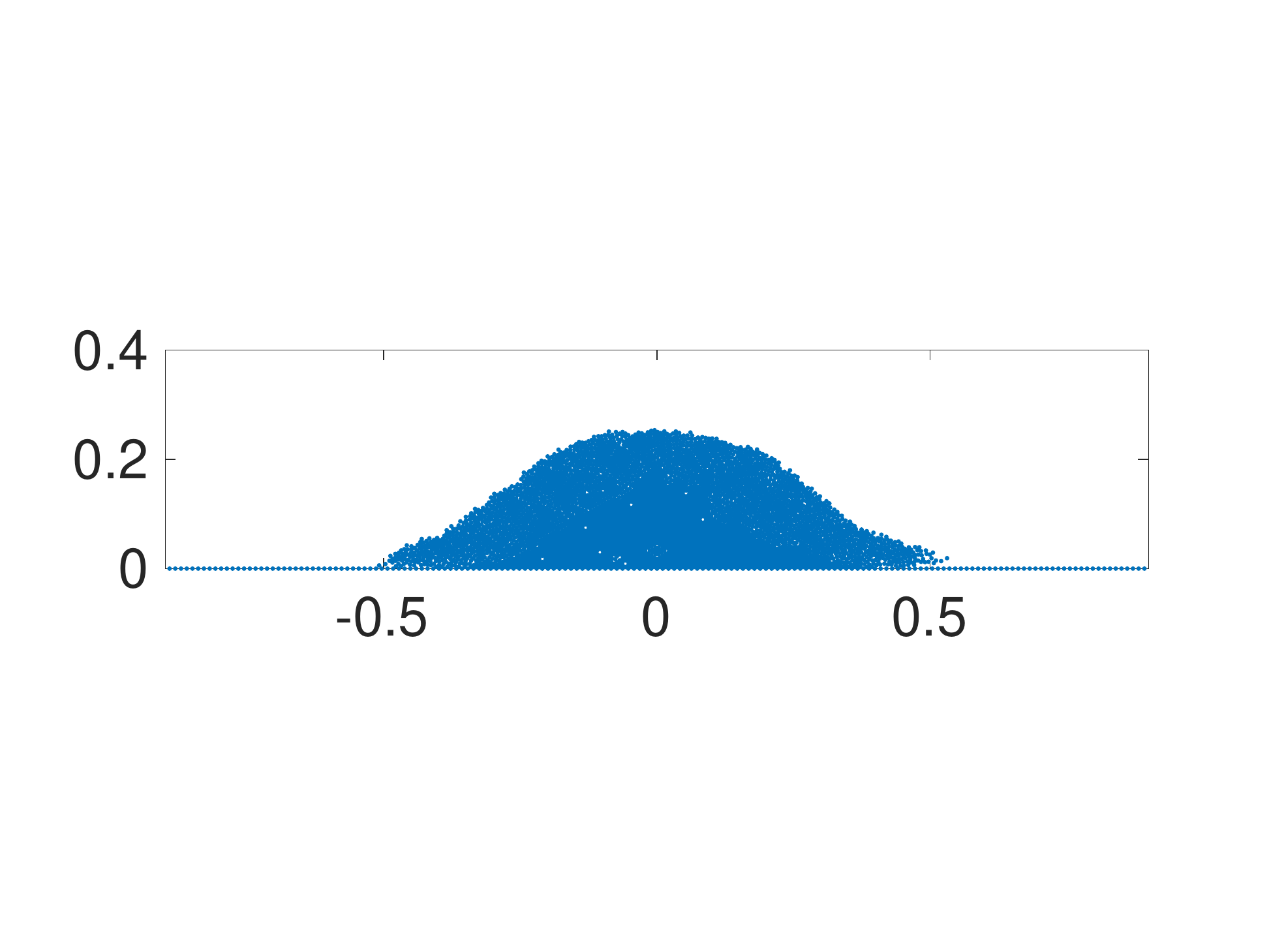} \\
	\vspace{-1cm}
	\includegraphics[keepaspectratio=true, angle=0, width=0.32\textwidth]{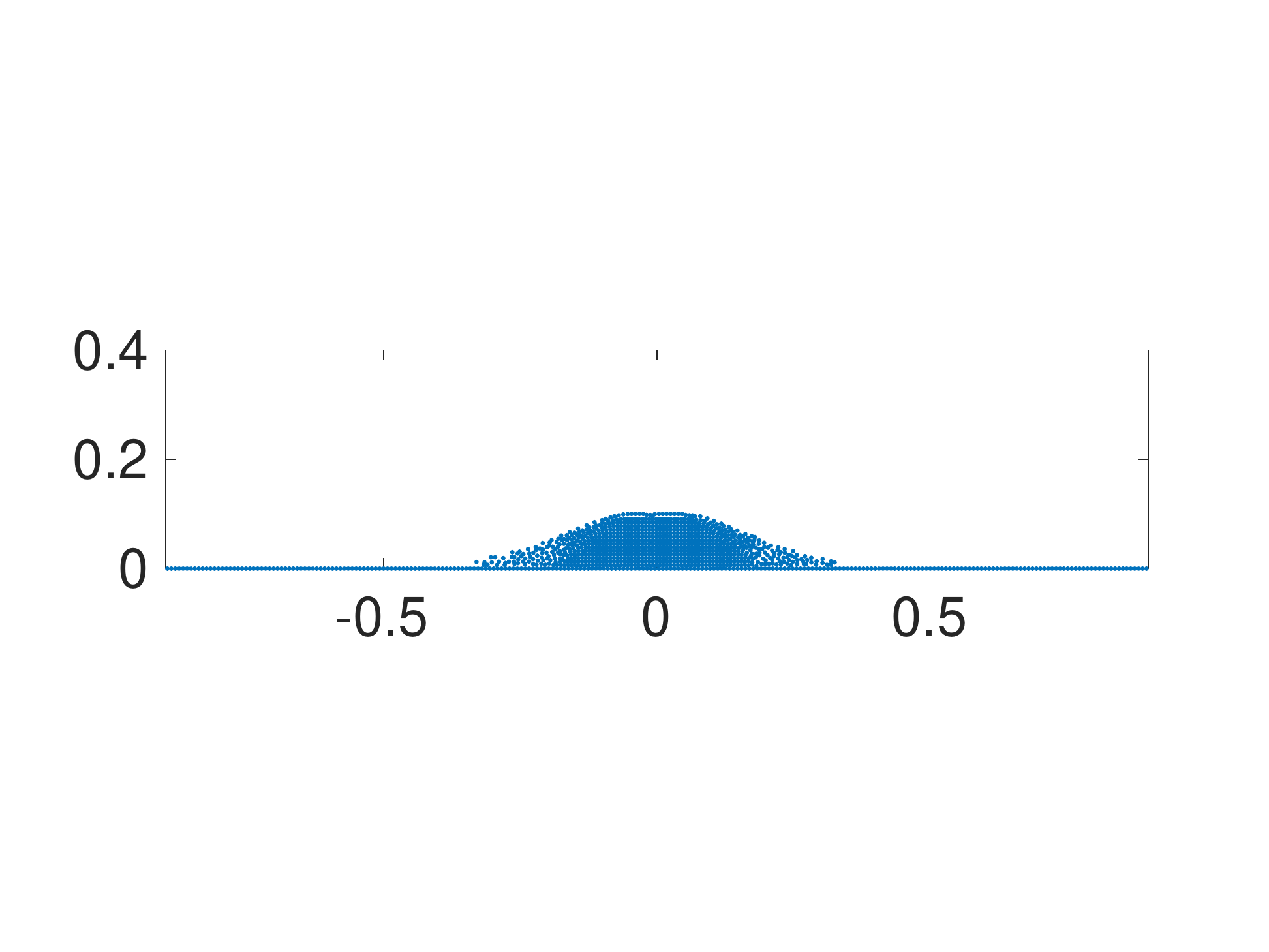} 
	\includegraphics[keepaspectratio=true, angle=0, width=0.32\textwidth]{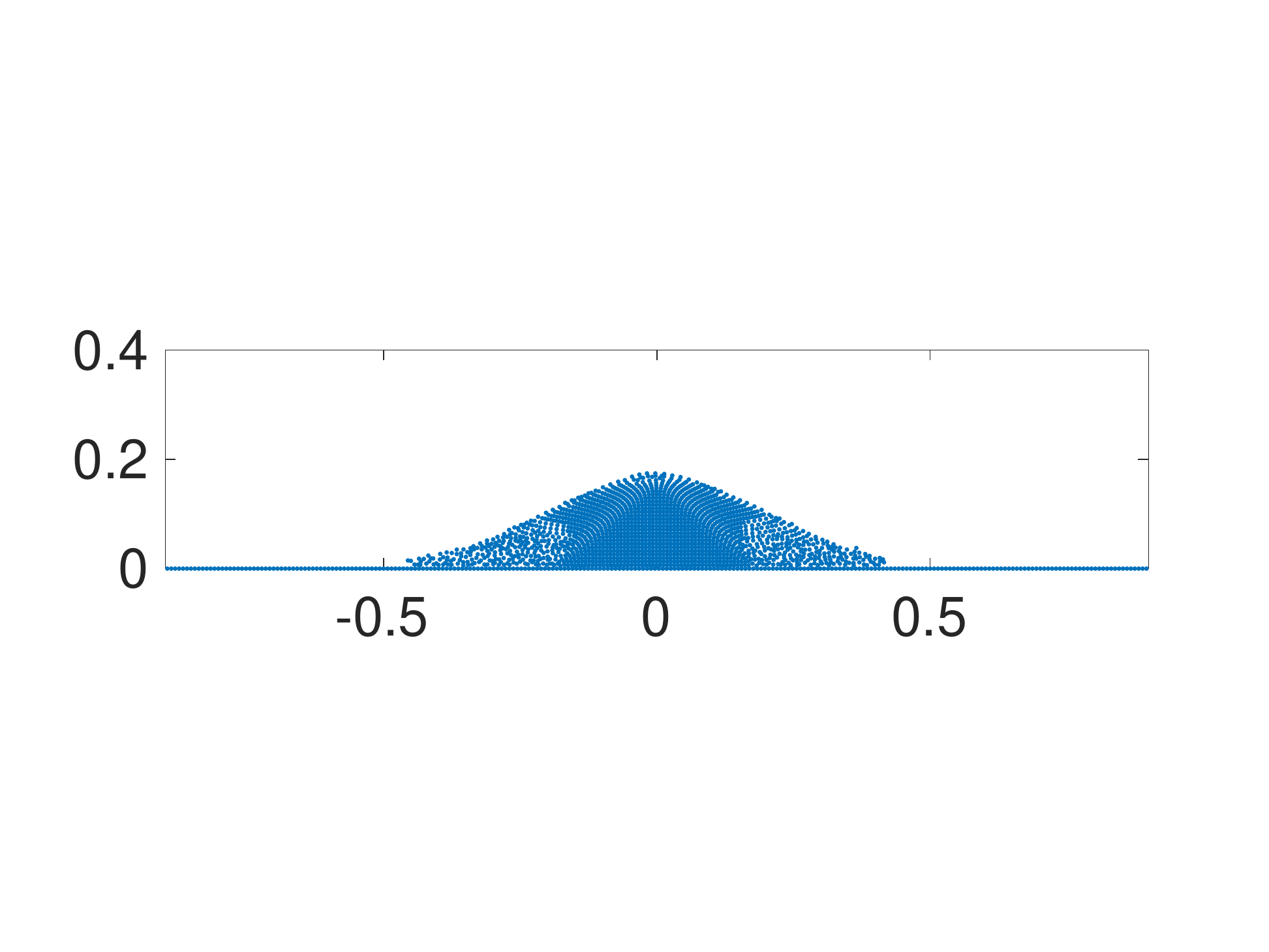} 
	\includegraphics[keepaspectratio=true, angle=0, width=0.32\textwidth]{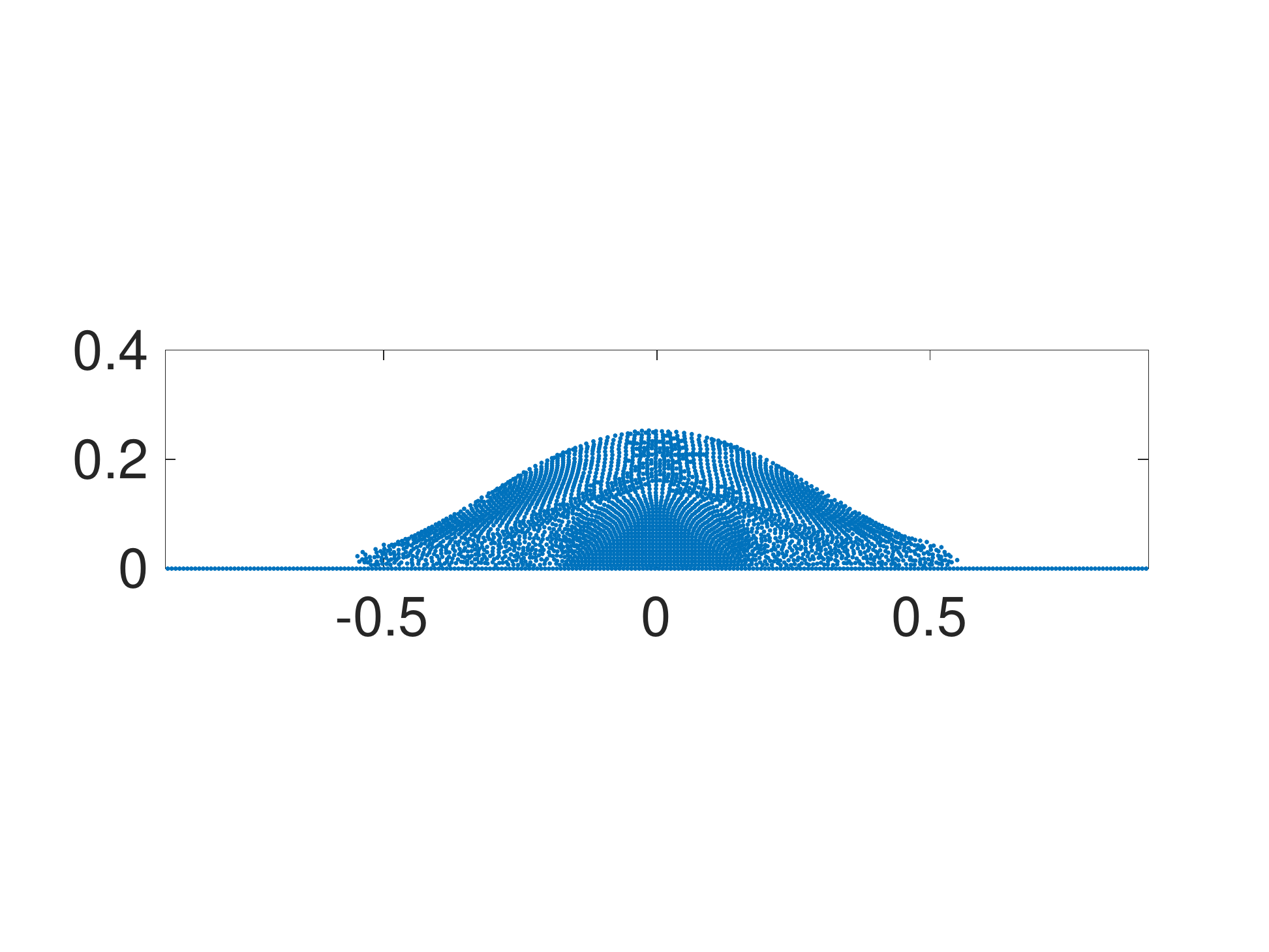} 	\\
	\vspace{-1cm}
	\includegraphics[keepaspectratio=true, angle=0, width=0.32\textwidth]{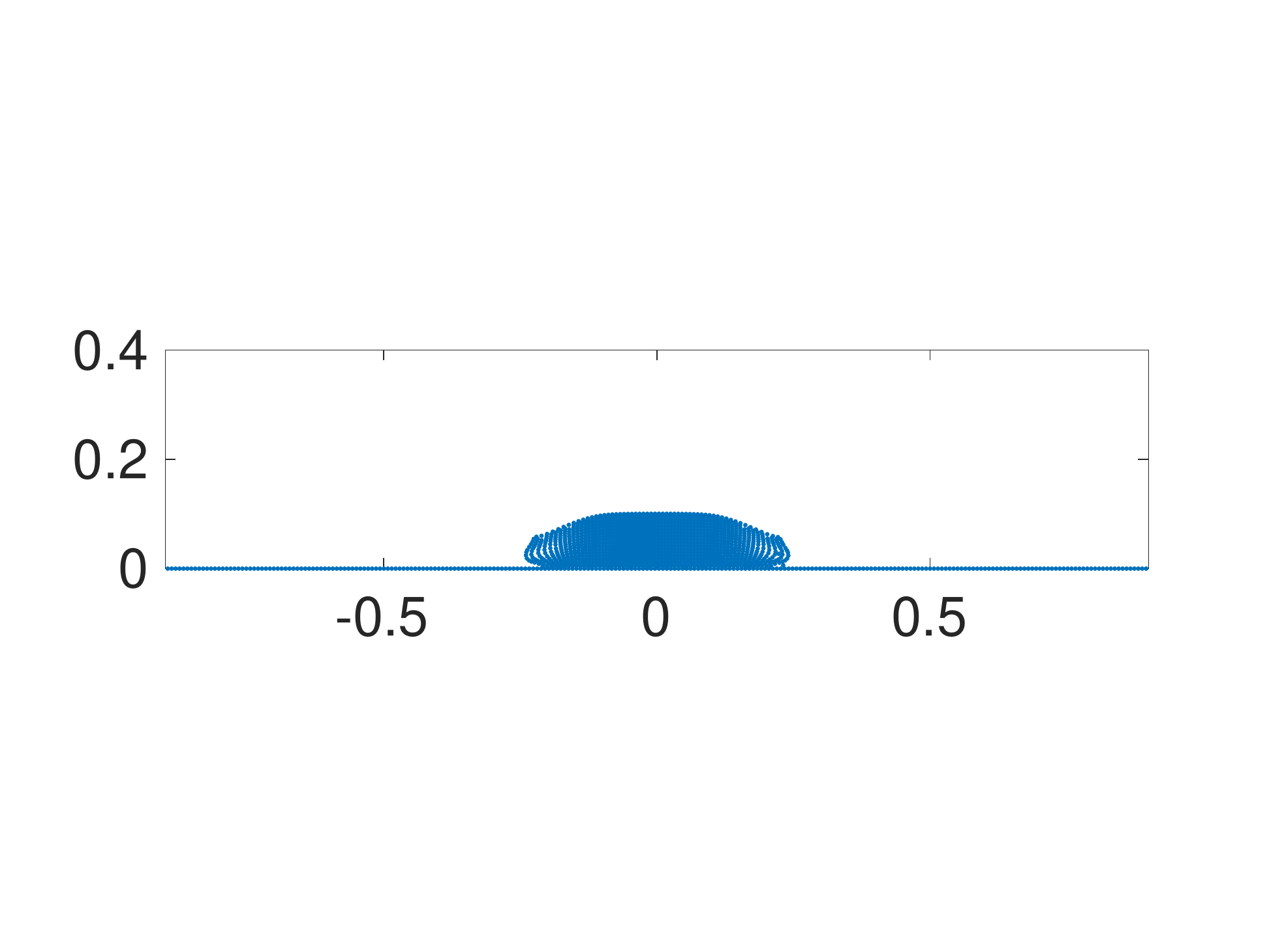} 
	\includegraphics[keepaspectratio=true, angle=0, width=0.32\textwidth]{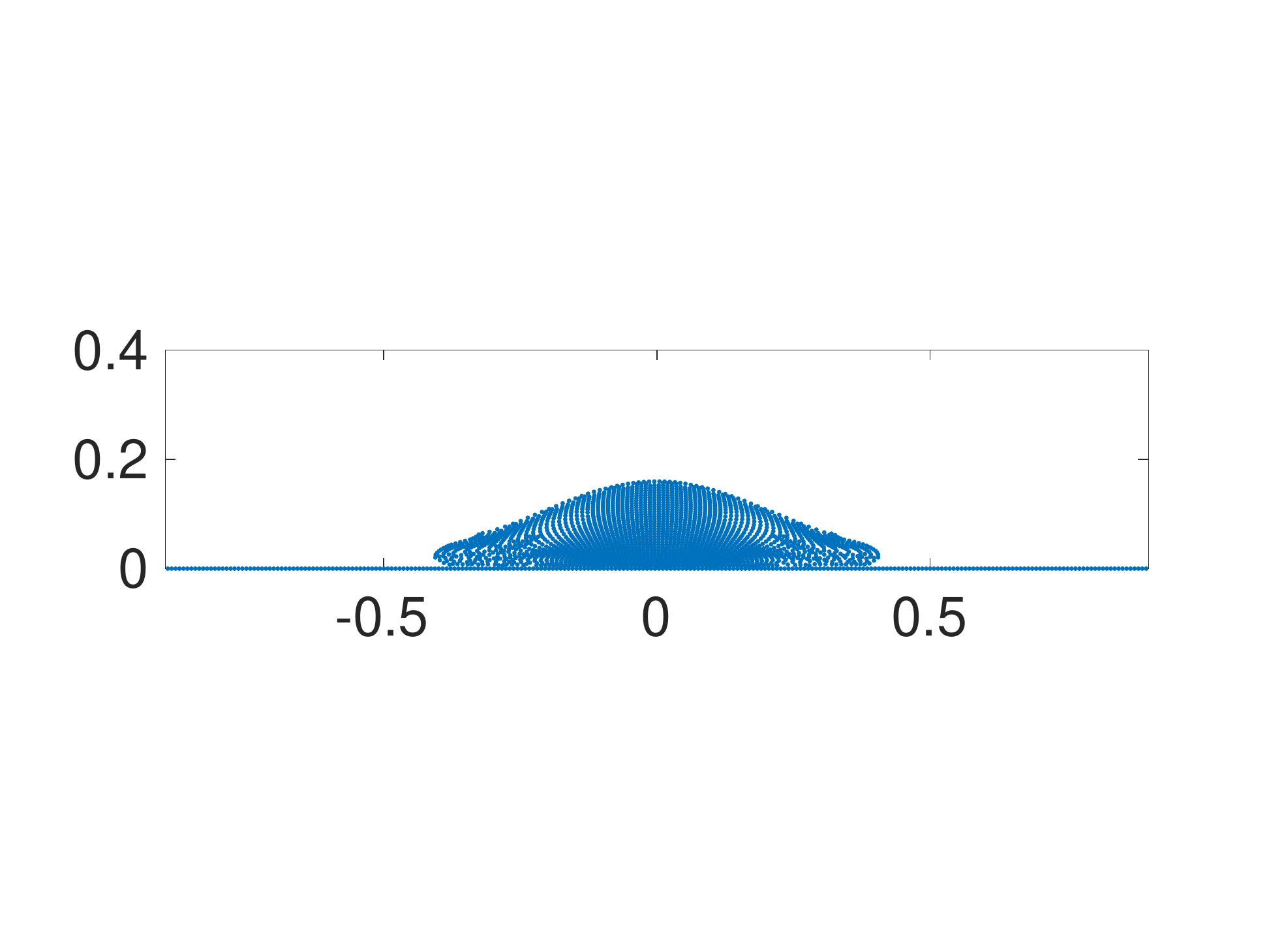} 
	\includegraphics[keepaspectratio=true, angle=0, width=0.32\textwidth]{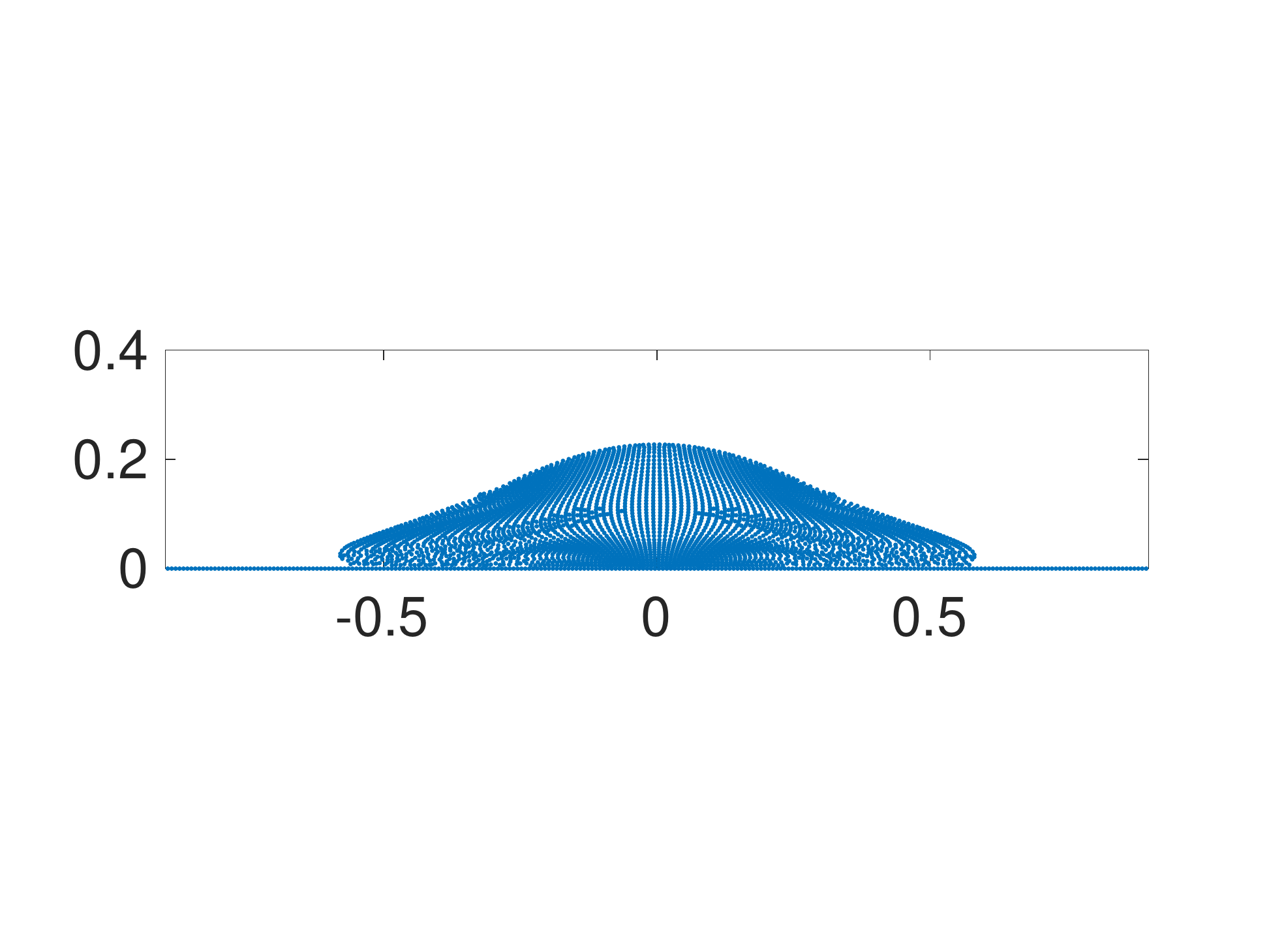} 	
		\caption{Time development at $ t = 0.3$ for $a = 0.5$ (left),  $a = 1$ ( middle), $a = 2$ (right). First row: micoscopic model (Section \ref{Lac}), second row: plasticity model (Section \ref{macroeq}),third row: Coloumb constitutive model (Section \ref{Jop})}	
	\label{column_compare}
	\centering
\end{figure}

\begin{figure}
	\centering
		\includegraphics[keepaspectratio=true, angle=0, width=0.32\textwidth]{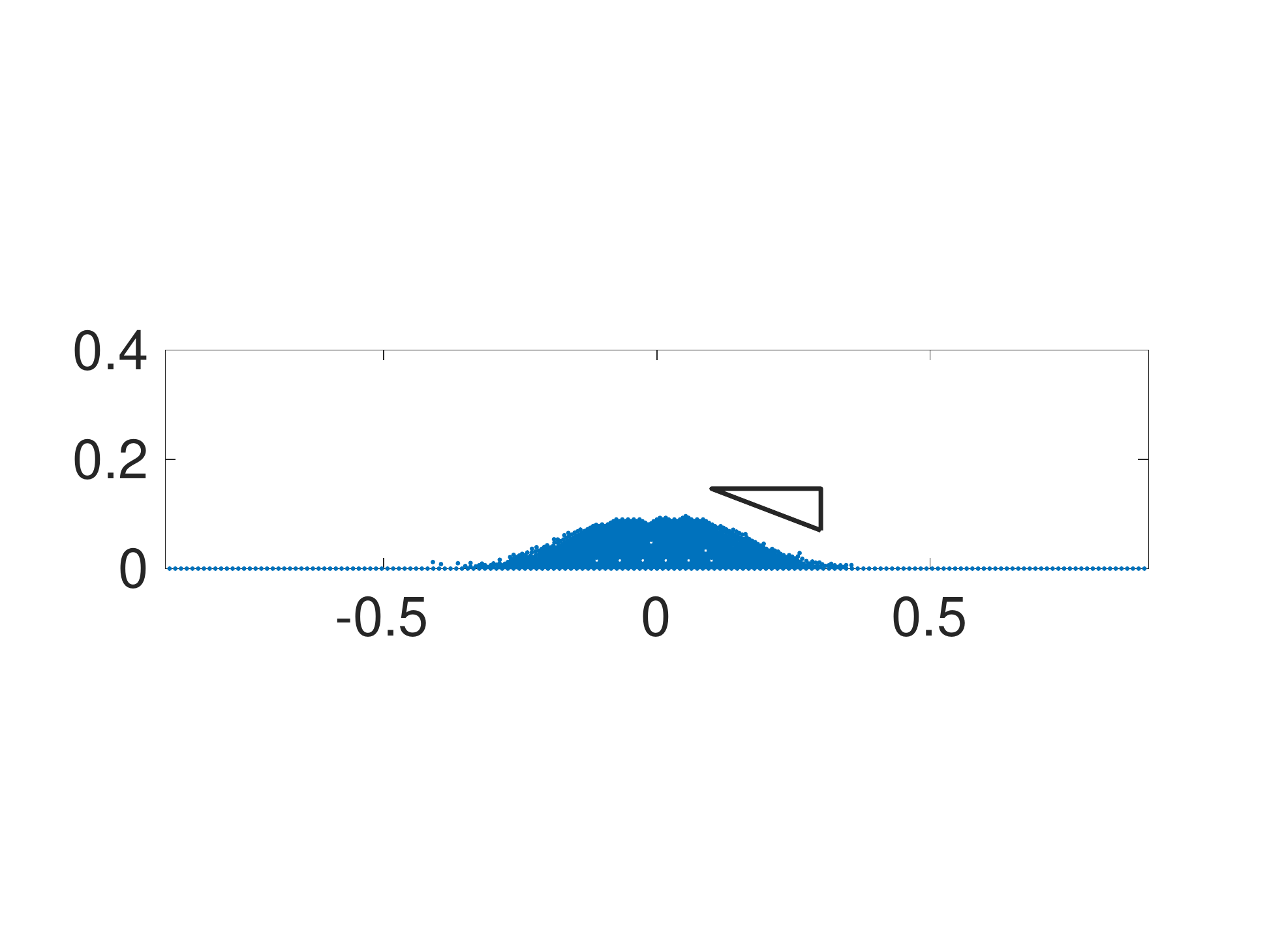} 
	\includegraphics[keepaspectratio=true, angle=0, width=0.32\textwidth]{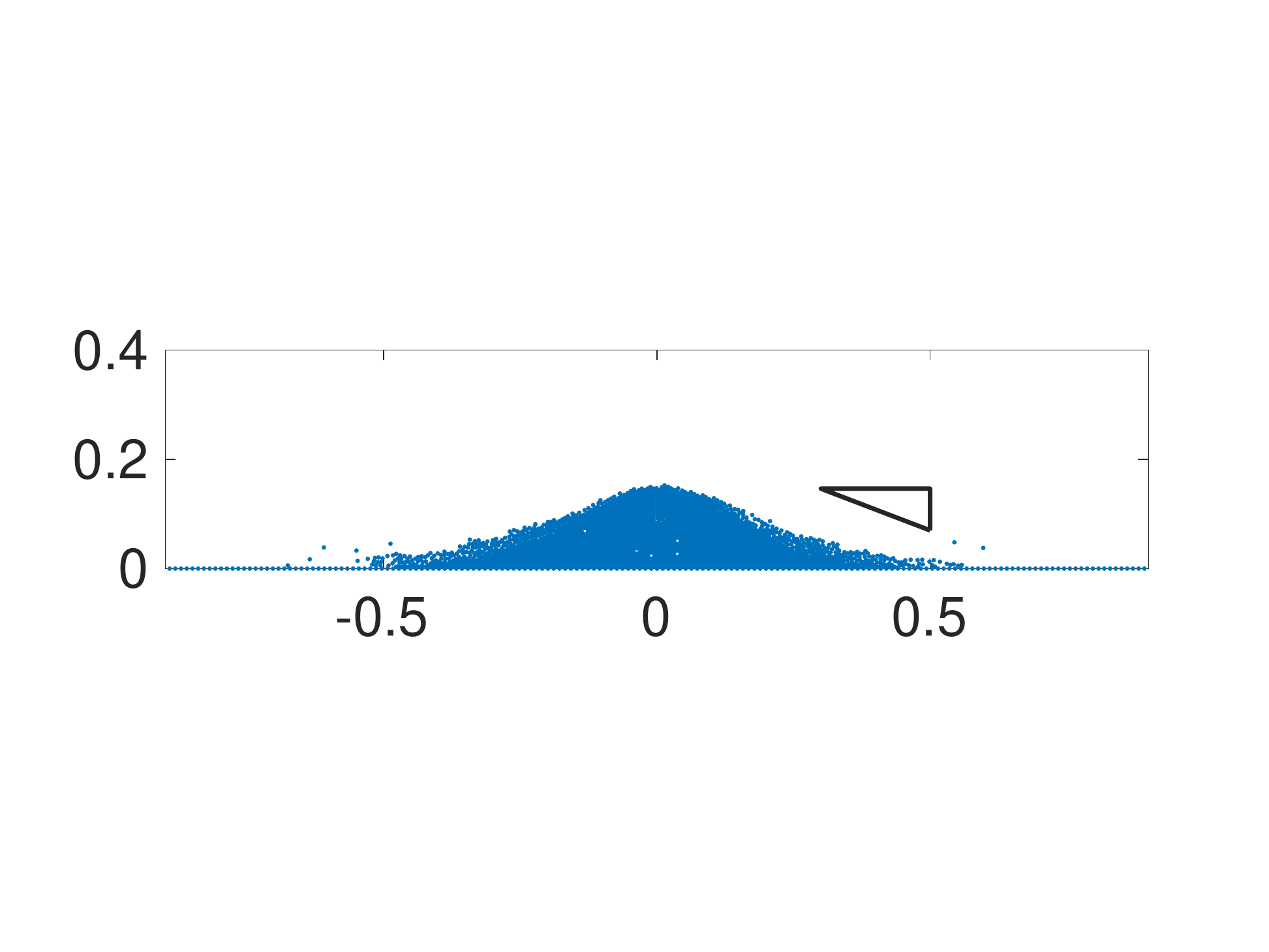} 
	\includegraphics[keepaspectratio=true, angle=0, width=0.32\textwidth]{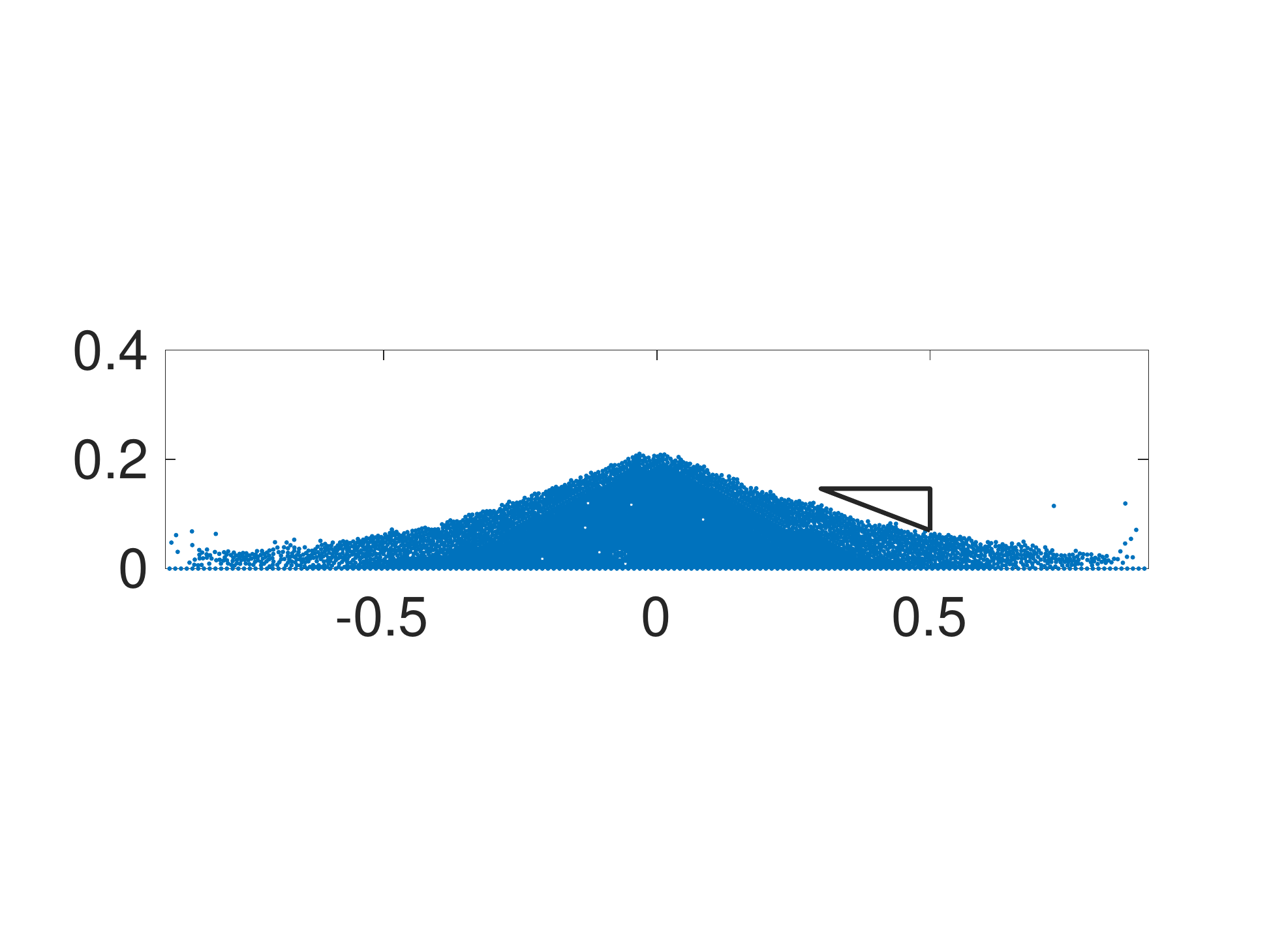} \\
	\vspace{-1cm}
	\includegraphics[keepaspectratio=true, angle=0, width=0.32\textwidth]{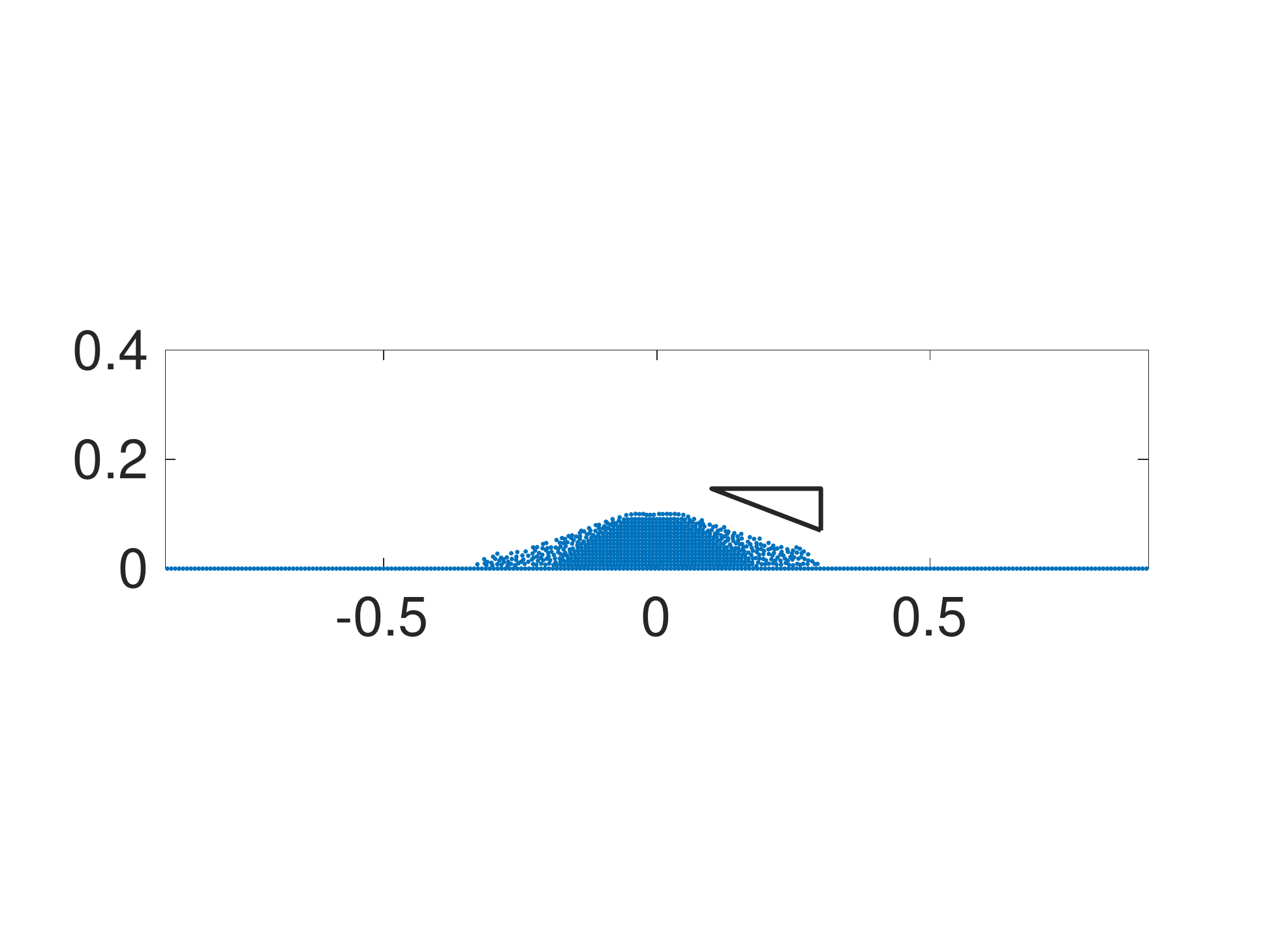} 
	\includegraphics[keepaspectratio=true, angle=0, width=0.32\textwidth]{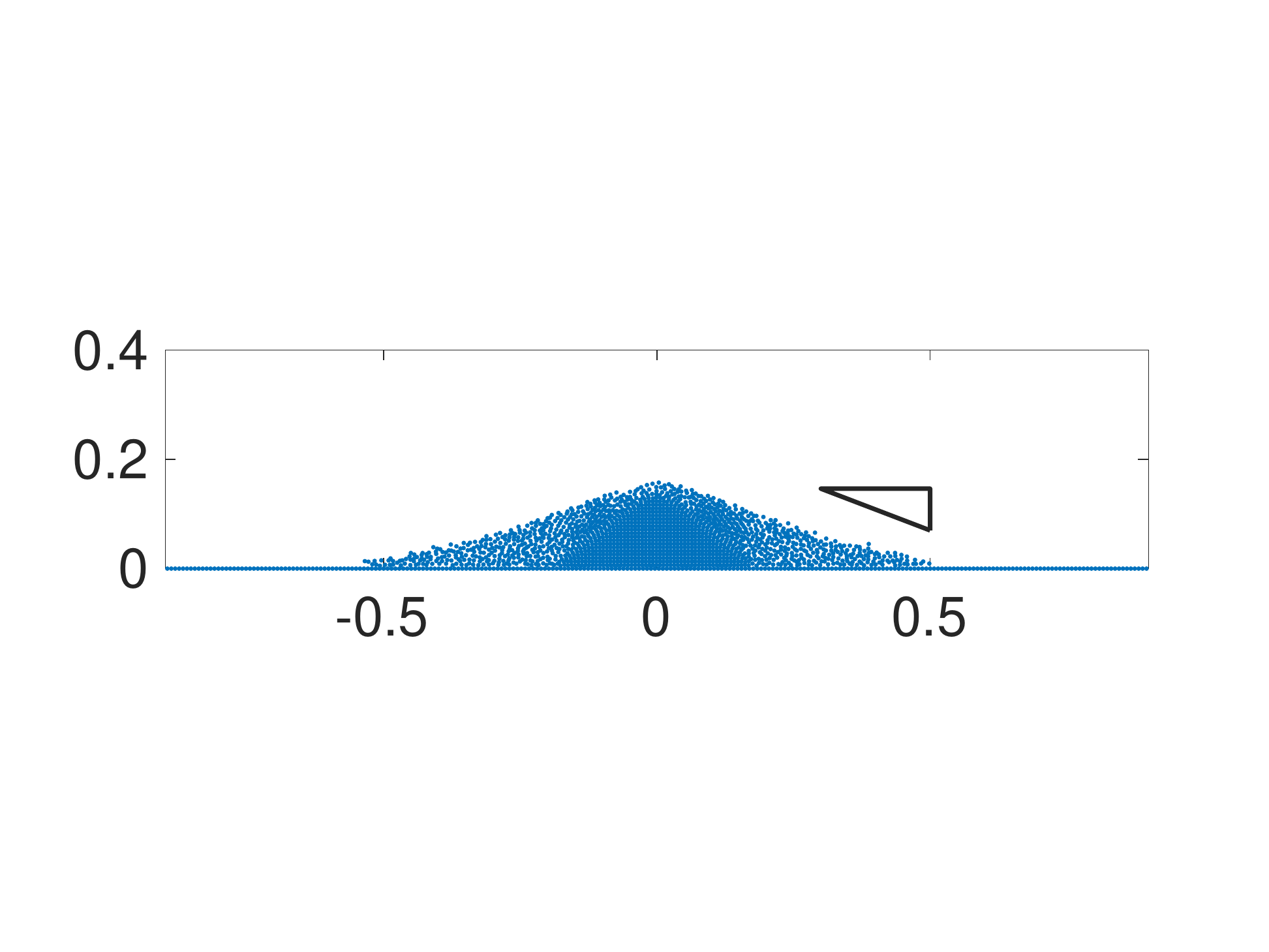} 
	\includegraphics[keepaspectratio=true, angle=0, width=0.32\textwidth]{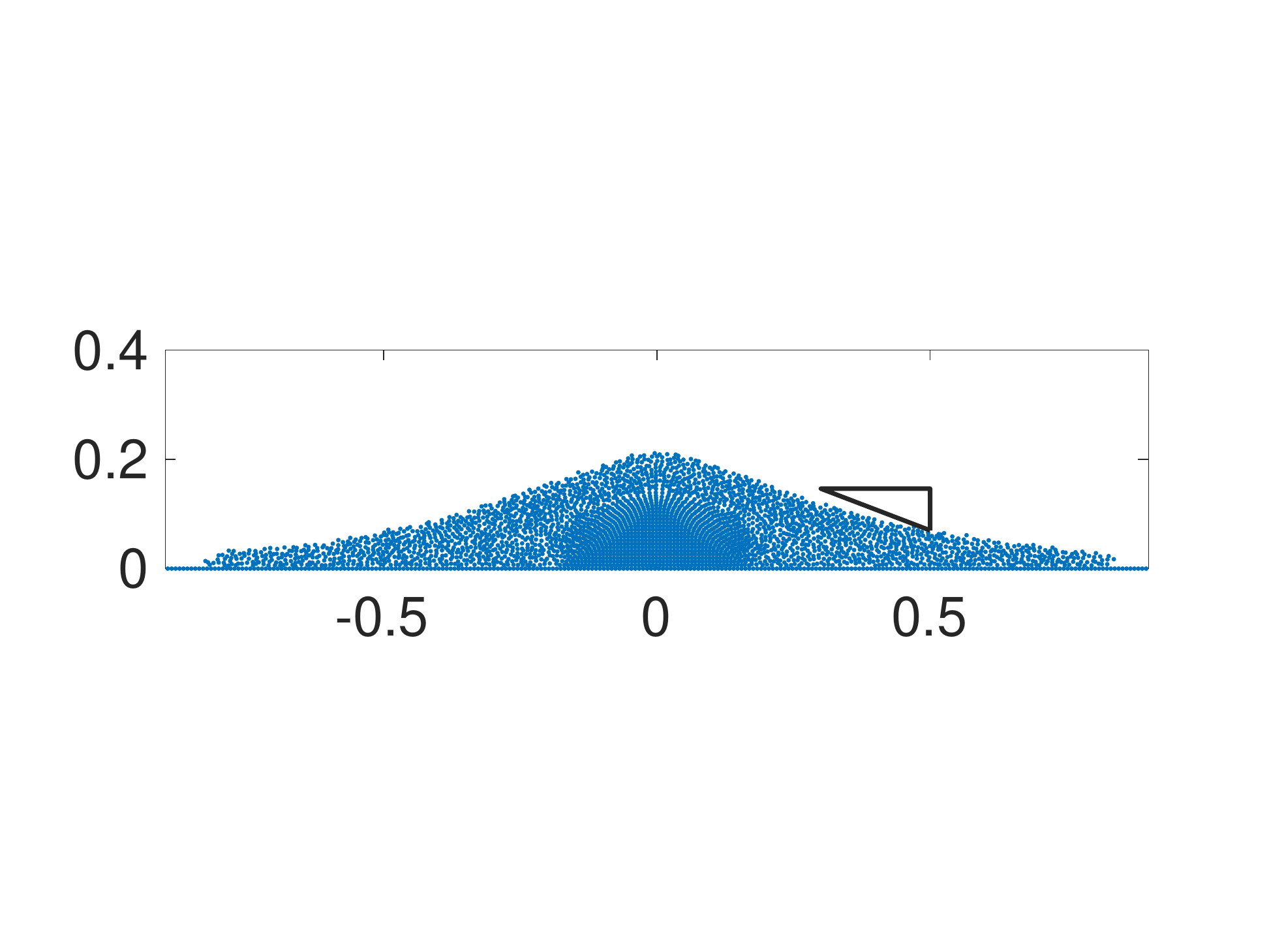} 	\\
	\vspace{-1cm}
	\includegraphics[keepaspectratio=true, angle=0, width=0.32\textwidth]{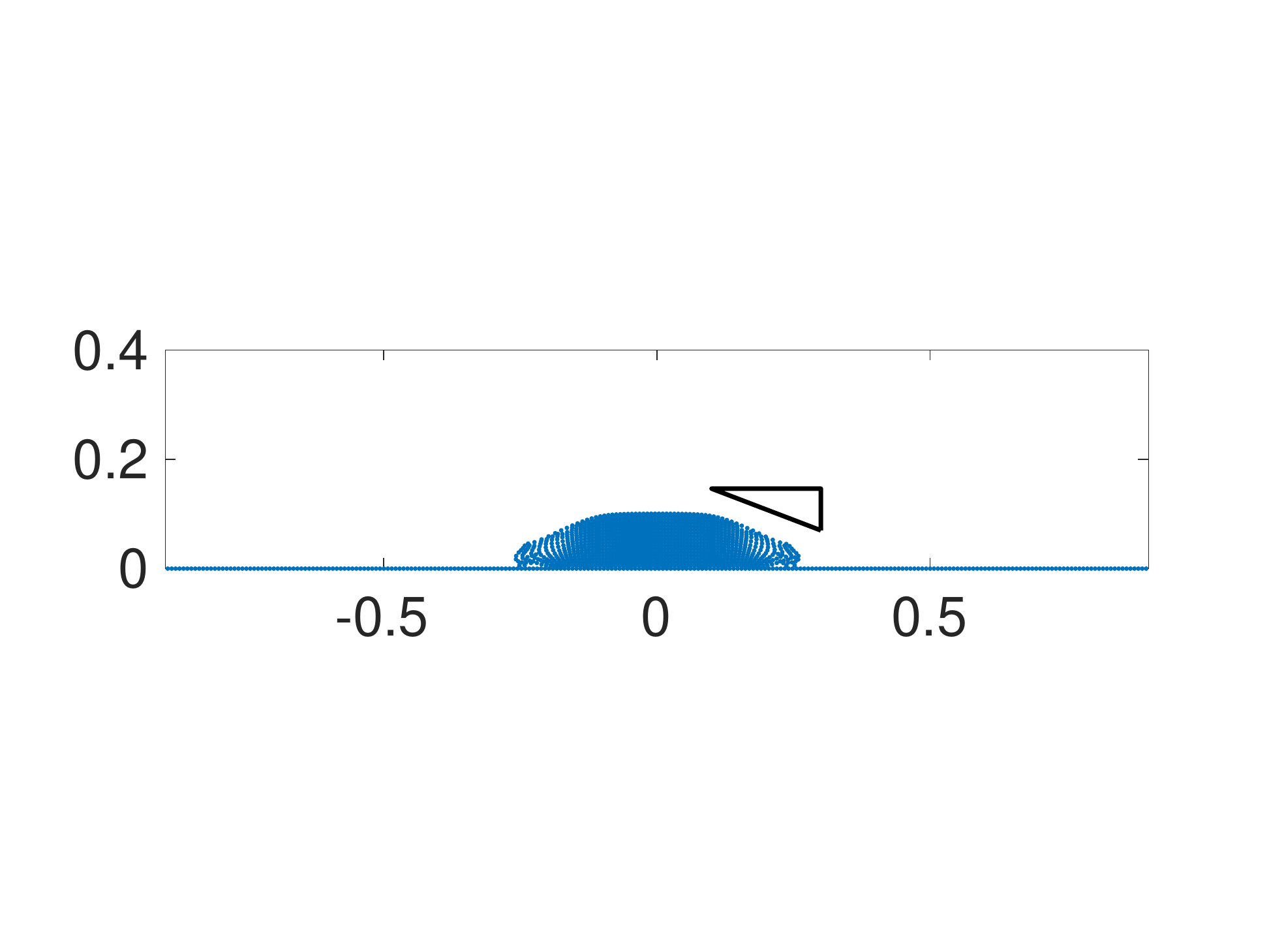} 
	\includegraphics[keepaspectratio=true, angle=0, width=0.32\textwidth]{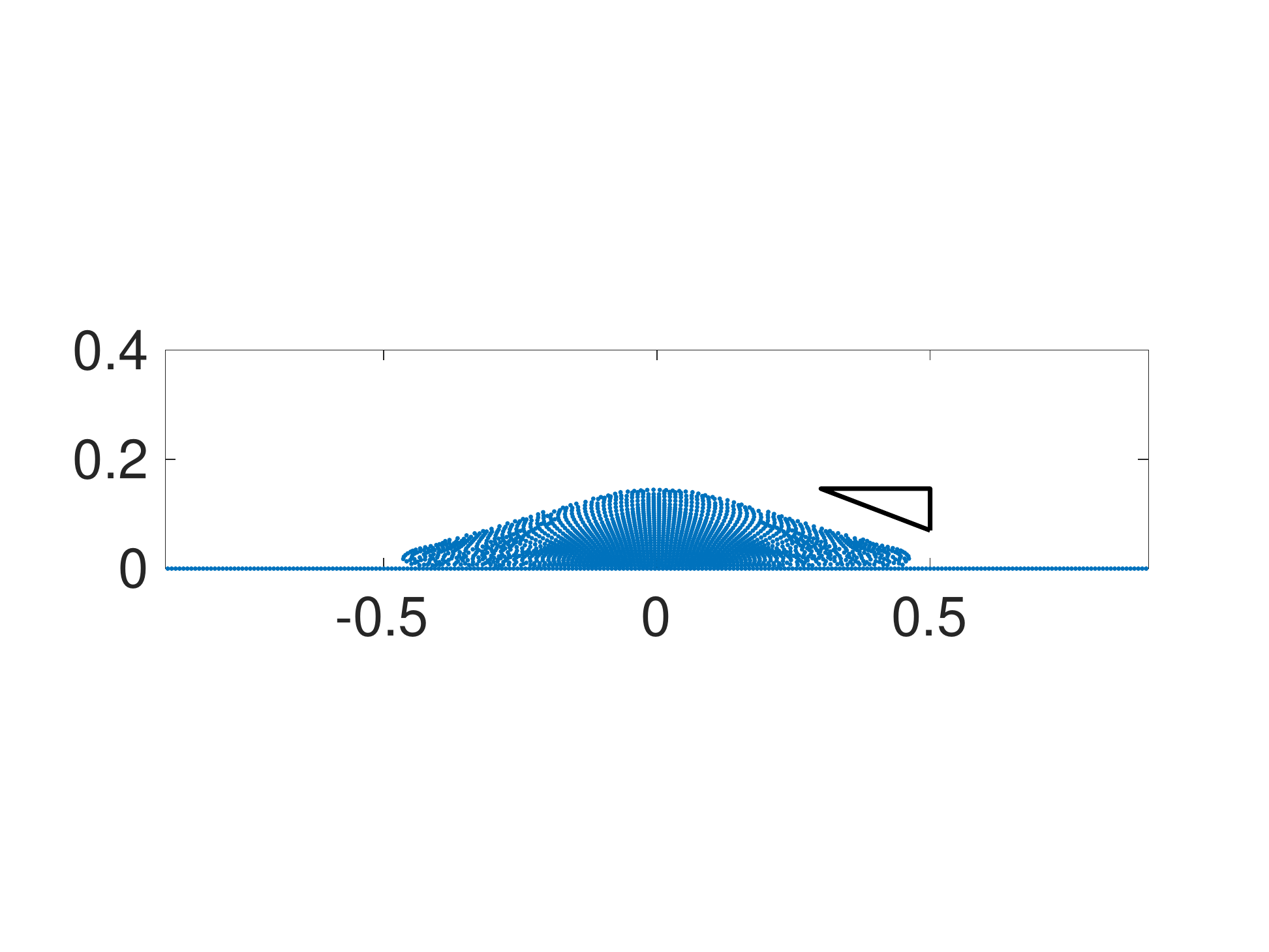} 
	\includegraphics[keepaspectratio=true, angle=0, width=0.32\textwidth]{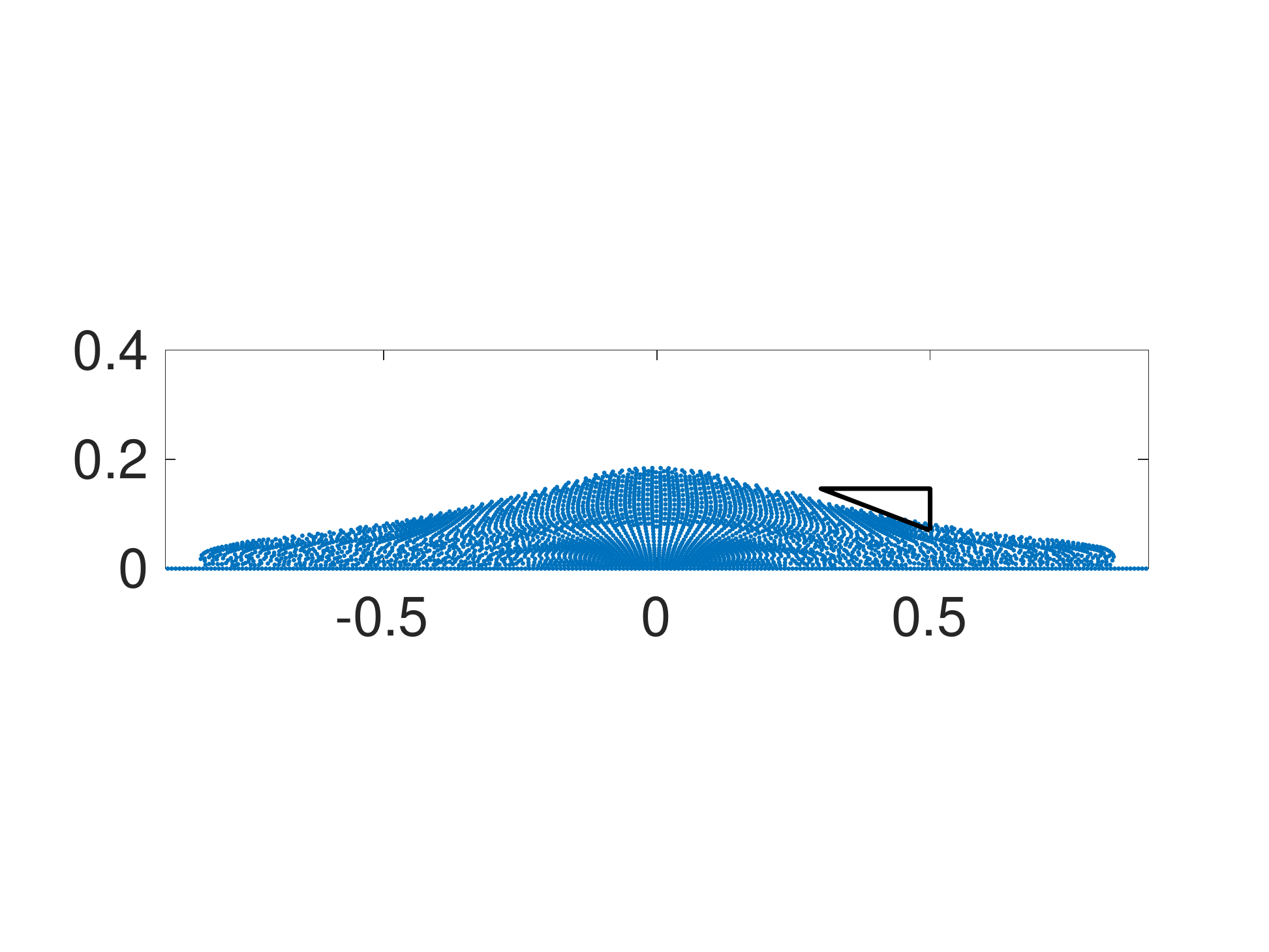} 	
		\caption{Time development at $ t = 0.5$ for $a = 0.5$ (left),  $a = 1$ ( middle), $a = 2$ (right). First row: micoscopic model (Section \ref{Lac}), second row: plasticity model (Section \ref{macroeq}),third row: Coloumb constitutive model (Section \ref{Jop})}	
	\label{column_compare_t0dot5}
	\centering
\end{figure}

For the first numerical experiment we fix all  parameters to the values given above. 
We determine first the results of the  microscopic model (Section \ref{Lac}) with $N = 2000$. 
Then, the results of a well resolved simulation using the  plasticity model  (Section \ref{macroeq}) and  the simplified model (Section \ref{Jop})
with  $N = 2000$ grid particles are computed.
In Figures \ref{column_t0} -- \ref{column_compare_t0dot5}  we have plotted the results  obtained from these three  models at times $ t = 0.0, 0.2, 0.3,  0.5$.  One observes a good coincidence between microscopic and plasticity model with some discrepancies for the  initial stages of the evolution. The simplified Coulomb constitutive model deviates further from the microscopic model.


\begin{figure}
	\centering
	\includegraphics[keepaspectratio=true, angle=0, width=0.8\textwidth]{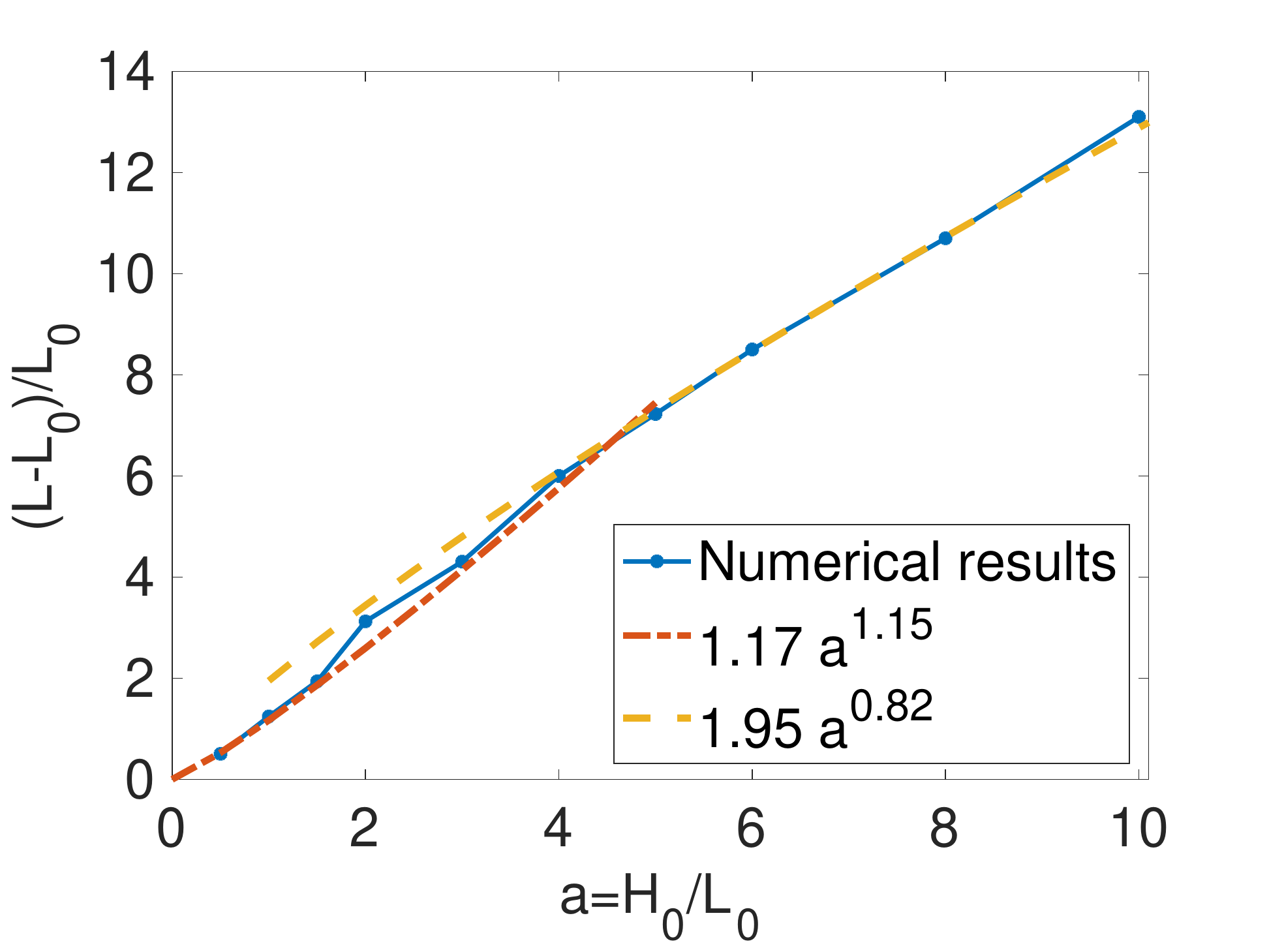}
	\caption{Non-dimensional run-out distance vs non-dimensional aspect ratio $a=H_0/L_0$.} 	
	\label{runout}
	\centering
\end{figure}

In Fig. \ref{runout} we have plotted the non-dimensional run out distance against the aspect ratio $a$. They are 
similar for all models.  We observed that for smaller values of aspect ratio $a$ a good approximation is given by 
$
\frac{L-L_0}{L_0} =
1.17 a^{1.15}.
$
For larger values of $a$ a good approximation is given by 
$
1.95^{0.82}.
$
Such values are consistent with  numerical and experimental results  discussed, for example,  in \cite{DK15}.

\begin{figure}
	\centering

	\includegraphics[keepaspectratio=true, angle=0, width=0.32\textwidth]{column_a_1_kamrin_t_0dot2} 
		\includegraphics[keepaspectratio=true, angle=0, width=0.32\textwidth]{column_a_1_kamrin_t_0dot3} 
		\includegraphics[keepaspectratio=true, angle=0, width=0.32\textwidth]{column_a_1_kamrin_t_0dot5} \\
\vspace{-1cm}		
	\includegraphics[keepaspectratio=true, angle=0, width=0.32\textwidth]{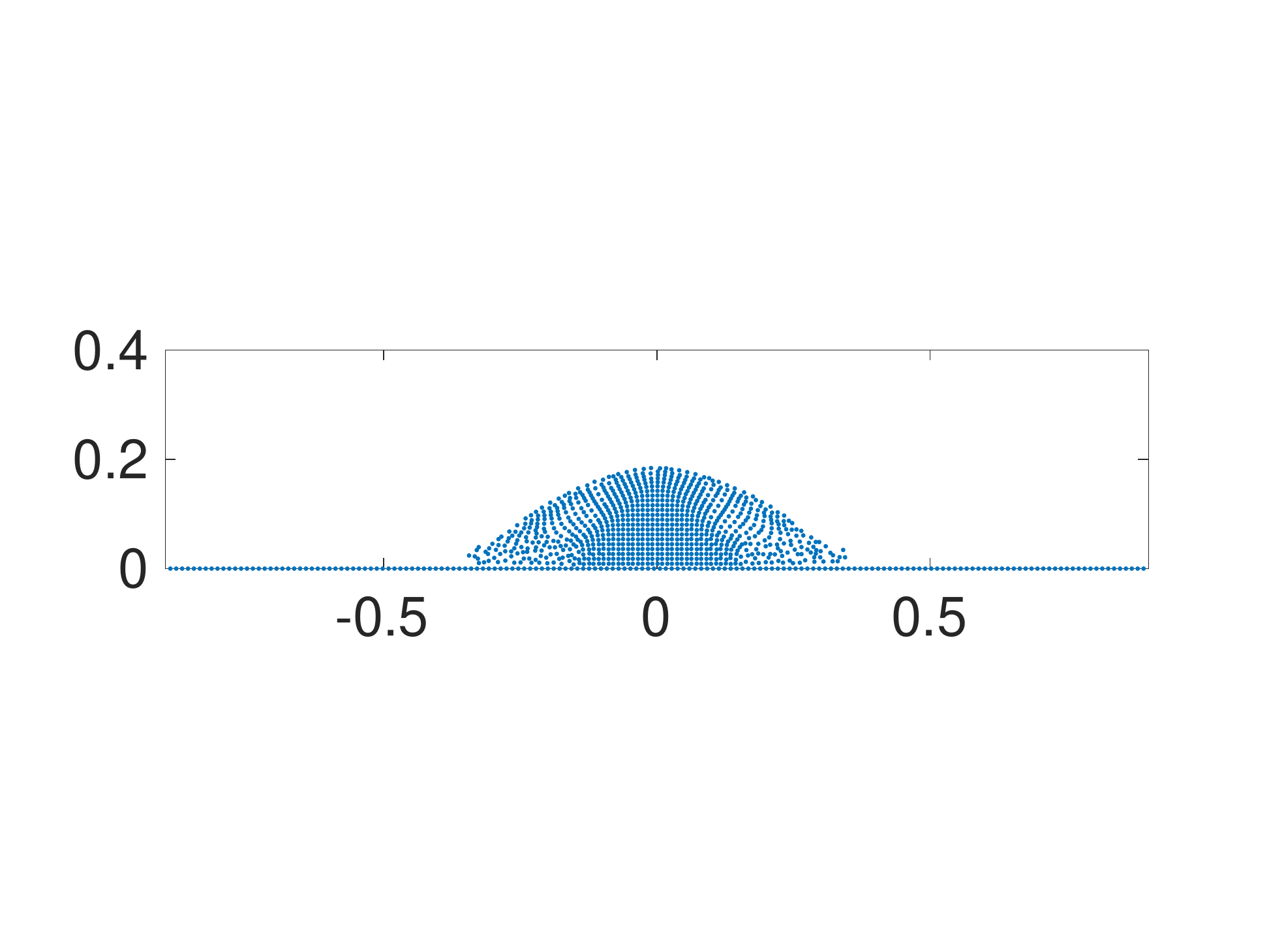} 
		\includegraphics[keepaspectratio=true, angle=0, width=0.32\textwidth]{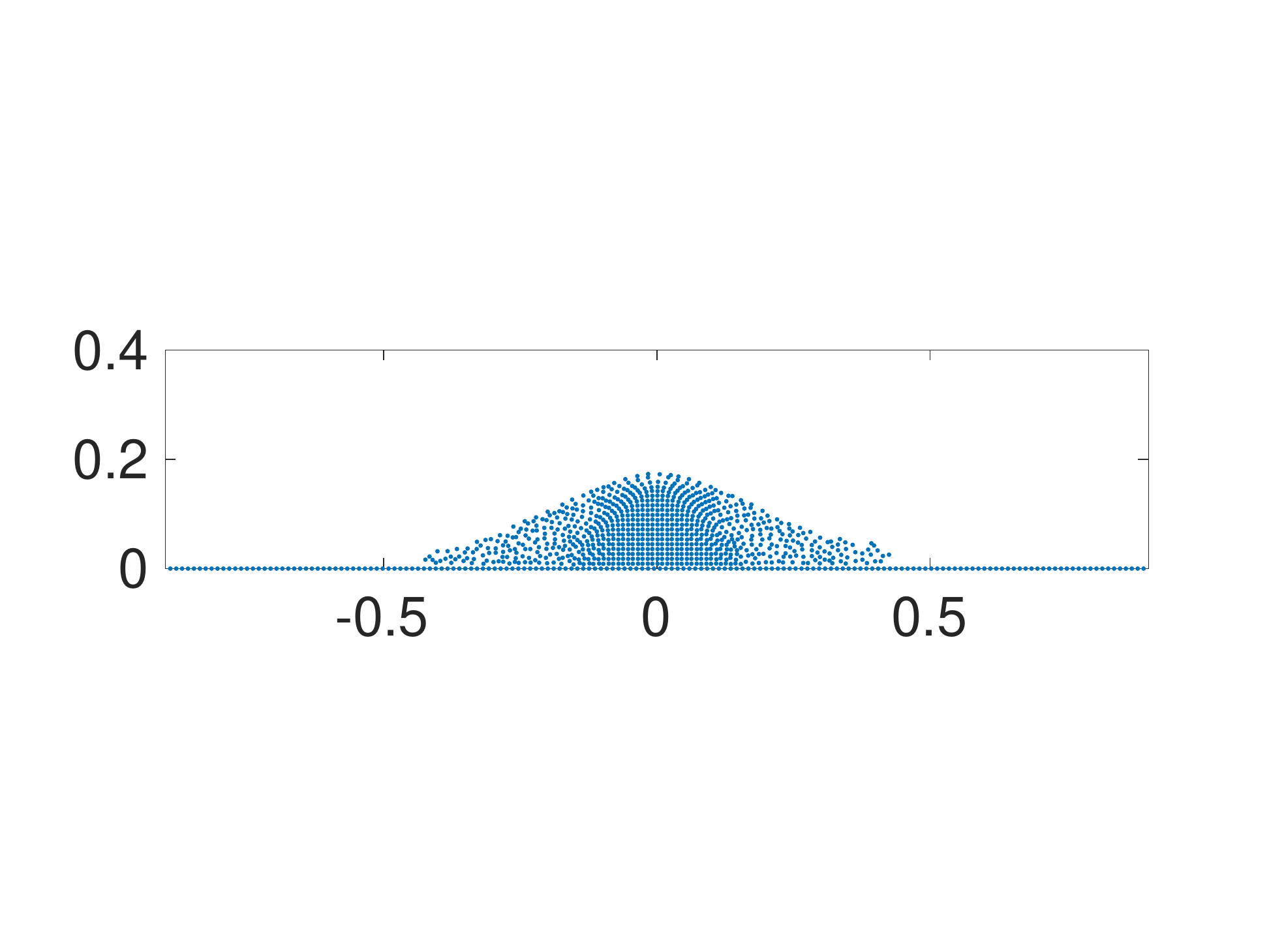} 
		\includegraphics[keepaspectratio=true, angle=0, width=0.32\textwidth]{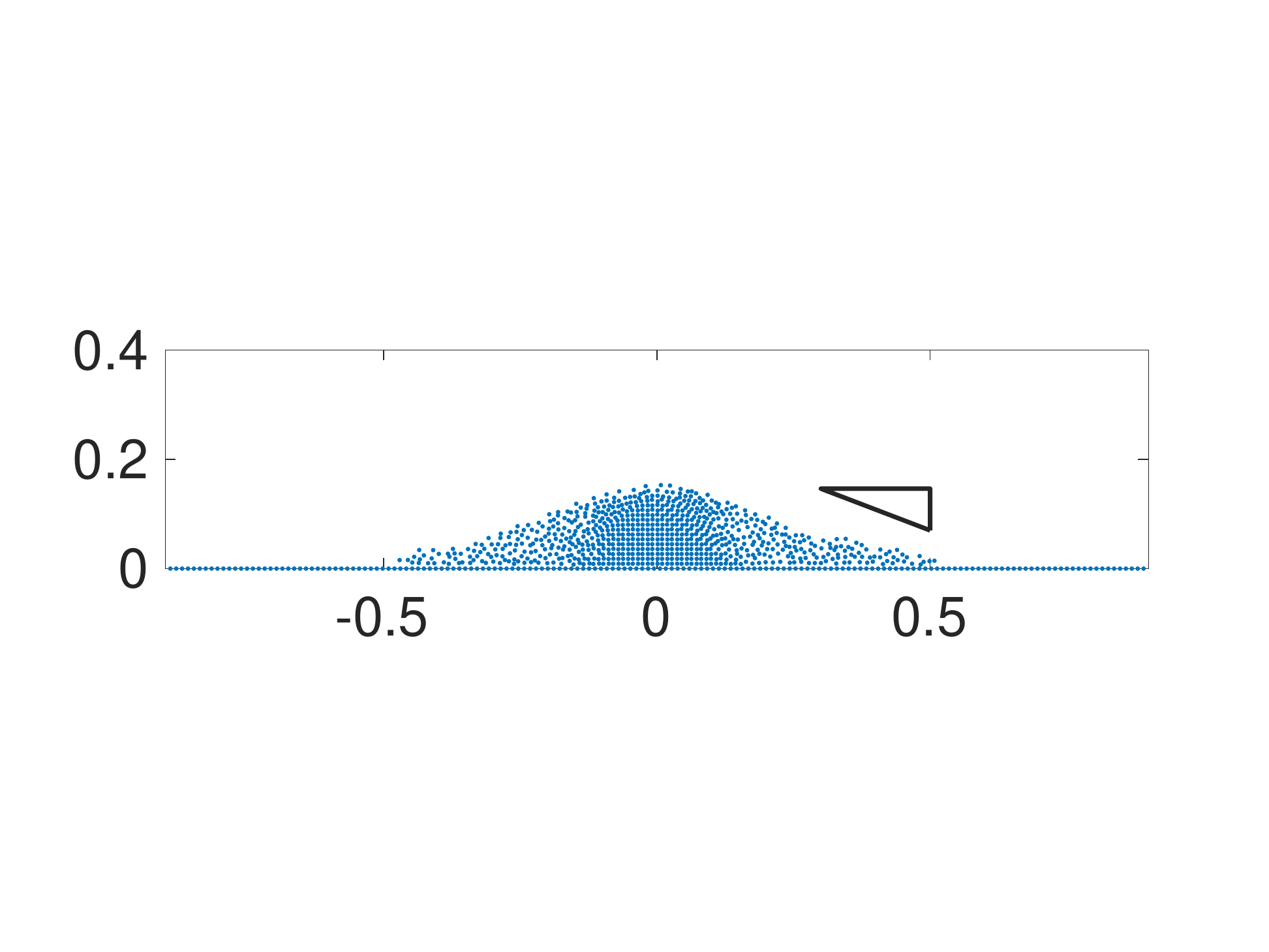} \\
\vspace{-1cm}	
		\includegraphics[keepaspectratio=true, angle=0, width=0.32\textwidth]{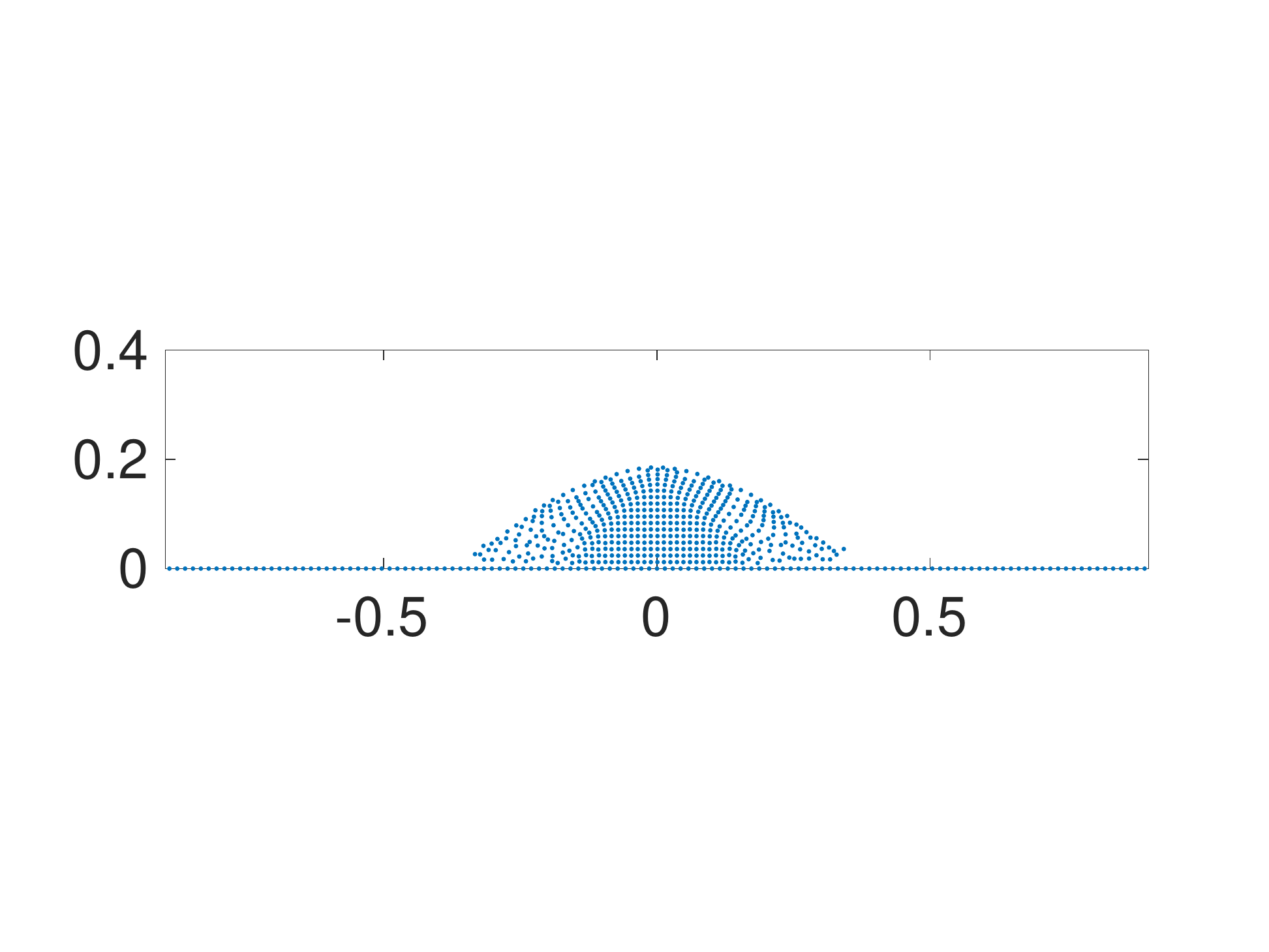}
	\includegraphics[keepaspectratio=true, angle=0, width=0.32\textwidth]{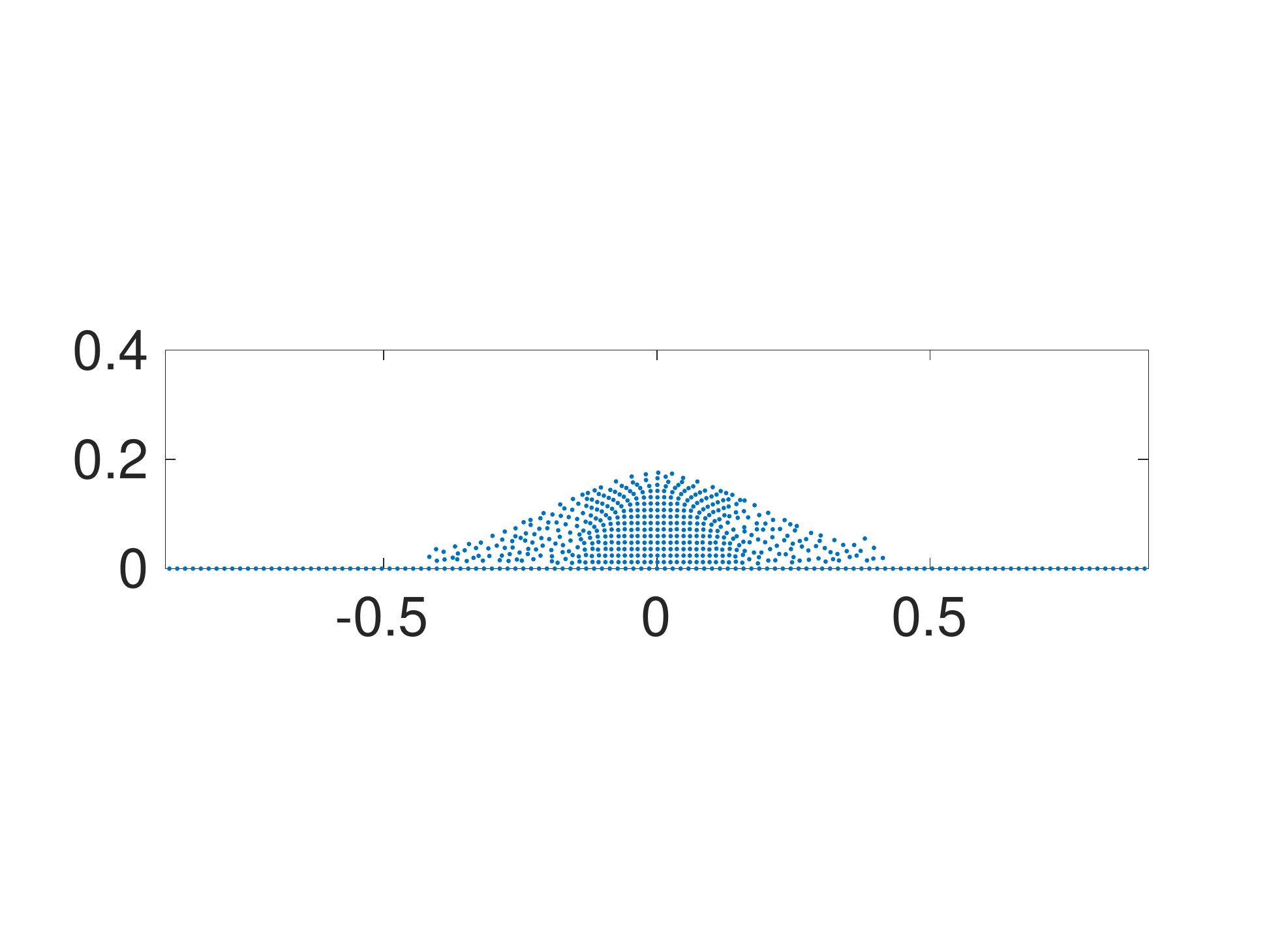}
	\includegraphics[keepaspectratio=true, angle=0, width=0.32\textwidth]{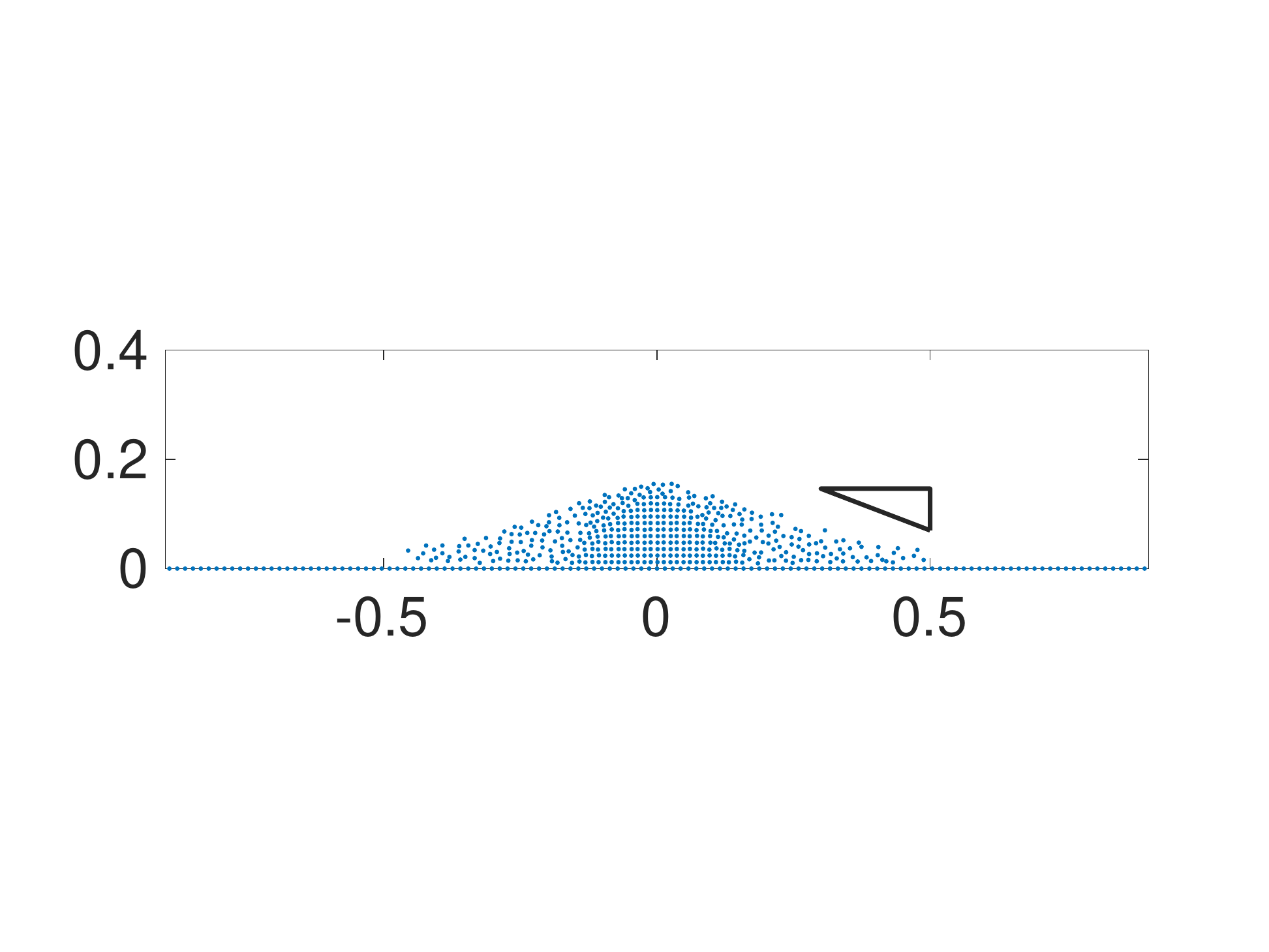}
	\caption{Time development of the solution of plasticity model for $a=1$ at times  $ t = 0.2, 0.3,  0.5$ (from top row to bottom row). First column:   $N=2000$, second column:  $N = 1000$, third column:  $N = 500$.}. 	
	\label{compare_resolution}
	\centering
\end{figure}

\begin{figure}
	\centering
		\includegraphics[keepaspectratio=true, angle=0, width=0.32\textwidth]{column_a_1_micro_t_0dot2} 
		\includegraphics[keepaspectratio=true, angle=0, width=0.32\textwidth]{column_a_1_micro_t_0dot3} 	
		\includegraphics[keepaspectratio=true, angle=0, width=0.32\textwidth]{column_a_1_micro_t_0dot5} 	\\
		\vspace{-1cm}
	\includegraphics[keepaspectratio=true, angle=0, width=0.32\textwidth]{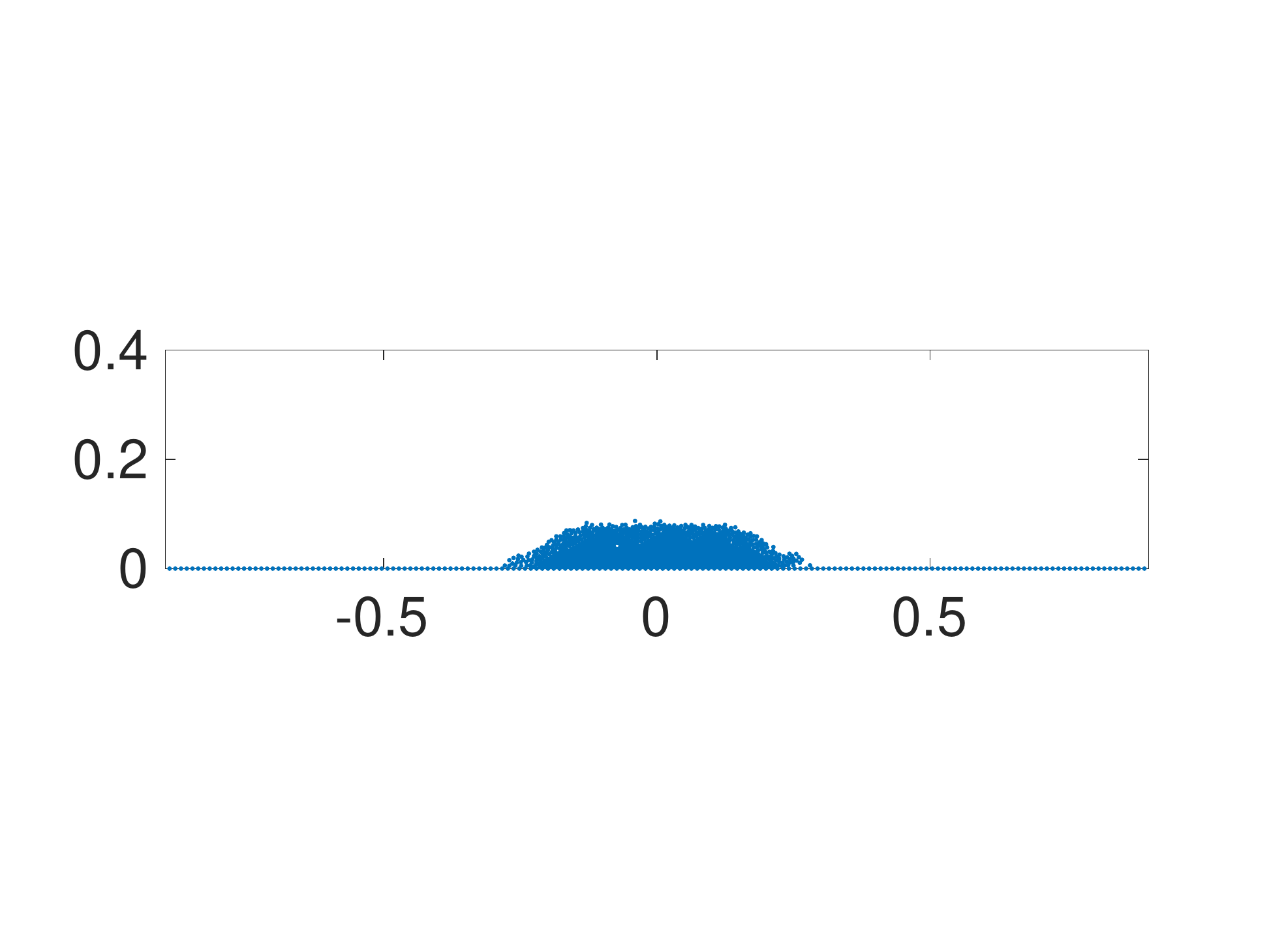} 
	\includegraphics[keepaspectratio=true, angle=0, width=0.32\textwidth]{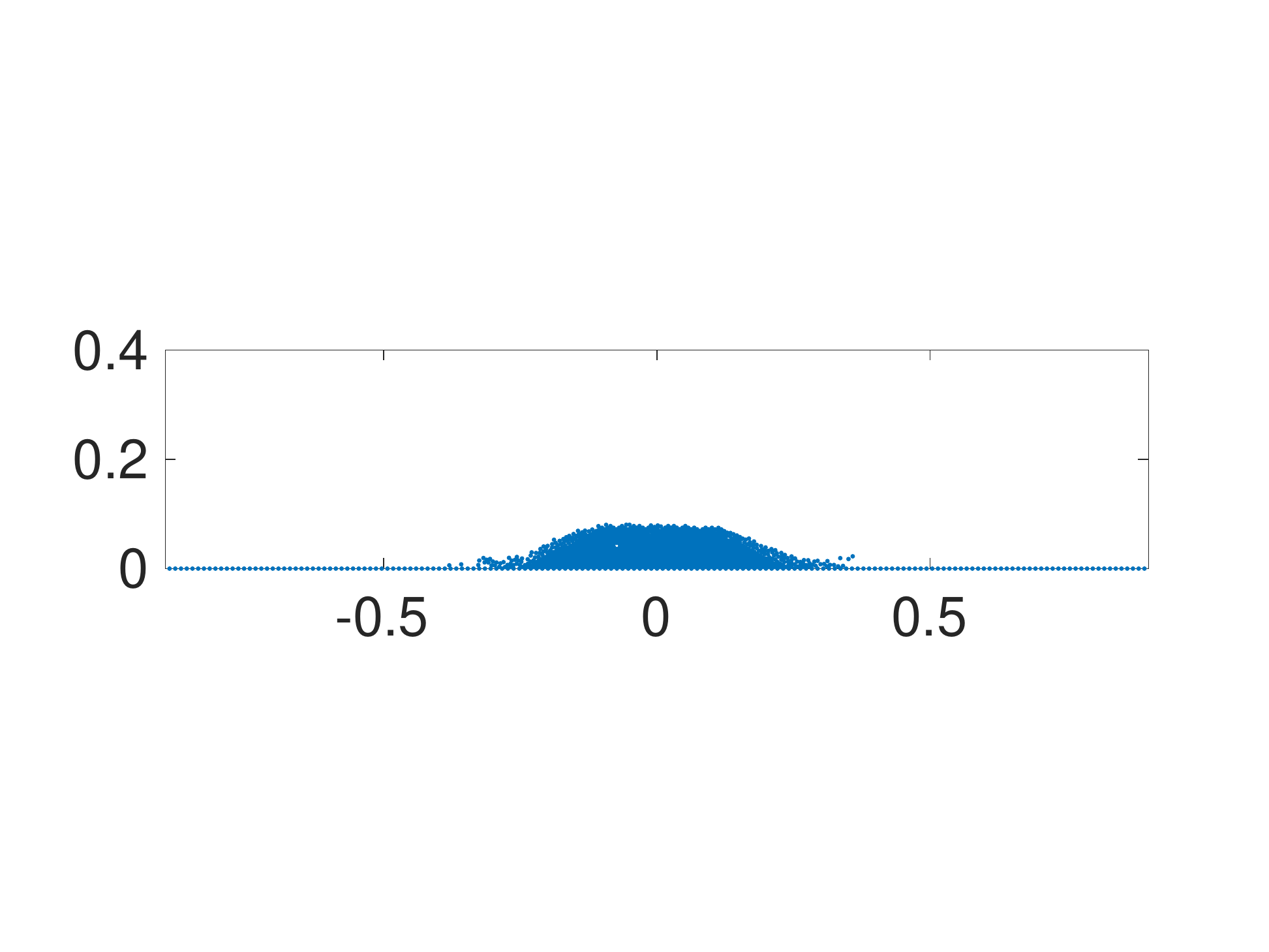} 
	\includegraphics[keepaspectratio=true, angle=0, width=0.32\textwidth]{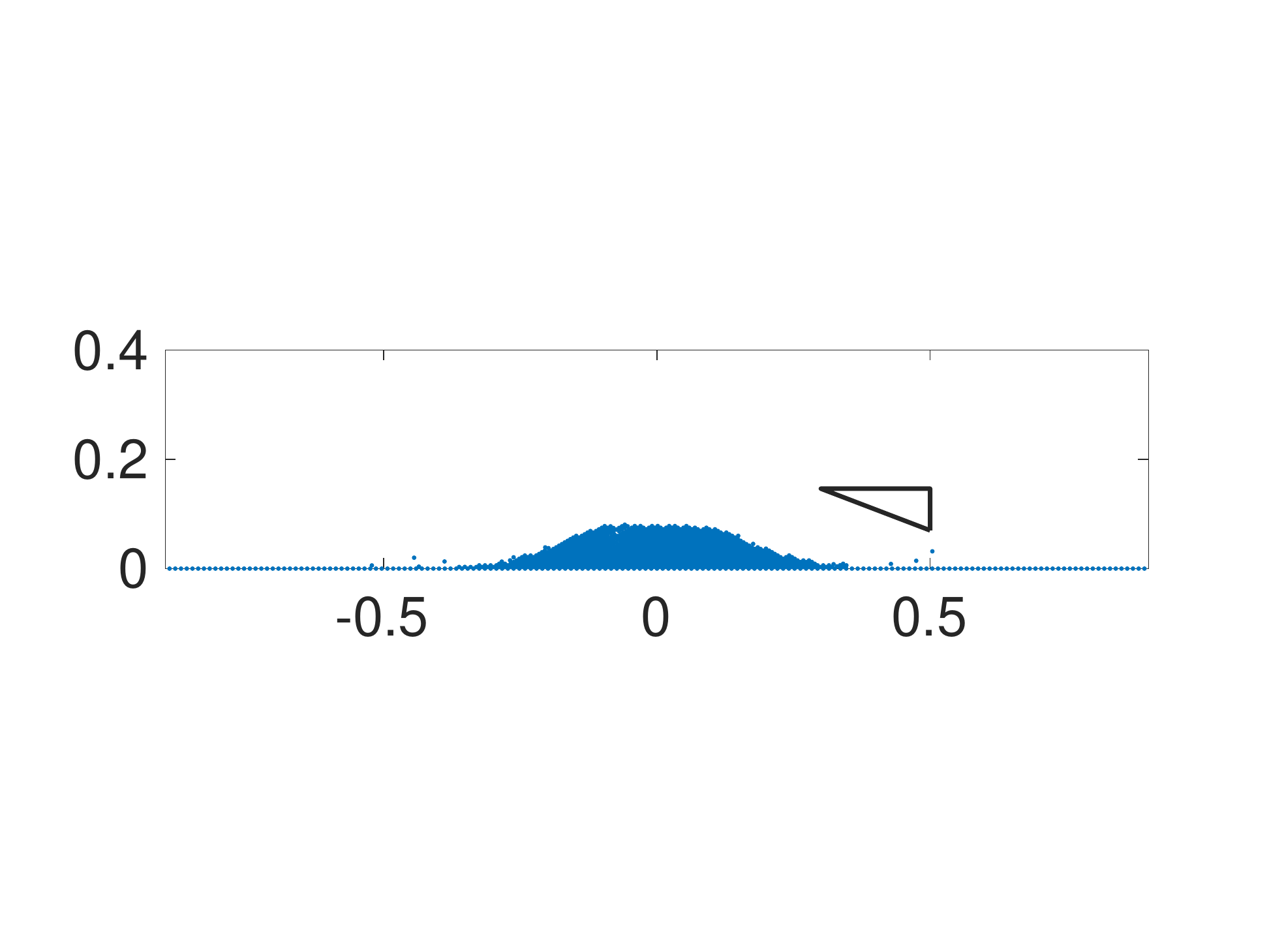} 
	\caption{Time development of the microscopic solution for $a=1$ at times  $ t = 0.2, 0.3,  0.5$ from left to right. First row: $N=2000$, second row:  $N = 1000$.}. 	
	\label{compare_resolution2}
	\centering
\end{figure}

\begin{figure}
	\centering
		\includegraphics[keepaspectratio=true, angle=0, width=0.32\textwidth]{column_a_1_DP_t_0dot2} 
		\includegraphics[keepaspectratio=true, angle=0, width=0.32\textwidth]{column_a_1_DP_t_0dot3} 	
		\includegraphics[keepaspectratio=true, angle=0, width=0.32\textwidth]{column_a_1_DP_t_0dot5} 	\\
		\vspace{-1cm}
	\includegraphics[keepaspectratio=true, angle=0, width=0.32\textwidth]{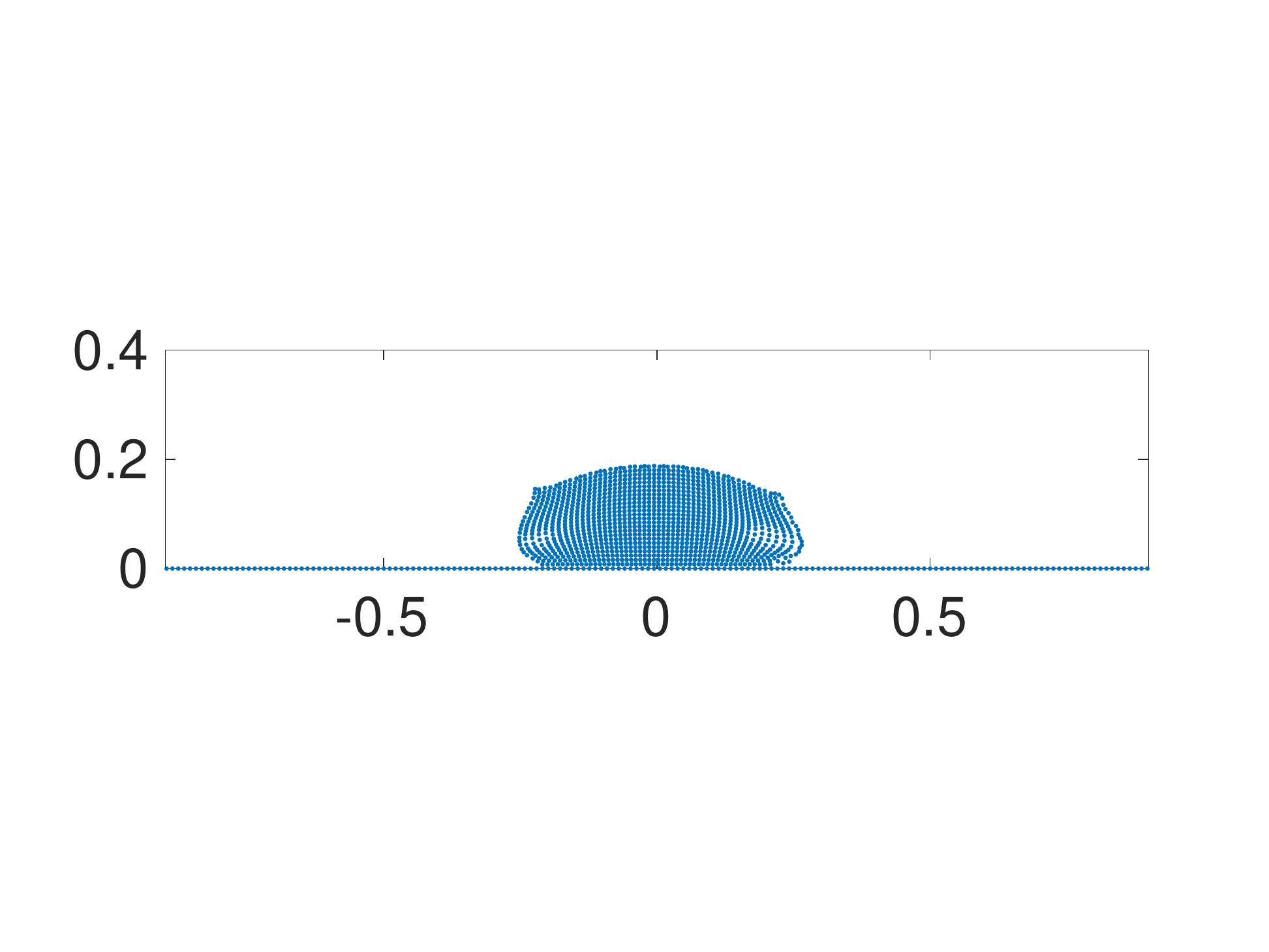} 
	\includegraphics[keepaspectratio=true, angle=0, width=0.32\textwidth]{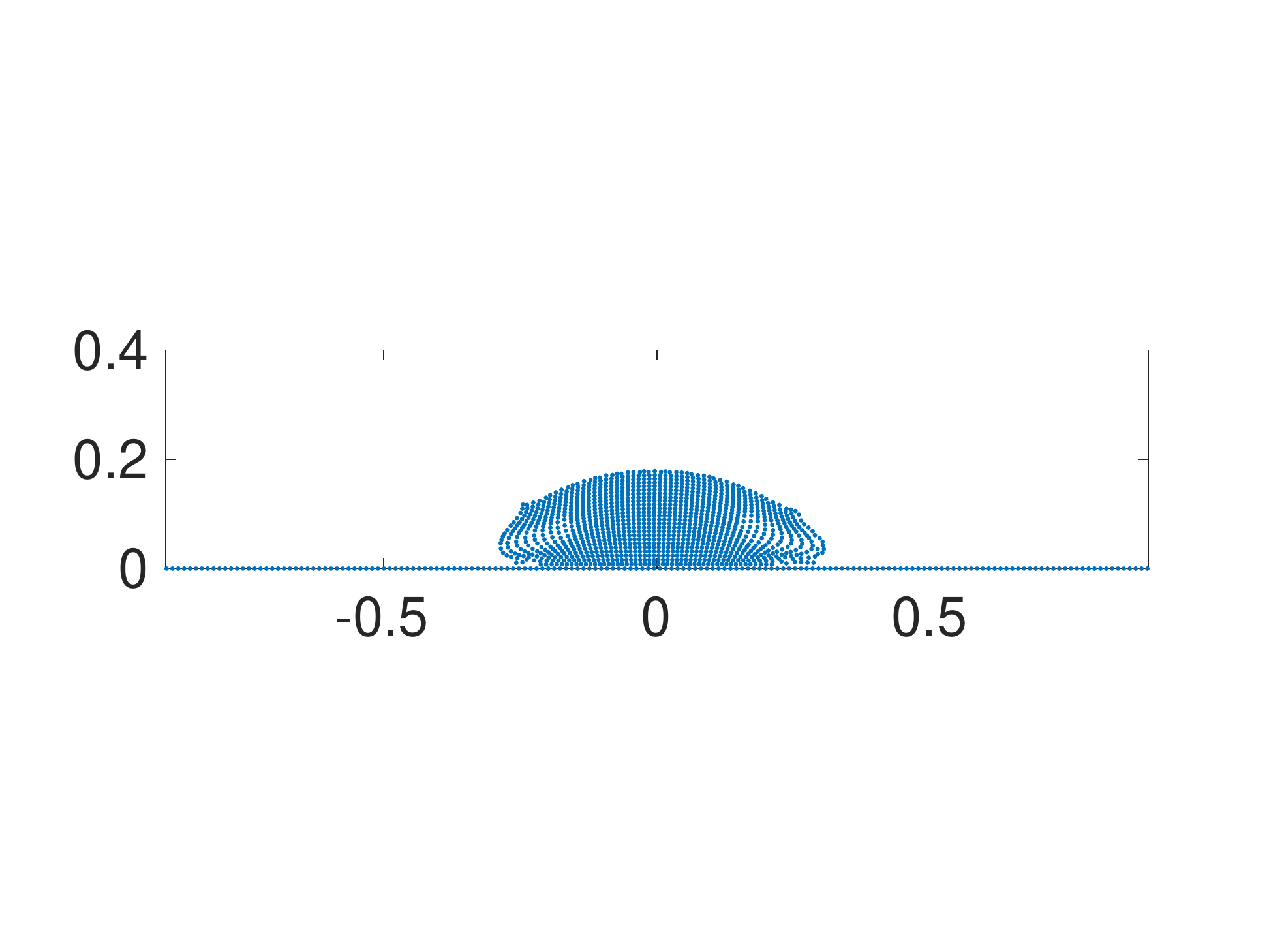} 
	\includegraphics[keepaspectratio=true, angle=0, width=0.32\textwidth]{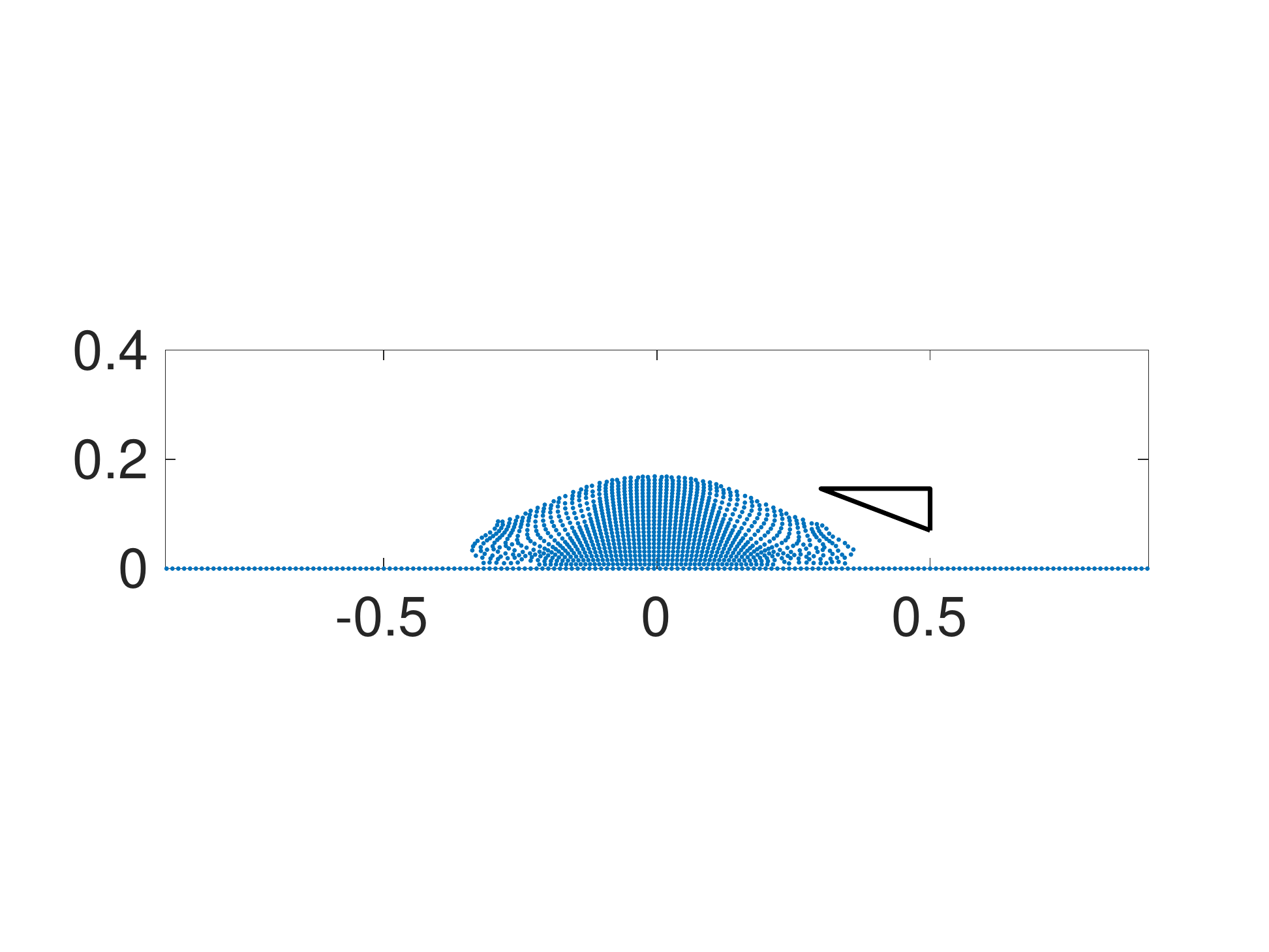} 
 
	\caption{Time development of the solutions of Coloumb constitutive model for $a=1$ at times  $ t = 0.2, 0.3,  0.5$ from left to right. First row: $N=2000$, second row:  $N = 1000$.}. 	
	\label{compare_resolution3}
	\centering
\end{figure}

Additionally, we have investigated the dependence of the results on the number of  microscopic and macroscopic grid particles used.
 In Figure \ref{compare_resolution} we have plotted the results of the plasticity model  with a number of particles $N = 2000, 1000, 500$ leading all to comparable solutions. 
 Similarly in Figure  \ref{compare_resolution2} the results of the microscopic equations  are obtained for $N=2000$ 
 and $N=1000$ particles. 
 In Figure  \ref{compare_resolution3} the results of the simplified macroscopic equations  are obtained for $N=2000$ 
 and $N=1000$ particles. Plotting   times are $ t = 0.2, 0.3,  0.5$ in all cases.  
 
 Note that the constant time step $\Delta t = 10^{-5} s$ is used for the microscopic model and the macroscopic plasticity model
 with $N=2000$.  For $N=500$ a smaller time step still gives stable solutions. We use $\Delta t = 2 \cdot 10^{-5} s$.
 In case of the simplified Coulomb constitutive  model a larger times step can be used due to the implicit time discretization. 
 
 From Figures \ref{compare_resolution} and \ref{compare_resolution2} one observes that the macroscopic 
 plasticity model gives accurate results also for a number of  particles as small as  $N=500$. 
 The results of the microscopic model strongly diverge for 
 $N=1000$ from the ones obtained from $N=2000$ particles.
In case of the simplified macroscopic models one obtains  a behaviour results, which is rather similar to the behaviour of a  viscous fluid.

The simulations is performed in MacOS Monterey Version $12.5$ with processor $3.4$ GHz Quad-Core Intel Core i5, memory $16$ GB $2400$ MHz DDR4, where a single processor is used.
The CPU time for the macroscopic plasticity simulation with $N = 2000$ is about $36$  and for $N = 500$ about $4$ minutes. 
The CPU time for  the microscopic  simulation with $N = 2000$ is around $7$ minutes.  
Thus,  using the macroscopic model with $N=500$ yields a gain of approximately a factor of $2 $ compared  to
the computation times of the microscopic model with $N=2000$.

 

We note that in case of  the simplified mode, due to the possibility of using a smaller time step,
the required computation time for $N=2000$ is an order of magnitude smaller than for the microscopic model.

%
\subsection{Free falling  disc in a $2D$ sand box}

\begin{figure}
	\centering
         \includegraphics[keepaspectratio=true, angle=0, width=0.32\textwidth]{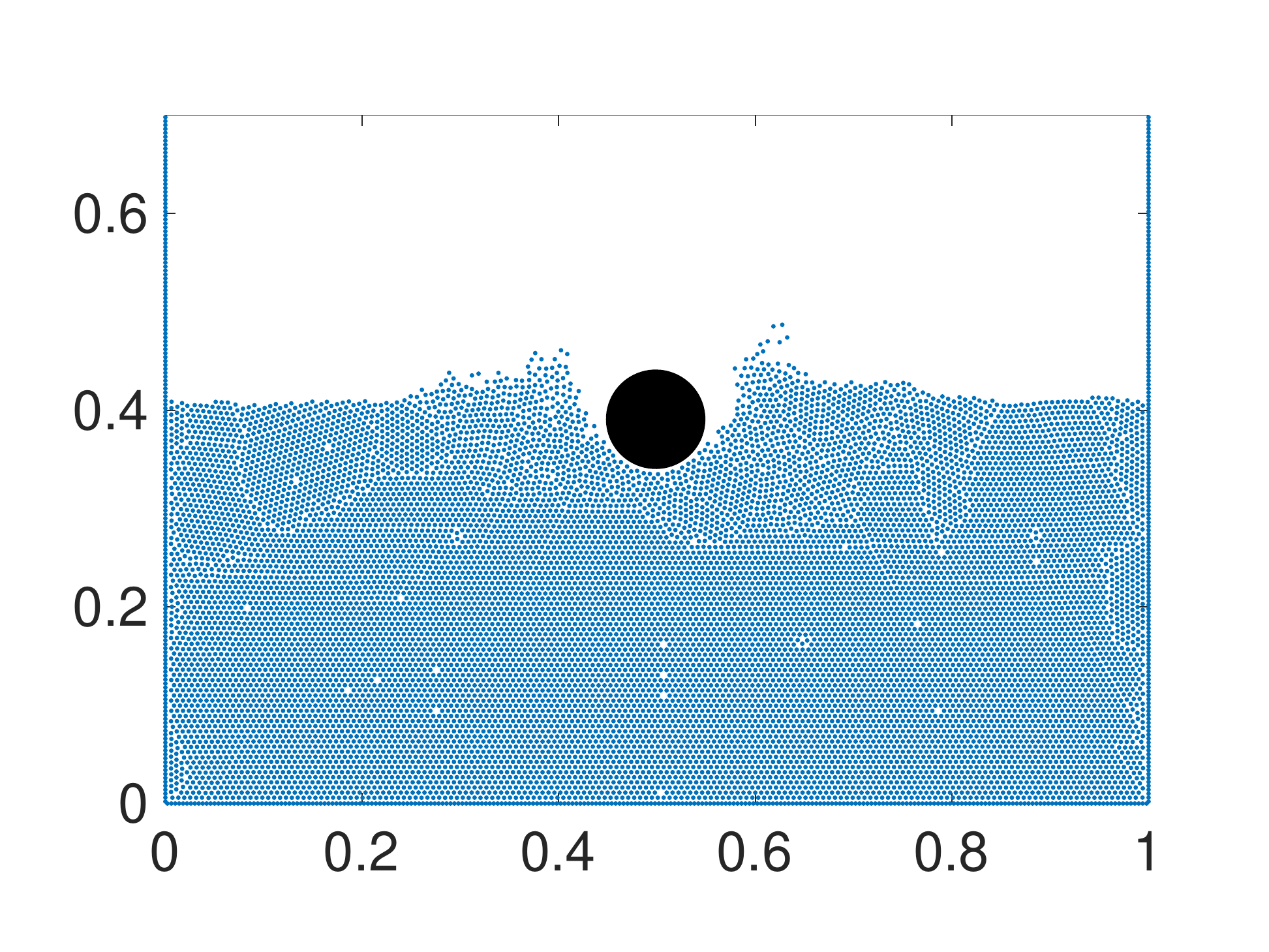} 
          \includegraphics[keepaspectratio=true, angle=0, width=0.32\textwidth]{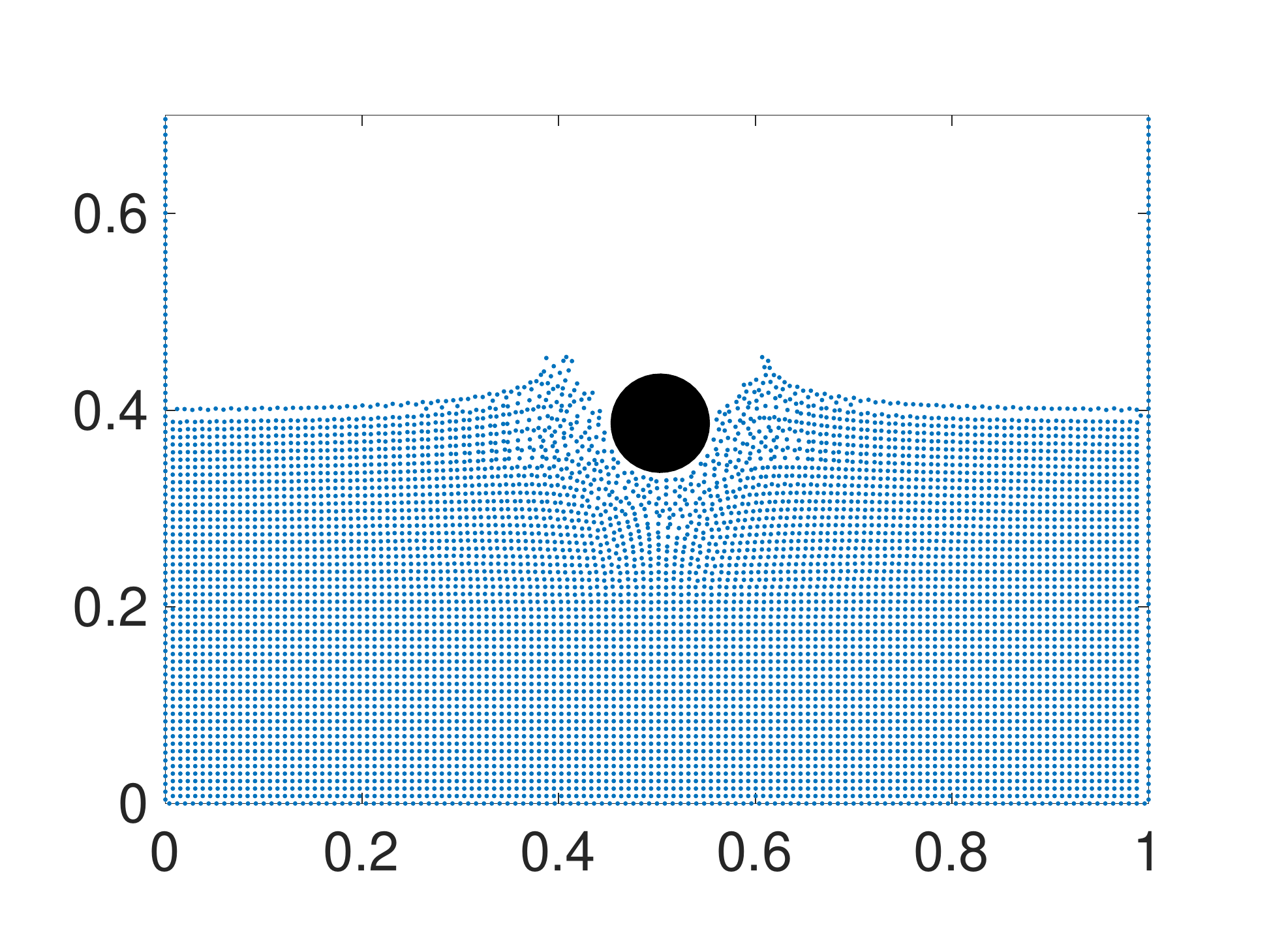} 
           \includegraphics[keepaspectratio=true, angle=0, width=0.32\textwidth]{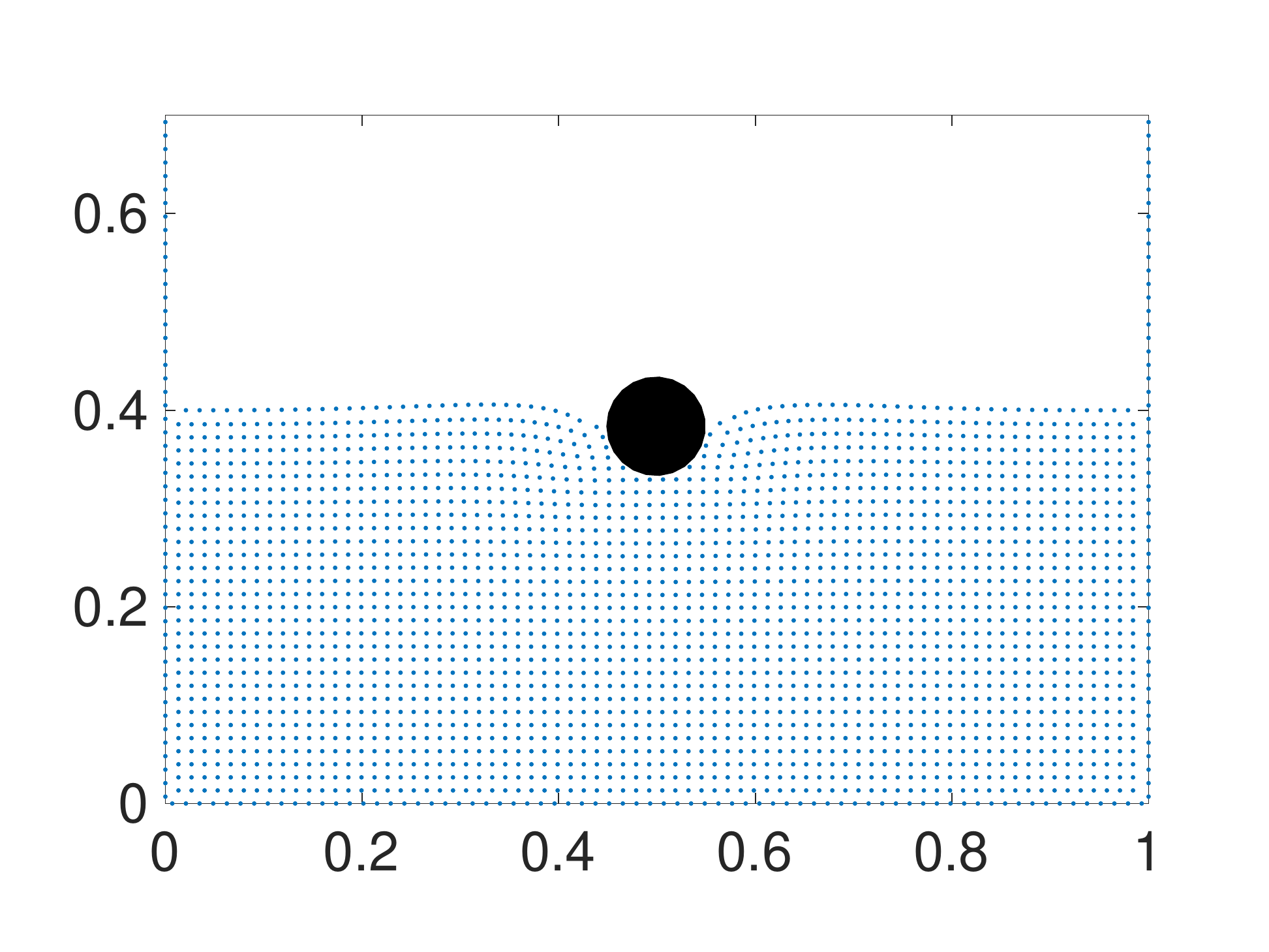} 
          \includegraphics[keepaspectratio=true, angle=0, width=0.32\textwidth]{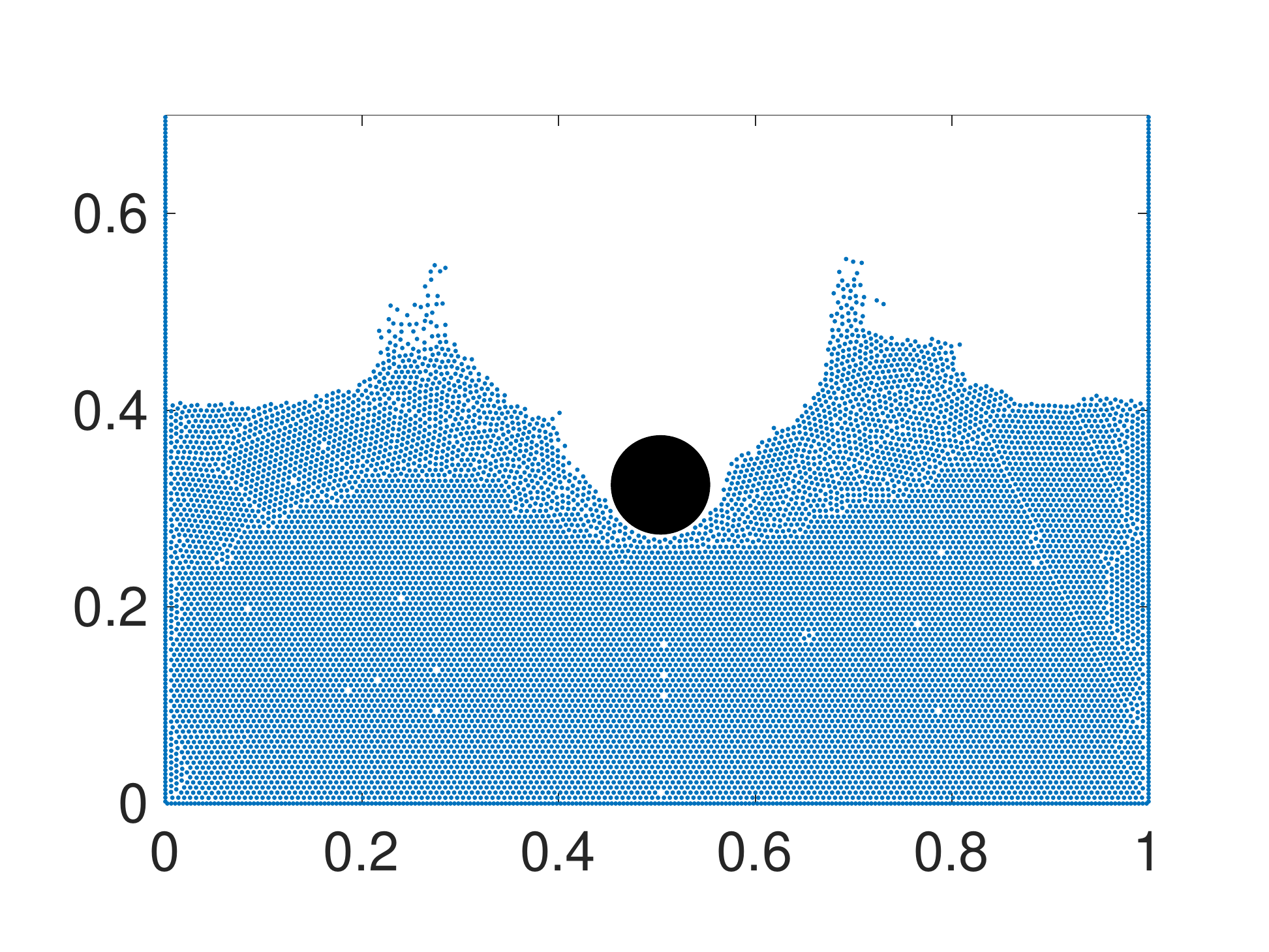} 
          \includegraphics[keepaspectratio=true, angle=0, width=0.32\textwidth]{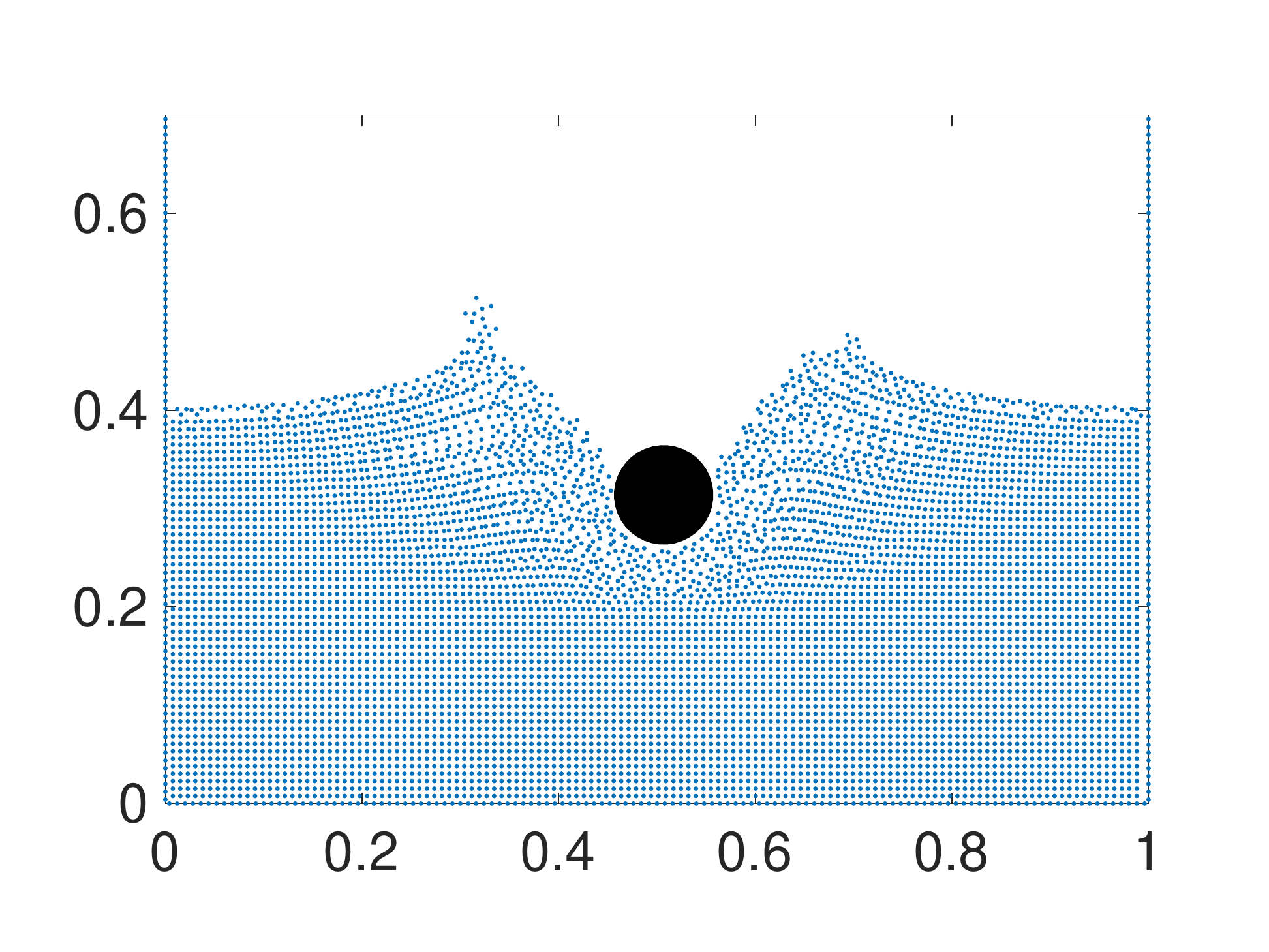} 
           \includegraphics[keepaspectratio=true, angle=0, width=0.32\textwidth]{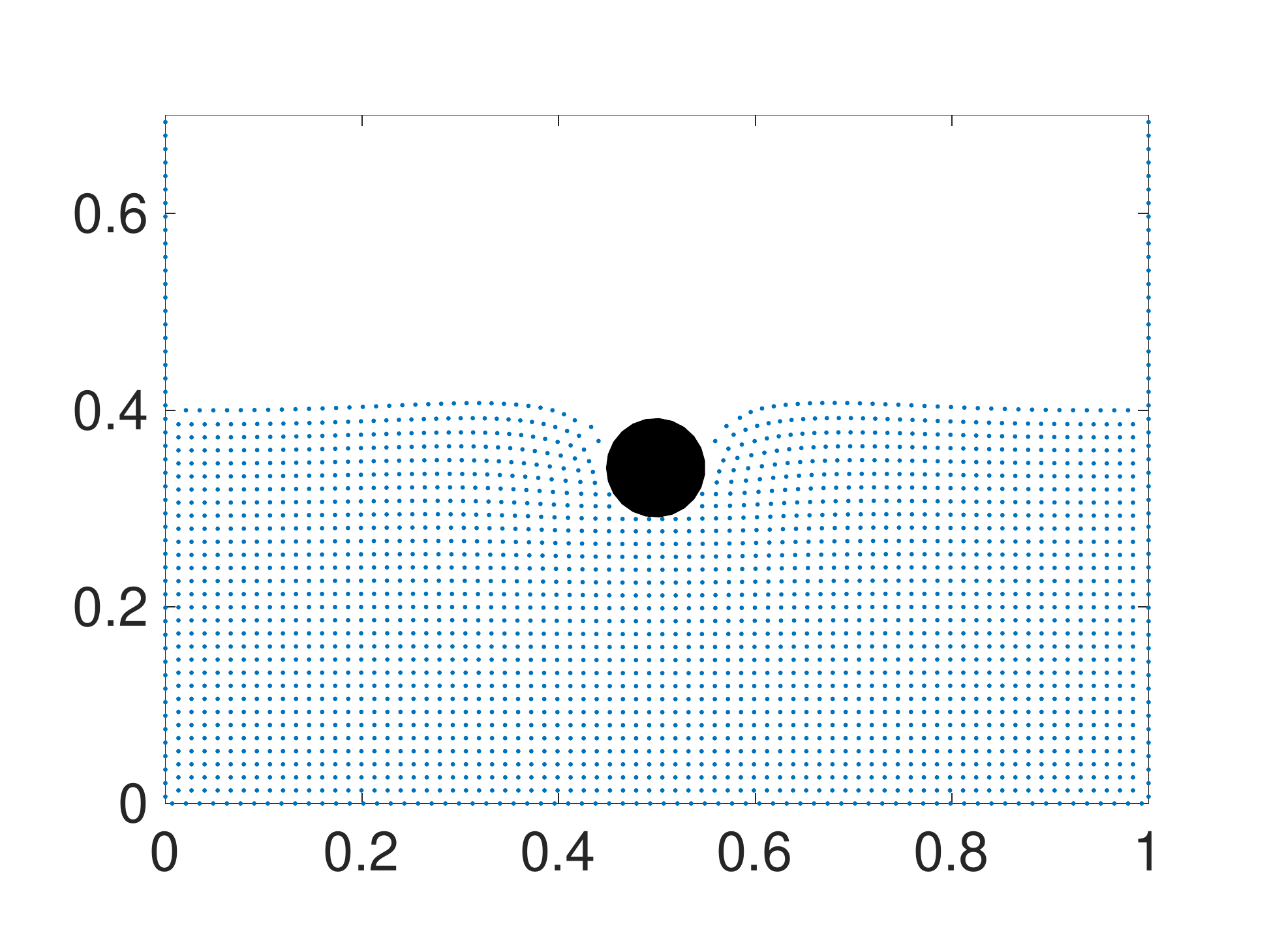} 
	  \includegraphics[keepaspectratio=true, angle=0, width=0.32\textwidth]{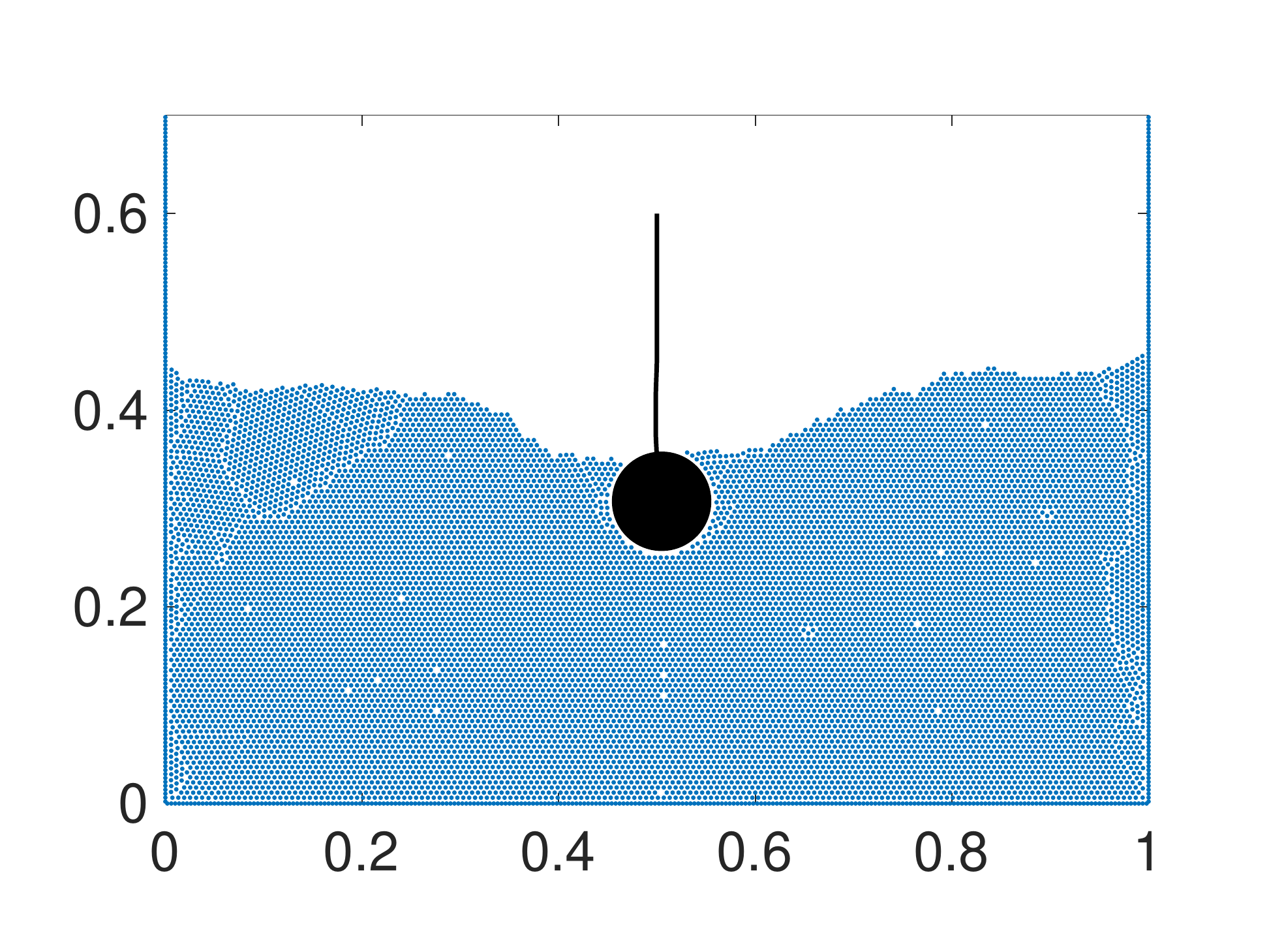} 
          \includegraphics[keepaspectratio=true, angle=0, width=0.32\textwidth]{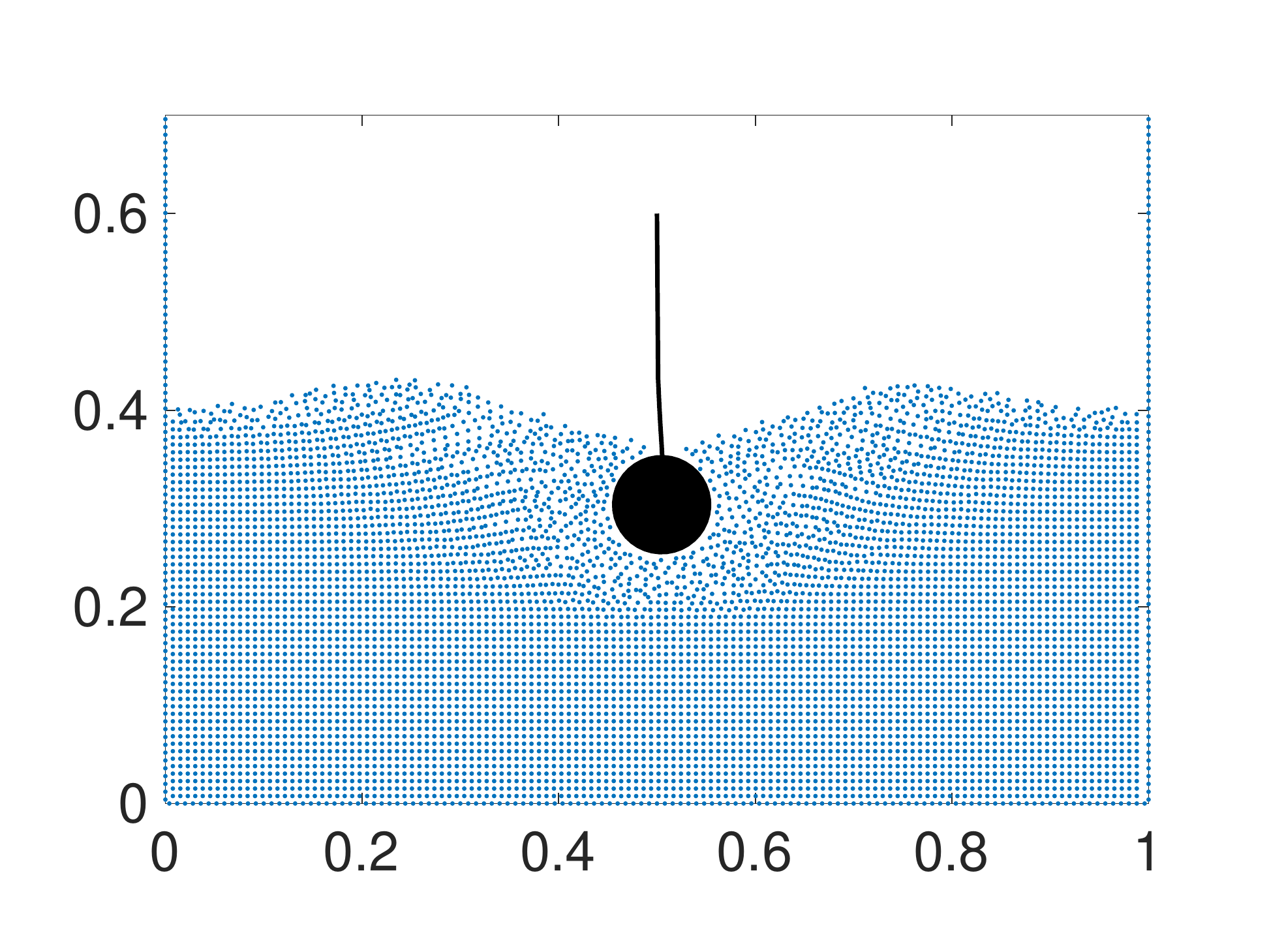} 
           \includegraphics[keepaspectratio=true, angle=0, width=0.32\textwidth]{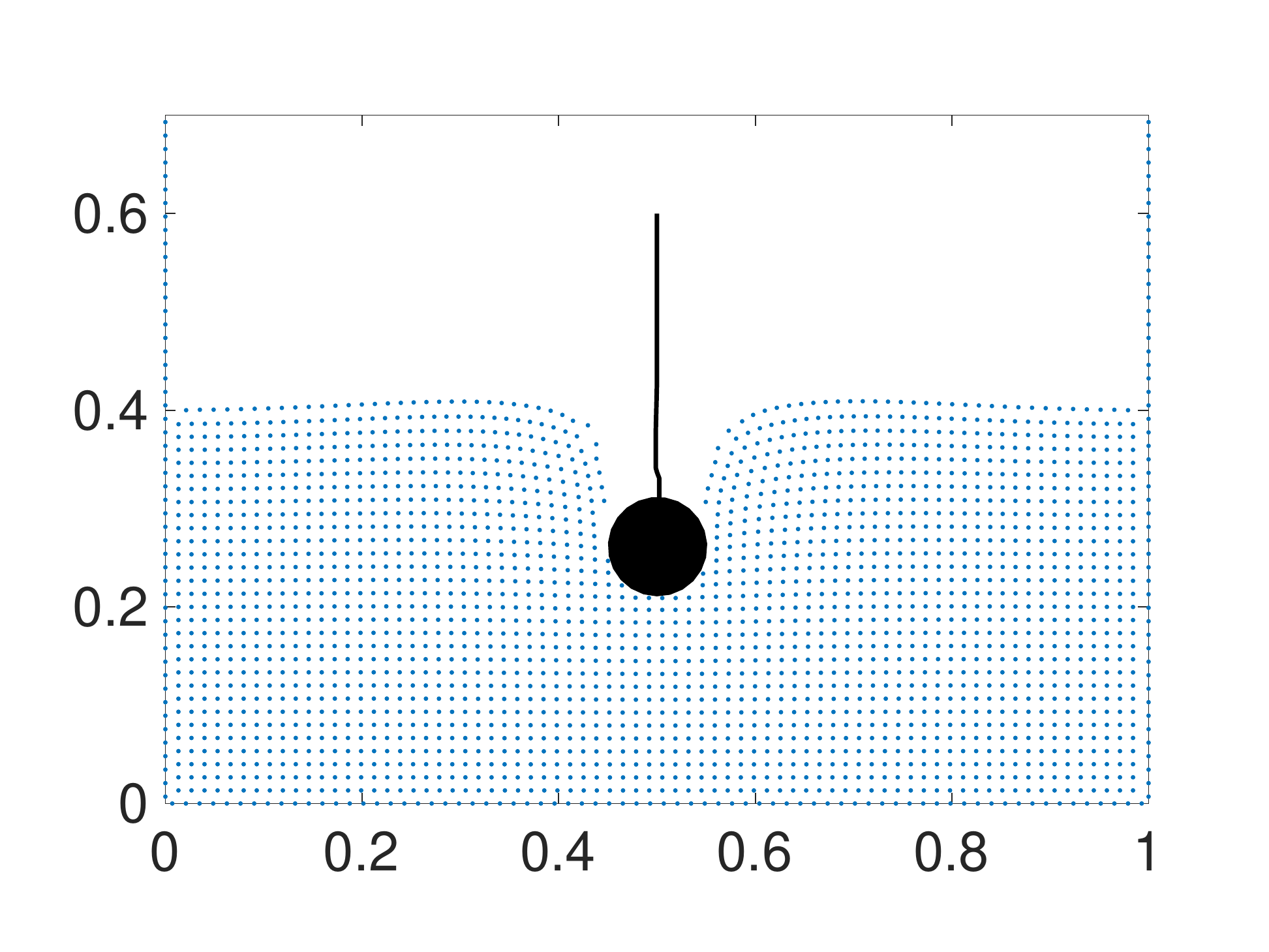} 
          \caption{Comparison of the stationary solutions of the microscopic model (left) plasticity model (middle) and simplified model (right) solutions. The first row is at time $t = 0.08s$, second row at $t = 0.16s$ and the third row is at $t = 1s$ seconds.} 	
	\label{sand_disc}
	\centering

\end{figure}

\begin{figure}
	\centering
         \includegraphics[keepaspectratio=true, angle=0, width=0.6\textwidth]{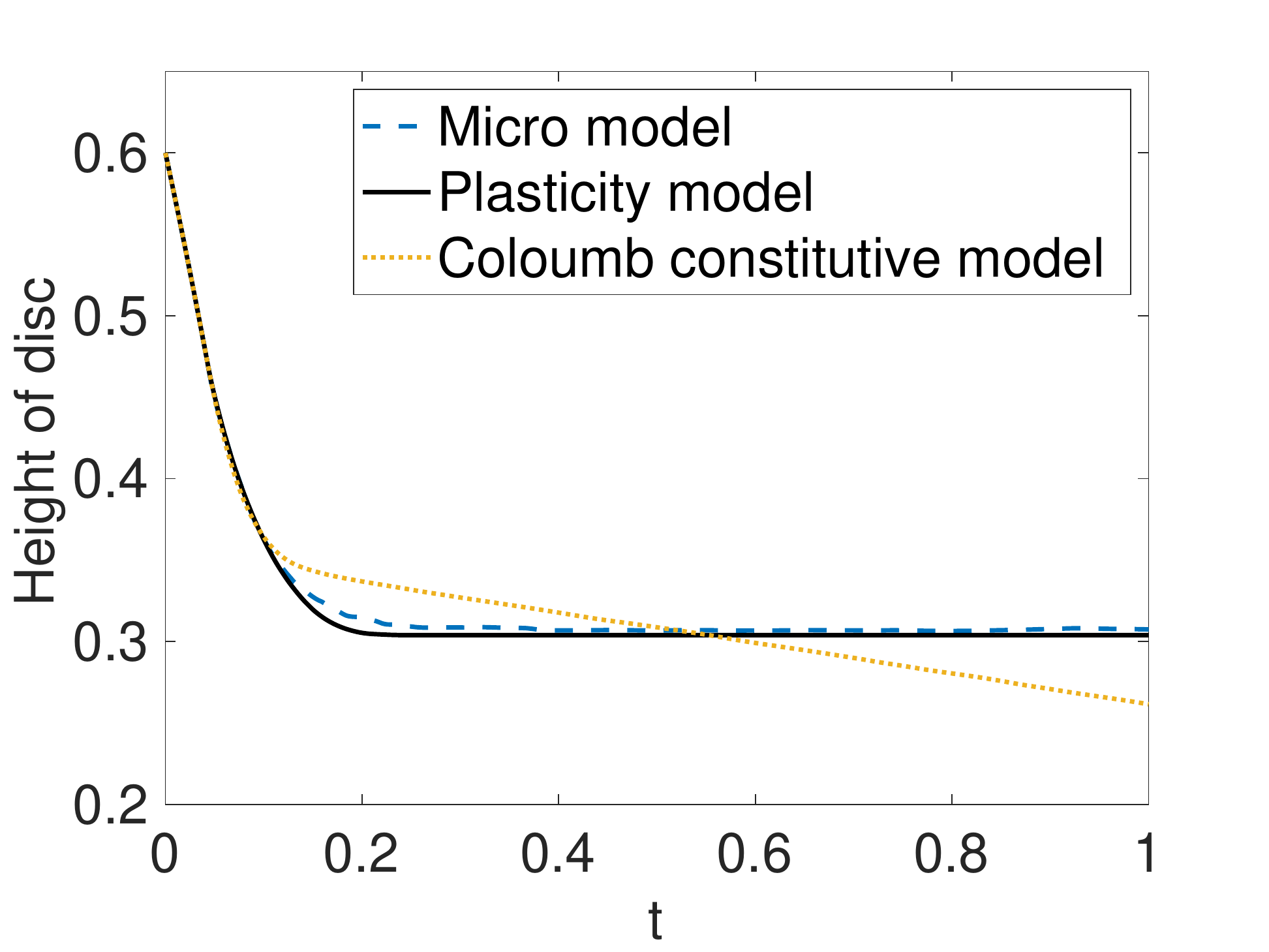}
          \caption{Comparison of the height of the disc with respect to time} 	
	\label{sand_disc_height}
	\centering

\end{figure}

In this subsection we present the simulation of  a falling disc into a  box filled with sand. Consider a box of size $[0,1]m\times[0,1]m$ with the sand initially filled up to the height $y=0.4m$. We consider a case, where a disc of radius $0.05m$  located in the center of the box falls from the height $y = 1m$.  
For the numerical simulation, initally, the center of the disc is located at $(0.5, 0.6)$ and the initial velocity of the disc is equal to $\sqrt{2g0.4}= 2.8014 m/s$, where $g = 9.81 m/{s^2}$. 
The density of the disc is $8050 kg/{m^3}$. The other parameters are the same as in the collapsing column case. 
The time step is chosen as $\Delta t = 5 \cdot 10^{-6}.$
 In figure \ref{sand_disc} we have plotted the positions of sand particles and disc computed  from the  microscopic, the plasticity and  the simplified model at time $0.08s, 0.16s$ and $1s$. The microscopic model is simulated with $10890$ particles and plasticity and simplified models are simulated with initially $2420$ grid points.  Figure  \ref{sand_disc_height} shows the height of the discs versus time. As can be seen in Figure \ref{sand_disc} and Figure  \ref{sand_disc_height}, the solutions obtained from microscopic and plasticity models are very close, however, the solution obtained from the simplified model deviates strongly in this case being 
similiar to a  highly viscous flow solution.

\subsection{Collapsing $3D$ sand column}

\begin{figure}
	\centering
			\includegraphics[keepaspectratio=true, angle=0, width=0.6\textwidth]{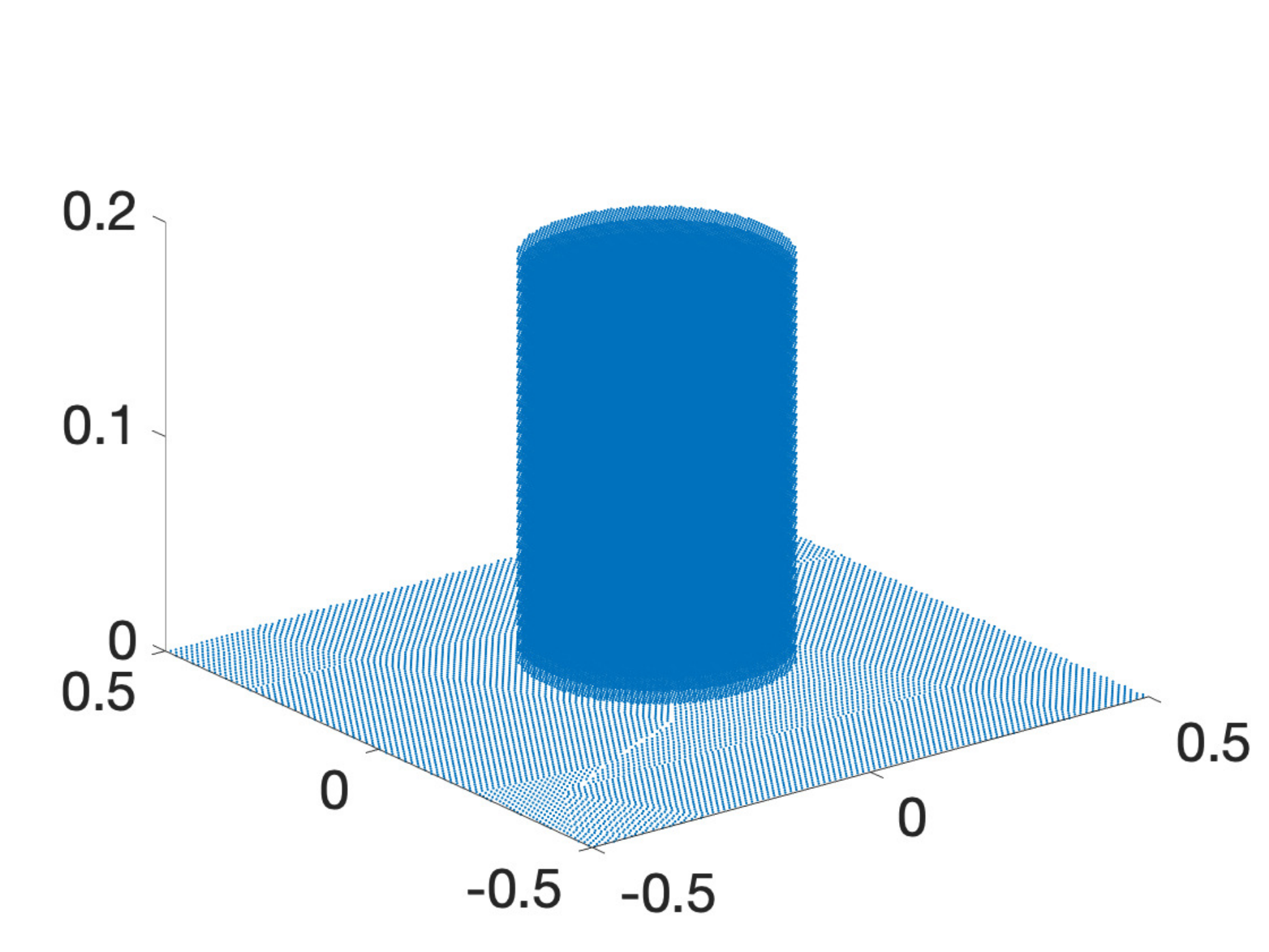} 
			\caption{Initial cylindrical sand column. }
\label{3d_column_t0}
\end{figure}

\begin{figure}
	\centering
	
		 \includegraphics[keepaspectratio=true, angle=0, width=0.32\textwidth]{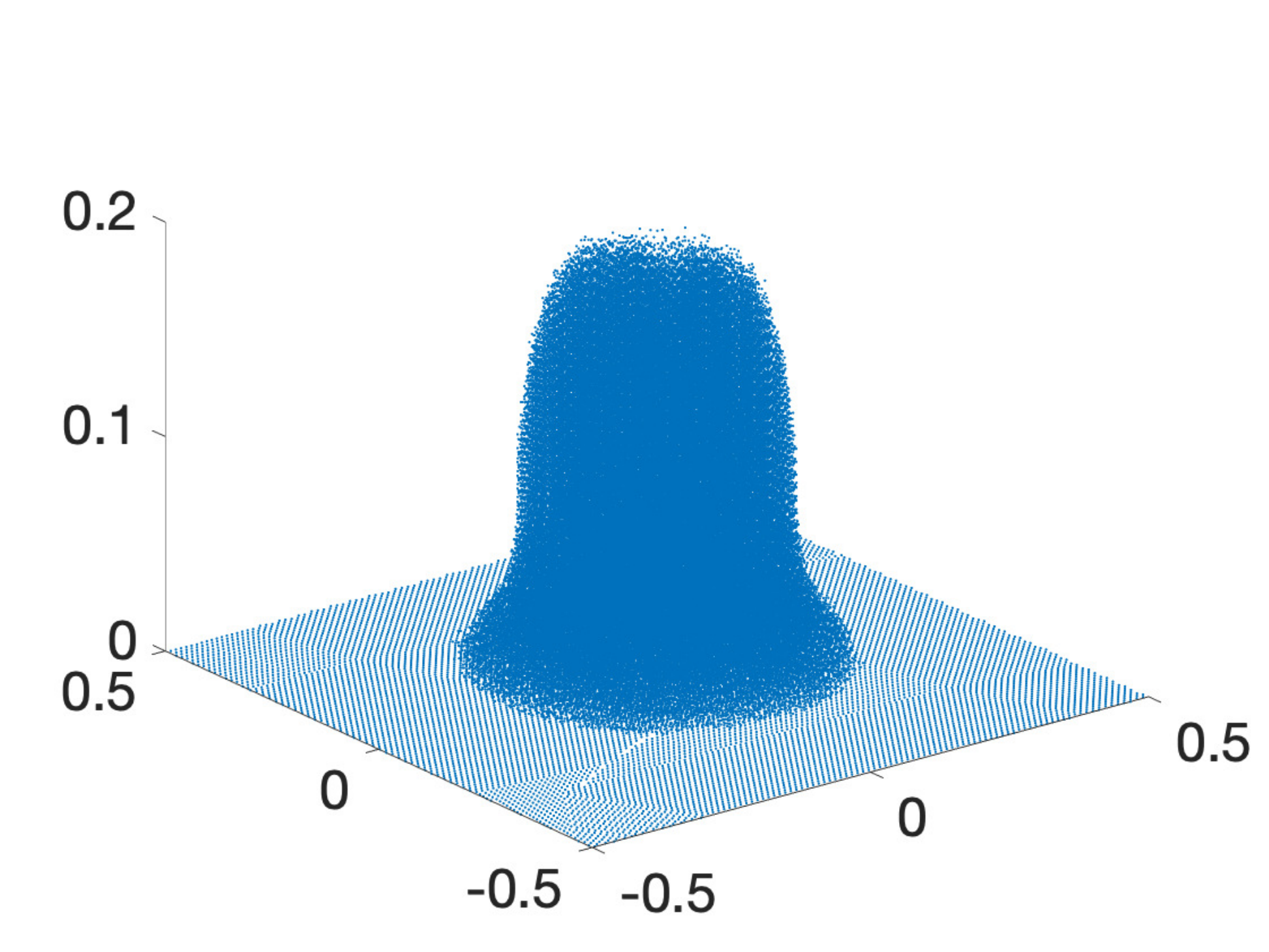} 
		 \includegraphics[keepaspectratio=true, angle=0, width=0.32\textwidth]{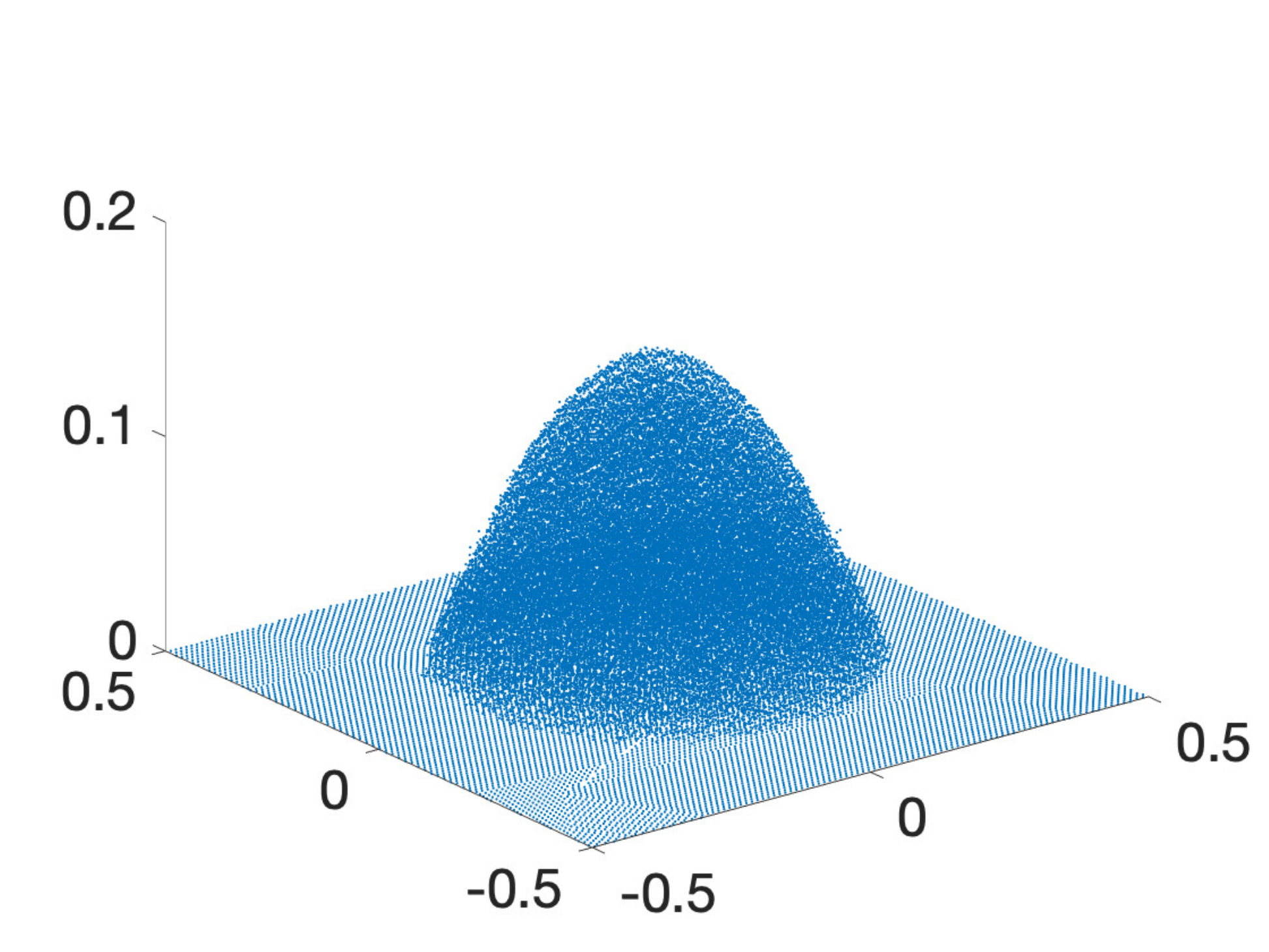} 
		 \includegraphics[keepaspectratio=true, angle=0, width=0.32\textwidth]{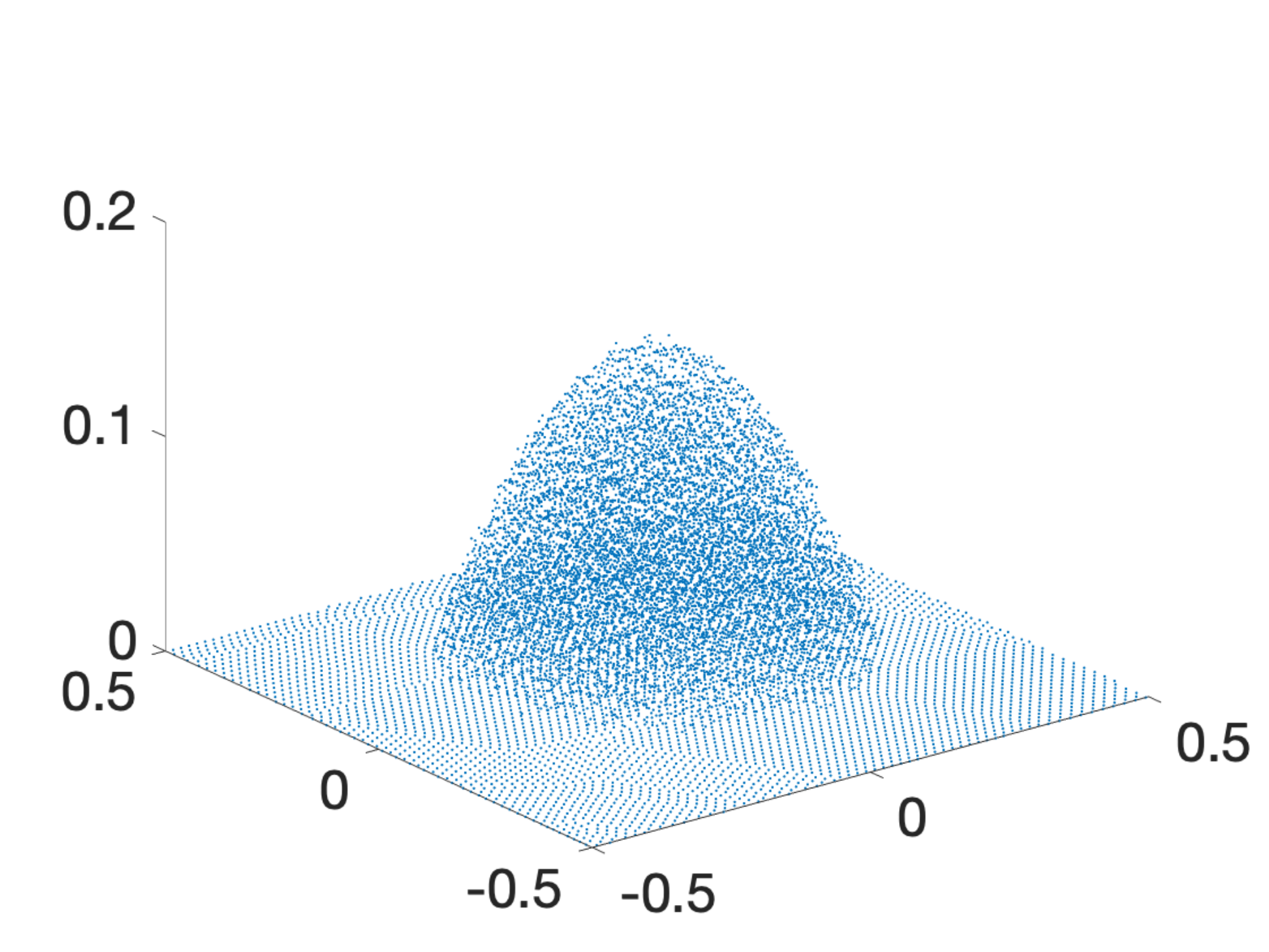} 
		 \includegraphics[keepaspectratio=true, angle=0, width=0.32\textwidth]{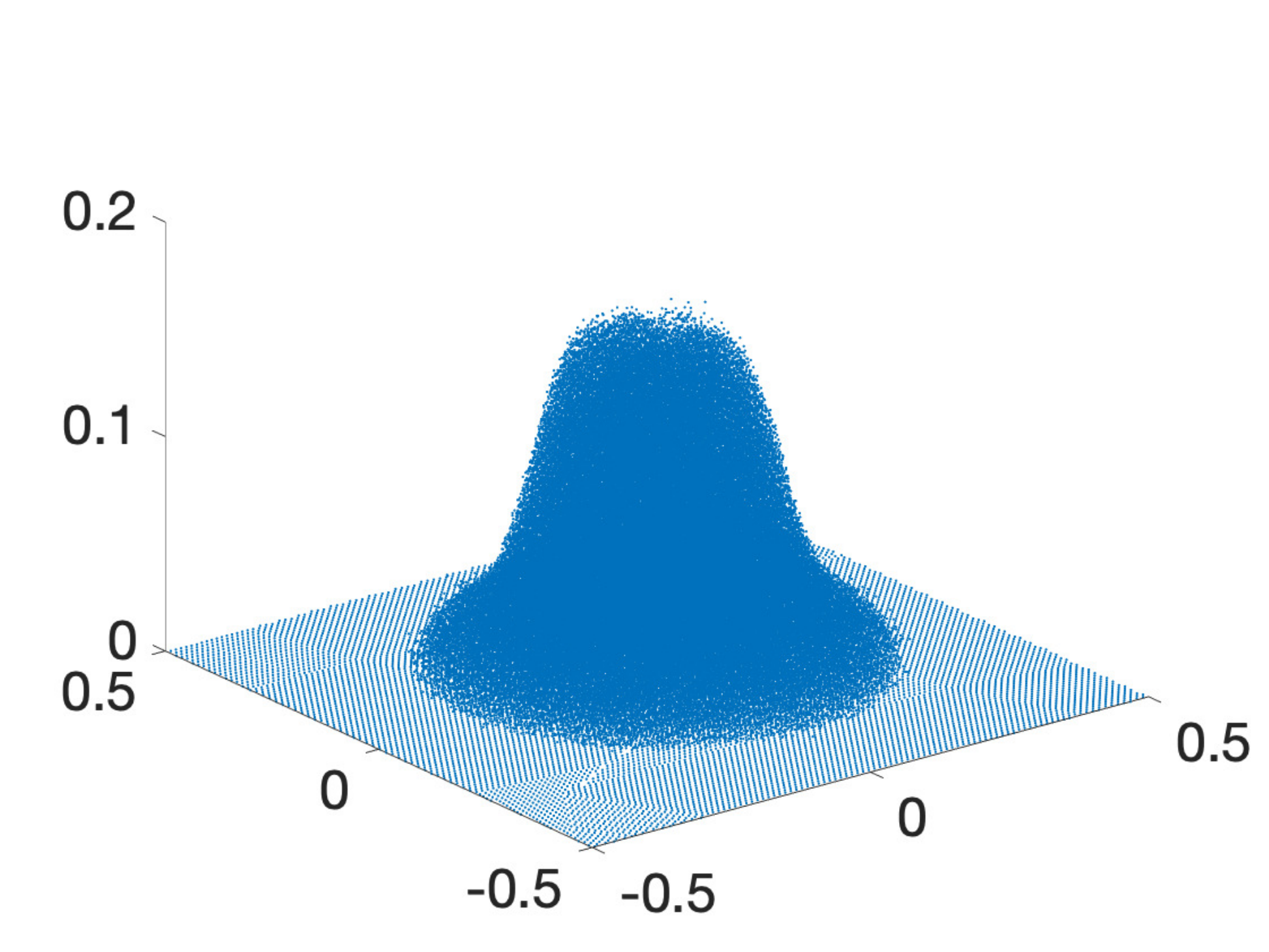} 
		 \includegraphics[keepaspectratio=true, angle=0, width=0.32\textwidth]{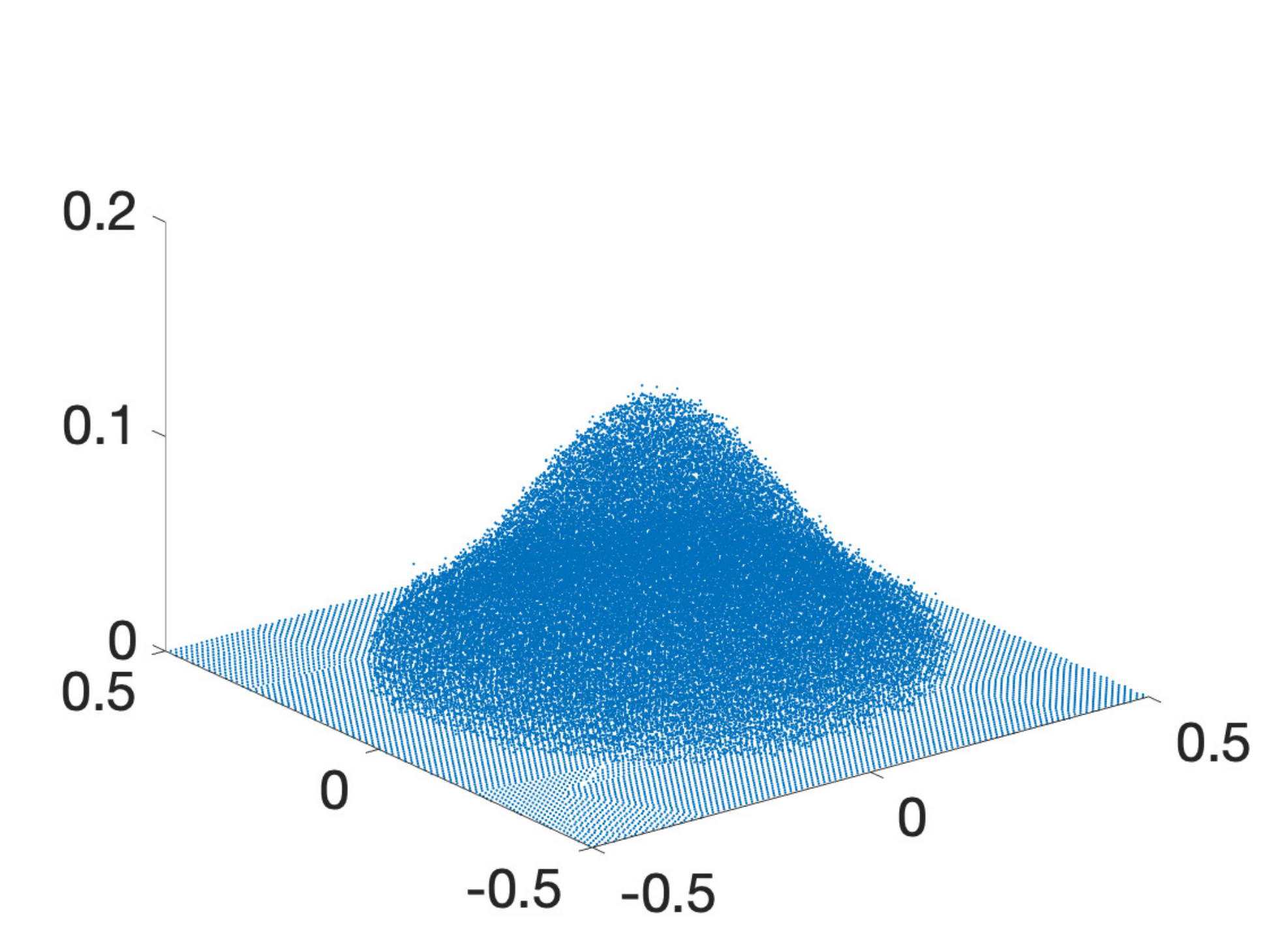} 
		  \includegraphics[keepaspectratio=true, angle=0, width=0.32\textwidth]{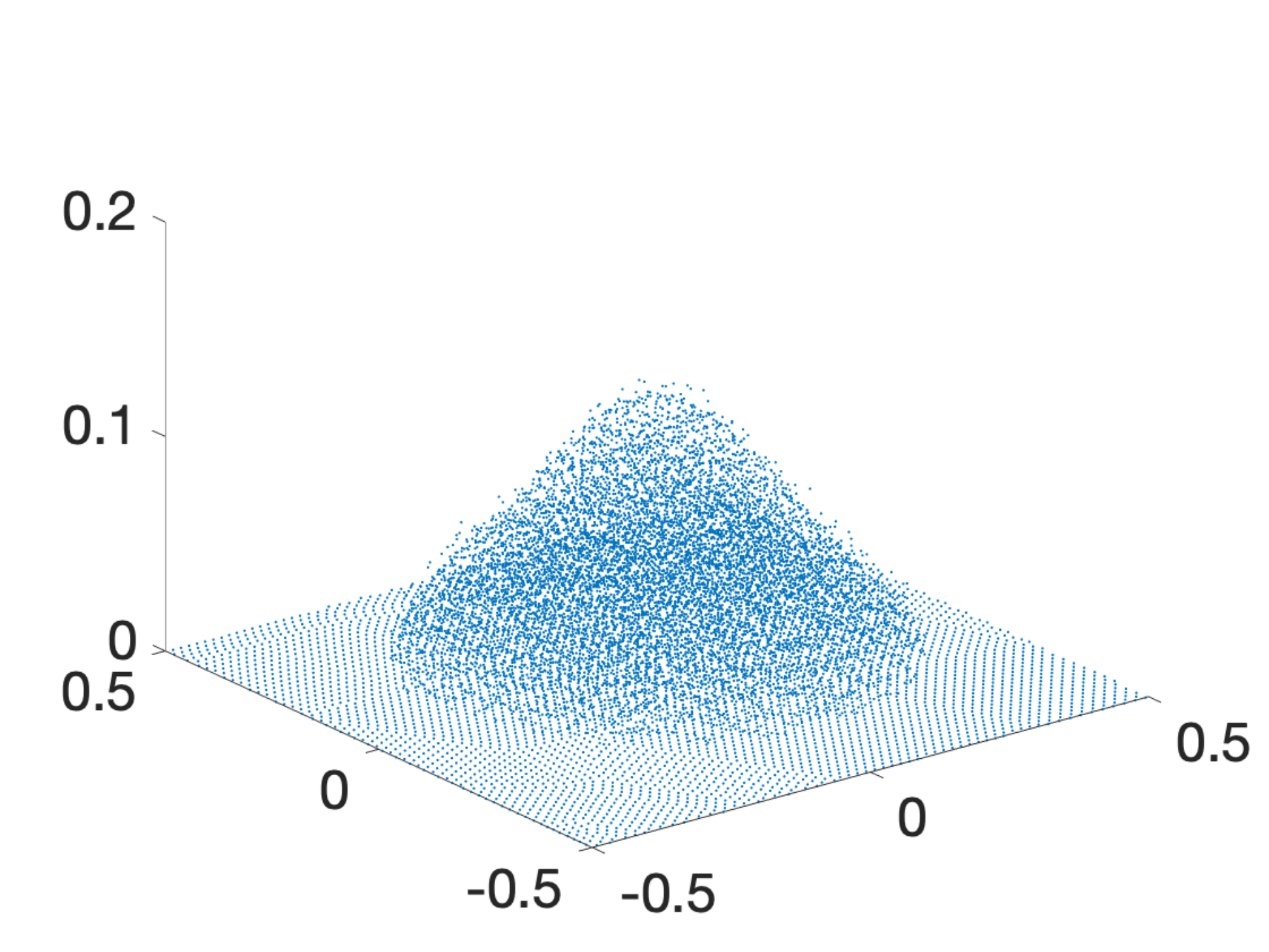} 
	         \includegraphics[keepaspectratio=true, angle=0, width=0.32\textwidth]{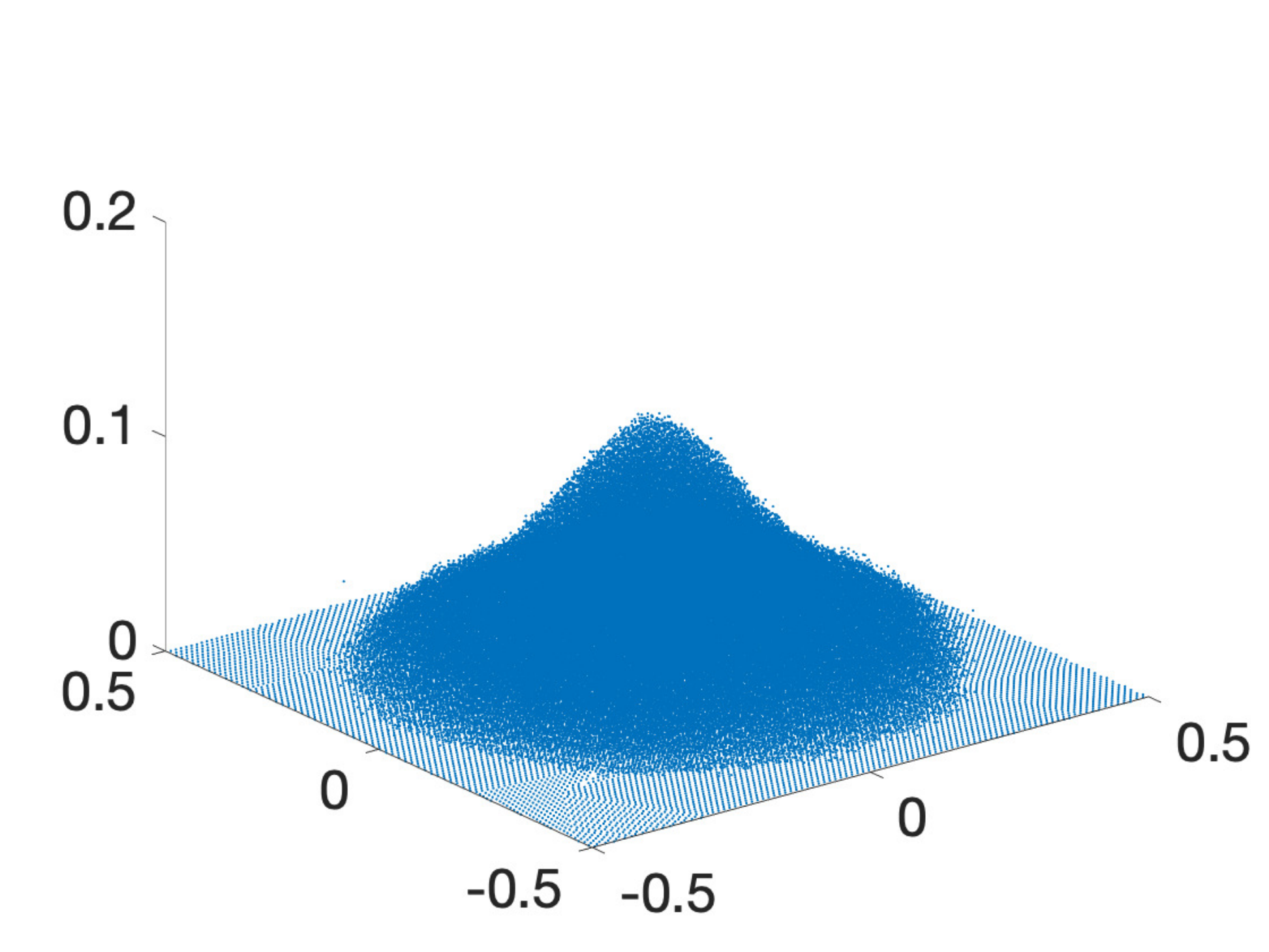} 
		 \includegraphics[keepaspectratio=true, angle=0, width=0.32\textwidth]{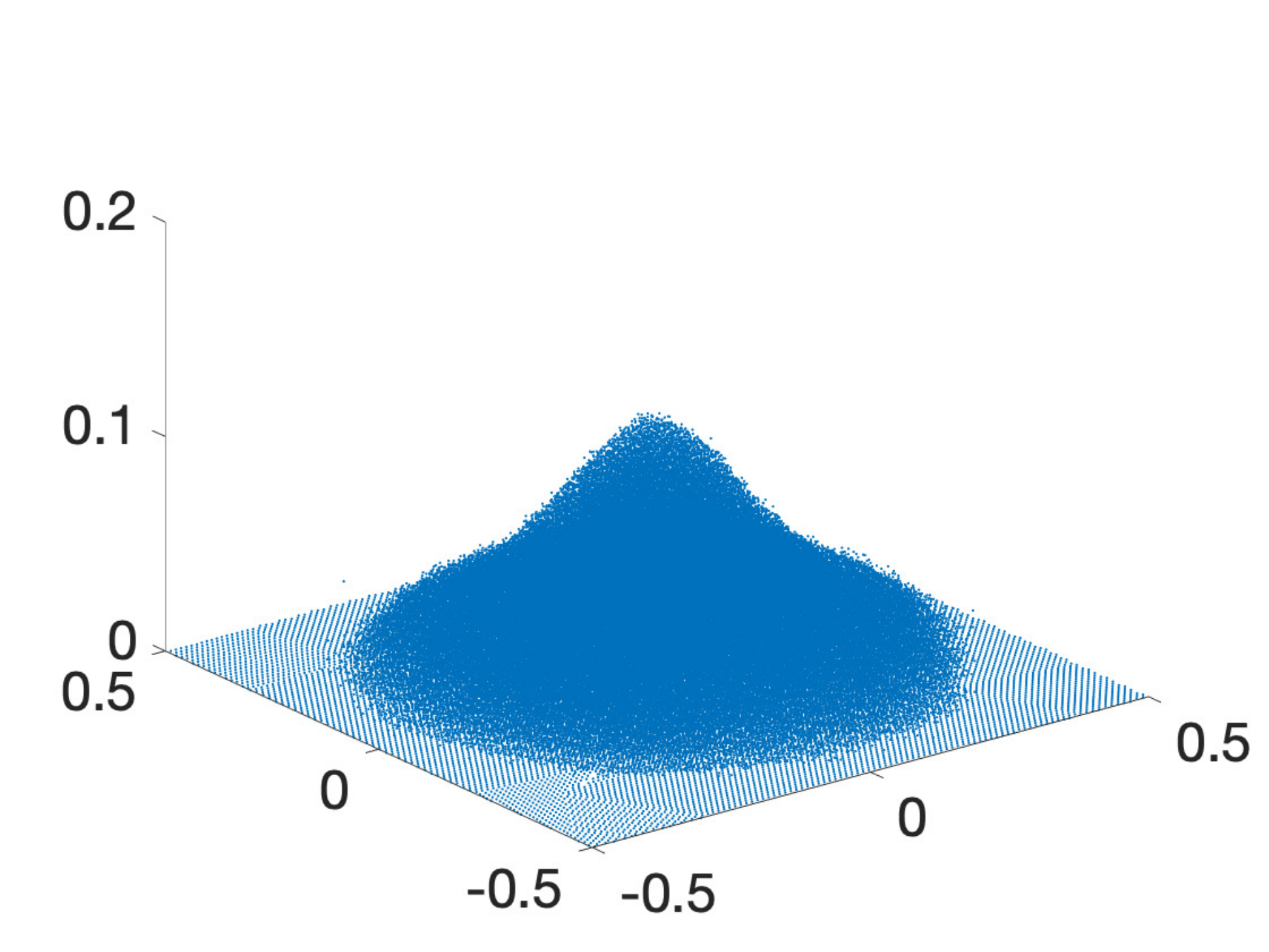} 
		 \includegraphics[keepaspectratio=true, angle=0, width=0.32\textwidth]{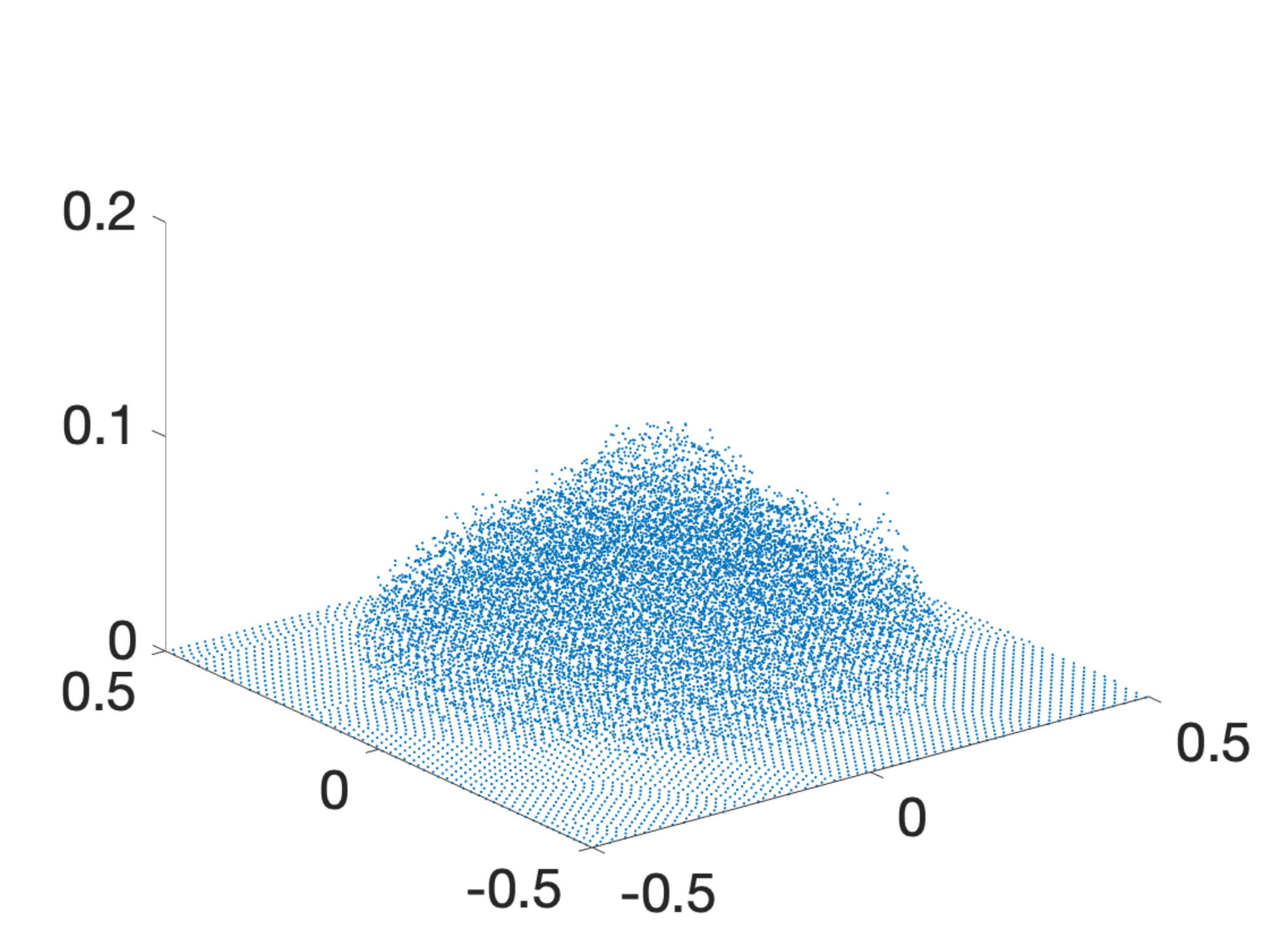} 
	         \caption{Time development of the microscopic model with $153317$ particles (left) and plasticity model with 128105 sand particles (middle) and 24624 sand particles (right) at times $t = 0.2, 0.3 0.5$ (from top to bottom). }
	         \label{3d_column_t0dot5}
\end{figure}

In this subsection we consider a collapsing 3 dimensional sand pile.
We have considered a cylindrical sand pile of radius 0.2 with height 0.2. The parameters are the same as in the $2D$ case. 
Like in the $2D$ case we have plotted in Fig. \ref{3d_column_t0dot5}  the sand position at $t = 0.2, 0.3, 0.5$ obtained from macroscopic plasticity model and from microscopic model.  The CPU times for microscopic and plasticity models are compared up to the final simulation time $t = 0.5$. 
For the microscopic model, the total number of  particles is $153317$ with interior and boundary particles are $121905$ and $31412$, respectively. 
The time step for the microscopic model is $\Delta t = 2 \times 10^{-5}$ is considered. The CPU time for the microscopic model is $8$ hours. 
In the case of the plasticity model, we have used different resolutions with a maximal number of particles $128105$ including $111360$ interior and $16745$ boundary particles. The time step for this resolution is $\Delta t = 10^{-5}$. The CPU time for such a simulation  is $38$ hours. However, using the macroscopic plasticity model  one can obtain the same results with a much coarser resolution. For example, we have considered the total number of particles $24624$ with $19440$ interior and $5184$ boundary particles. Again, the time step is 
$\Delta t = 2\times 10^{-5}$. In this case the CPU time is  only $3$ hours, which is considerably smaller than in the microscopic case. 
The results are plotted in Fig. \ref{3d_column_t0}, where the first column shows the results from the microscopic model,  the second column the ones from the plasticity model with the fine resolution and the third column shows again the  plasticity model, but this  time with the coarse  resolution. The time evolutions of $2D$ and $3D$ columns of microscopic and macroscopic models are similar. In the early stage, for example, at times $t = 0.2$ and $0.3$ the sand column obtained from the microscopic model differs from the macroscopic model, where as in the stationary case, they have very good agreements.

\section{Conclusion}
This work shows the results of an implementation of the elasto-viscoplastic model developed in \cite{DK15} in the framework 
of a GFDM particle method. The case of a collapsing sand column and of a sphere falling into a box filled with sand has been considered. Moreover, a comparison with a microscopic model of nonlinear Hertz type with parameters, which are physically consistent with those of the macroscopic model,  is presented in both cases, showing 
very good coincidence of the two approaches. Additionally a simplified model based on the $\mu(I)$ rheology in the framework of  incompressible balance equations is investigated.

 In particular, for the case of the falling disc advantages of the approach in \cite{DK15} compared to a simplified use of the $\mu(I)$ rheology in the framework of  incompressible balance equations are 
clearly observed.
Concerning the computational times required for the simulations, the results for the  collapsing sand columns  shows, that  the number of particles can be  considerably  reduced in the  case of the macroscopic plasticity model compared to the microscopic model. This  leads to  a reduction in  CPU times between microscopic and macroscopic simualtions by a factor of 3.

\subsection*{Acknowledgments} 
This  work is supported by the BMBF-Program  'Mathematik f\"ur Innovationen', Project HYDAMO
and by the DFG (German research foundation) under Grant No. KL 1105/30-1.

%


\end{document}